\documentclass[11pt]{amsart}
\usepackage{amssymb}
\usepackage{tikz} 
\usetikzlibrary{arrows,backgrounds,fit,positioning,shapes.symbols,chains,cd,decorations.markings}
\usepackage{hyperref}

\numberwithin{equation}{section}
\newtheorem{theorem}{Theorem}[section]
\newtheorem{proposition}[theorem]{Proposition}
\newtheorem{defprop}[theorem]{Definition-Proposition}
\newtheorem{lemma}[theorem]{Lemma}
\newtheorem{corollary}[theorem]{Corollary}
\newtheorem{conj}[theorem]{Conjecture}
\theoremstyle{definition}
\newtheorem{definition}[theorem]{Definition}
\newtheorem{remark}[theorem]{Remark}
\newtheorem*{remark*}{Remark}

\newcommand{\Z}{\mathbb{Z}}
\newcommand{\R}{\mathbb{R}}
\newcommand{\C}{\mathbb{C}}
\newcommand{\CP}{\mathbb{CP}}
\newcommand{\PP}{\mathbb{P}}
\newcommand{\FF}{\mathbb{F}}
\newcommand{\K}{\mathbb{K}}
\renewcommand{\O}{\mathcal{O}}

\newcommand{\dbar}{\bar\partial}
\newcommand{\m}{\mathfrak{m}}
\renewcommand{\b}{\mathfrak{b}}
\newcommand{\CC}{\mathfrak{C}}
\newcommand{\val}{\mathrm{val}}
\newcommand{\Mbar}{\overline{\mathcal{M}}}

\title[Discs of negative Maslov index and extended deformations]
{Holomorphic discs of negative Maslov index and extended deformations in 
mirror symmetry}
\author{Denis Auroux}
\address{Harvard University, Department of Mathematics, 1 Oxford St., Cambridge 
MA 02138, USA}
\email{auroux@math.harvard.edu}

\begin{document}

\begin{abstract}
The SYZ approach to mirror symmetry for log Calabi-Yau
manifolds starts
from a Lagrangian torus fibration on the complement of an anticanonical 
divisor. A mirror space is constructed 
by gluing local charts (moduli spaces of local systems on generic torus
fibers) via wall-crossing transformations which account for corrections to
the analytic structure of moduli spaces of objects of the Fukaya category
induced by bubbling of Maslov index 0 holomorphic discs,
and made into a Landau-Ginzburg model by equipping it with a 
regular function (the superpotential) which enumerates Maslov index 2 
holomorphic discs.

When they occur, holomorphic discs of negative Maslov index deform this
picture by introducing inconsistencies in the wall-crossing transformations,
so that the mirror is no longer an analytic space;
the geometric features of the corrected mirror can be understood in the
language of extended deformations of Landau-Ginzburg models. 
We illustrate this phenomenon (and show
that it actually occurs) by working
through the construction for an explicit example (a log Calabi-Yau 4-fold 
obtained by blowing up a toric variety), and discuss a family Floer approach
to the geometry of the corrected mirror in this setting. Along the way, we
introduce a Morse-theoretic model for family Floer theory which may be of
independent interest.
\end{abstract}

\maketitle


\section{Introduction}

\subsection{SYZ mirror symmetry relative to a nef anticanonical divisor}
The Strominger-Yau-Zaslow (SYZ) approach to mirror symmetry
gives a geometric construction of mirror spaces from Lagrangian torus
fibrations on Calabi-Yau manifolds: roughly speaking, a mirror Calabi-Yau
is obtained as a dual torus fibration, modified by ``instanton corrections''
in the presence of singular fibers \cite{SYZ}. 
A more modern interpretation of the SYZ conjecture describes the mirror as a moduli space of 
objects of the Fukaya category of $X$ supported on the torus fibers; this
viewpoint leads naturally to Fukaya's family Floer program
\cite{Ffamily, AbICM, AbFFF, Tu, Yuan}, which produces a rigid analytic
mirror space out of a fibration by unobstructed Lagrangian tori of vanishing
Maslov class (as well as a functor from the Fukaya category of $X$ to
coherent sheaves on the rigid analytic mirror, which one may then try to use
to prove homological mirror symmetry).

The SYZ approach was subsequently extended
to the setting of log Calabi-Yau pairs $(X,D)$, where $X$ is a smooth
K\"ahler manifold and $D$ is a (reduced, normal crossings) complex hypersurface
in $X$ representing the anticanonical class $-K_X$. Given a suitable
Lagrangian torus fibration on the complement of $D$, one 
first constructs an SYZ mirror to the open Calabi-Yau $X^0=X\setminus D$, 
before analyzing the manner in which the divisor $D$ deforms the
Lagrangian Floer theory of the torus fibers (and hence the geometry of the mirror). 
This deformation is typically described by a regular function $W\in \O(X^\vee)$ called {\em superpotential}, so
that the SYZ mirror of $X$ (or more accurately, of the pair $(X,D)$) is a {\em Landau-Ginzburg model} 
$(X^\vee,W)$. The superpotential $W$ records the fact that the torus fibers, while
unobstructed in $X\setminus D$, are only {\em weakly unobstructed} as
objects of the Fukaya category of $X$, i.e.\ the Floer-theoretic
obstruction $\m_0\in CF(L,L)$ is a scalar multiple $W\cdot 1_L$ of the
identity, where $W$ is a weighted count of Maslov index 2
holomorphic discs with boundary on $L$. See
e.g.\ \cite{Au07} for an informal overview, and \cite{AAK,Yuan} for a more
up-to-date perspective. (We briefly review the main ingredients in 
\S \ref{ss:syzreview} below.)

The situation is simplest when the Lagrangian torus fibers do not bound any
holomorphic discs in $X\setminus D$, and $D$ is numerically effective (nef).
The fibers are assumed to have vanishing Maslov class in $X\setminus D$, 
so the Maslov
index of a disc is equal to twice its intersection number with $D$, and
the simplest holomorphic discs (intersecting $D$ just once) have Maslov index 2.
The prototypical setting where these assumptions are satisfied is when $X$ is 
toric Fano and $D$ is the toric
anticanonical divisor. The mirror $X^\vee$ is then an algebraic torus
(parametrizing rank 1 local systems on the fibers of the toric moment map), and it
follows from an explicit classification of Maslov index 2 discs bounded 
by the fibers that $W\in\O(X^\vee)$ is a Laurent polynomial determined
combinatorially by the moment polytope \cite{ChoOh,Au07,FO3toric}.
The next simplest case is that of semi-Fano toric varieties,
when the toric anticanonical divisor $D$ is nef but not necessarily
ample. In this case, the coefficients of the Laurent polynomial
$W$ are modified by the contributions of
nodal configurations consisting
of a Maslov index 2 disc in $X$ together with one or more rational curves 
with $c_1(X)\cdot C=0$ contained in the toric divisor $D$. 
The first example in which these contributions were determined explicitly is the Hirzebruch surface $\FF_2$,
i.e.\ the total space of the $\CP^1$-bundle $\PP(\O_{\CP^1}\oplus \O_{\CP^1}(-2))$
over $\CP^1$ \cite{Au09,FO3F2}. General results were subsequently obtained by Chan et
al.\ using comparisons between open and closed Gromov-Witten invariants; see e.g.\ \cite{Chan,CLL,ChanLau,CLLT}.

Outside of the toric setting, the Lagrangian torus fibration $\pi:X^0\to B$
typically has singular fibers, and
the geometric picture is complicated by the presence of holomorphic discs of 
Maslov index 0. (Still assuming $D$ to be nef, these are precisely the
discs which do not intersect $D$). The fibers of $\pi$ which bound such discs typically
lie along (a small neighborhood of) a union of {\em walls} of 
codimension 1 in $B$. There is a discontinuity in the Floer-theoretic
behavior of the fibers of $\pi$ on 
either side of a wall, due to
bubbling of Maslov index 0 discs; nonetheless, it follows from deep
results of Fukaya et al.\ \cite{FO3book} that, across each wall,
the moduli spaces of local systems on the fibers can be glued together via
a suitable analytic coordinate change (the {\em wall-crossing
transformation}) to construct
a moduli space of objects of the Fukaya category of $X^0$ supported on the
fibers of $\pi$, i.e.\ the mirror $X^\vee$ \cite{Au07,AAK,Tu,Yuan}. 
In general there may be an infinite collection of walls,
possibly covering a dense subset of $B$, so that $X^\vee$ cannot be
described explicitly but rather arises as the limit of an inductive
construction \cite{KS,GS}. 

While the positions of the walls in $B$ depend on the choice of 
complex structure, near the {\em large complex structure limit} 
(also known as {\em tropical limit}) the whole process
can be understood combinatorially in terms of tropical geometry:
$B$ carries an integral affine structure (outside of the locus
$B^{sing}$ of singular fibers of $\pi$), and the
{\em scattering diagram}, i.e.\ the set of walls and the corresponding 
wall-crossing transformations, can be determined
via an inductive process first proposed by Kontsevich and 
Soibelman \cite{KS}, based on {\em consistency} of the scattering diagram,
i.e.\ the requirement that the wall-crossing transformations must satisfy
the cocycle property around each codimension 2 locus where walls intersect.
This approach allows one to bypass symplectic geometry
altogether: the Gross-Siebert approach to mirror symmetry starts from a
toric degeneration to construct a tropical manifold $B$, its scattering 
diagram, and a mirror; see e.g.\ \cite{GS}, \cite{GHK}, etc.

Under the assumption that $D$ is nef, discs (or stable discs, i.e.\ nodal
unions of discs and spheres)
which intersect $D$ have Maslov index at least two, so $X$ and $X^0$
have the same scattering diagram, and the mirror of $(X,D)$ is a
Landau-Ginzburg model $(X^\vee,W)$ where $X^\vee$ is entirely determined
by the geometry of $X^0$.
The main point of this paper is to show
that {\em this generally fails to hold 
when $D$ is not nef,} even in examples where the geometry of SYZ
fibrations is well understood. Namely, if $D$ contains rational curves with
$c_1(X)\cdot C<0$, then:

\begin{enumerate}
\item the scattering diagram for $(X,D)$ may contain additional walls
compared to that for $X^0$, or the wall-crossing
transformations for $(X,D)$ may differ from those of $X^0$;
\item in the presence of discs of negative Maslov index, the scattering
diagram for $(X,D)$ may be {\em inconsistent}, so that the wall-crossing
transformations defining $X^\vee$ no longer satisfy the cocycle condition.
\end{enumerate}
While the example we give below is mostly a proof of concept, this has
significant implications. For instance, the construction of
Landau-Ginzburg mirrors for general hypersurfaces in toric varieties
given in \cite{AAK} may require modifications when the stated
assumptions about Chern numbers of rational curves do not hold; more
generally, it is not quite clear which classes of varieties 
should be expected to admit genuine Landau-Ginzburg B-model 
mirrors, rather than deformed LG models of the sort we discuss below. 
By contrast, this 
issue does not seem to affect the other 
direction of homological mirror symmetry: forthcoming work of the author with
Abouzaid (the sequel to \cite{AA}) 
is expected to prove that the (suitably defined) Fukaya categories of the Landau-Ginzburg
A-models given by the construction in \cite{AAK} are indeed equivalent to
the derived categories of the corresponding hypersurfaces, without
Chern class restrictions.

\begin{remark*}
Our results do not contradict in any way the recent work of Gross and
Siebert \cite{GSnew} (see also Keel and Yu \cite{KeelYu}) constructing a
canonical scattering diagram for log Calabi-Yau pairs $(X,D)$ and proving
its consistency. Namely, Gross-Siebert's scattering diagram
only includes Maslov index 0 discs which are contained in $X\setminus D$,
and determines the SYZ mirror of $X\setminus D$; whereas we are studying 
SYZ mirror symmetry for $X$, whose scattering diagram also involves Maslov 
index 0 configurations consisting of a disc in $X$ together with one or more
rational curves in $D$.
\end{remark*}

\subsection{A log Calabi-Yau 4-fold with an inconsistent scattering diagram}
Our main example is the following. Let $K_{\CP^1}=\O_{\CP^1}(-2)$ be the
total space of the canonical bundle of $\CP^1=\C\cup\{\infty\}$.
The toric mirror Landau-Ginzburg
model of $\C^2\times K_{\CP^1}$ is (a domain in) the algebraic torus $(\K^*)^4$ (where
$\K$ is the nonarchimedean field over which we define the Fukaya category, 
say the Novikov field over $\C$ for concreteness) with coordinates
$(z_1,\dots,z_4)$, equipped with the
superpotential $$W=z_1+z_2+(1+q^2+qz_3+qz_3^{-1})z_4,$$ where $q\in \K^*$ is a
constant determined by the choice of K\"ahler form (namely, $q^2$ is
the Novikov weight of the zero section $C_0\subset K_{\CP^1}$).

\begin{theorem}\label{thm:main}
Let $X$ be the blowup of $\C^2\times K_{\CP^1}$ at 
$H_0=\C\times\{1\}\times L_0$ and $H_\infty=\{1\}\times\C\times L_\infty$,
where $L_0$ and $L_\infty$ are the fibers of $K_{\CP^1}$ over $0$ and
$\infty\in \CP^1$ respectively, and let $D$ be the proper transform
of the toric anticanonical divisor of $\C^2\times K_{\CP^1}$. Equip $X$
with a suitable K\"ahler form. Then
$X\setminus D$ carries a fibration by Lagrangian tori of vanishing Maslov
class, whose SYZ mirror consists of four charts which are domains
in $(\K^*)^4$, with superpotentials
\begin{eqnarray}\label{eq:Wmain}
\nonumber W_{--}&=&z_1+z_2+(1+q^2+qz_3+qz_3^{-1})z_4,\\
\nonumber W_{-+}&=&z_1+z_2(1+qq'z_4+q'z_3z_4)+(1+q^2+qz_3+qz_3^{-1})z_4,\\
W_{+-}&=&z_1(1+qq''z_4+q''z_3^{-1}z_4)+z_2+(1+q^2+qz_3+qz_3^{-1})z_4,\\
\nonumber W_{++}&=&z_1(1+qq''z_4+q''z_3^{-1}z_4)+z_2(1+qq'z_4+q'z_3z_4)+q'q''z_1z_2z_4+\\
\nonumber &&\quad +(1+q^2+qz_3+qz_3^{-1})z_4,
\end{eqnarray}
where $q',q''\in\K^*$ are suitable constants. These charts are
glued pairwise by coordinate transformations which preserve $z_3,z_4$ and
act on $z_1,z_2$ by
\begin{align}\label{eq:wallcrossmain}
\nonumber \varphi_{-0}(z_1,z_2)&=(z_1,z_2(1+qq'z_4+q'z_3z_4)),& \varphi_{-0}^*(W_{--})=W_{-+},\\
\nonumber \varphi_{+0}(z_1,z_2)&=(z_1,z_2(1+qq'z_4+q'z_3z_4+q'q''z_1z_4)),& \varphi_{+0}^*(W_{+-})=W_{++},\\
\varphi_{0-}(z_1,z_2)&=(z_1(1+qq''z_4+q''z_3^{-1}z_4),z_2),& \varphi_{0-}^*(W_{--})=W_{+-},\\
\nonumber \varphi_{0+}(z_1,z_2)&=(z_1(1+qq''z_4+q''z_3^{-1}z_4+q'q''z_2z_4),z_2),& \varphi_{0+}^*(W_{-+})=W_{++}.
\end{align}
\end{theorem}

The wall-crossing transformations \eqref{eq:wallcrossmain} are {\em
inconsistent}, in the sense that $$\varphi_{-0}\circ \varphi_{0+}\neq
\varphi_{0-}\circ \varphi_{+0}.$$ This inconsistency arises from the presence of
a codimension 2 locus in the base of the SYZ fibration over which the
fibers bound stable nodal discs of Maslov index $-2$. Indeed, the cocycle property
for wall-crossing transformations is equivalent to the statement that Maslov index 0 discs can only break into 
unions of Maslov index 0 discs; whereas in our example they can also
degenerate to the union of discs of Maslov indices 2 and $-2$. We expect
this to be a general feature of mirror symmetry in settings where the
non-negativity of Maslov index cannot be guaranteed.

The proof of Theorem \ref{thm:main} is given in Section \ref{s:4dex}; the
main new ingredient compared to previous calculations 
on blowups of toric varieties (see in particular \cite{AAK}) is a
study of the contributions of stable nodal configurations consisting of
a holomorphic disc in $X$ together with a rational curve in $D$. 

\begin{remark}
By contrast, the construction of the SYZ mirror of $X\setminus D$ involves 
the same four charts, but the wall-crossing transformations have simpler
expressions:
\begin{eqnarray}\label{eq:wallcrossopen}
\varphi^o_{-0}=\varphi^o_{+0}:&&(z_1,z_2)\mapsto (z_1,z_2(1+q'z_3z_4))\\
\nonumber
\varphi^o_{0-}=\varphi^o_{0+}:&&(z_1,z_2)\mapsto (z_1(1+q''z_3^{-1}z_4),z_2).
\end{eqnarray}
(These are determined by Maslov index 0 discs in $X\setminus D$, whereas the additional terms in \eqref{eq:wallcrossmain} 
correspond to Maslov index 0 configurations with sphere components in $D$.)
The formulas \eqref{eq:wallcrossopen} match the consistent scattering diagram constructed 
by Gross-Siebert \cite{GSnew} for the mirror of $X\setminus D$.
\end{remark}

We also give in \S \ref{ss:compactified} the analogous formulas for the
mirror of a compact example, namely the projective log Calabi-Yau $(\bar{X},\bar{D})$
obtained from $(X,D)$
by compactifying $\C$ to $\CP^1$ and $K_{\CP^1}$ to the Hirzebruch
surface $\FF_2$: namely $\bar{X}$ is the blowup of $\CP^1\times\CP^1\times \mathbb{F}_2$ 
at $\bar{H}_0=\CP^1\times \{1\}\times \bar{L}_0$ and
$\bar{H}_\infty=\{1\}\times \CP^1\times \bar{L}_\infty$, 
where $\bar{L}_0$ and $\bar{L}_\infty$ are the fibers of the projection
from $\FF_2$ to $\CP^1$ over $0$ and $\infty$,  and
$\bar{D}$ is the proper transform of the toric anticanonical divisor of
$\CP^1\times\CP^1\times\FF_2$.

\subsection*{The mirror as a deformed Landau-Ginzburg model}
Even though the mirror in Theorem \ref{thm:main} is no longer a
Landau-Ginzburg model, setting $q'q''=0$ in the formulas
\eqref{eq:Wmain}--\eqref{eq:wallcrossmain}
(i.e., discarding the terms $q'q''z_1z_2z_4$ in $W_{++}$, $\varphi_{+0}$, and
$\varphi_{0+}$)
cures the inconsistency in \eqref{eq:wallcrossmain} and gives a well-defined 
Landau-Ginzburg model $(X^\vee,W)$, of which the mirror in Theorem
\ref{thm:main} can be viewed as a {\em deformation}.
By a result of Lin and Pomerleano \cite[Theorem 3.1]{LP},
the Hochschild cohomology of the category of matrix factorizations of 
$(X^\vee, W)$ is the hypercohomology of the complex of
sheaves $(\Lambda^* T_{X^\vee},\iota_{dW})$ on $X^\vee$.  
We claim that the first-order deformation in Theorem \ref{thm:main} can be viewed
as a class in $$HH^*(MF(X^\vee,W))=\mathbb{H}^*(
X^\vee,(\Lambda^* T_{X^\vee},\iota_{dW}))$$
determined by the contributions of holomorphic discs of Maslov index $-2$ in $X$.

The Maslov index $-2$ discs bounded by the Lagrangian
torus fibers sweep a complex codimension 2 locus in $X$, namely
$\{1\}\times\{1\}\times K_{\CP^1}$, with a Floer-theoretic weight equal to
$q'q''z_4$ mod higher order terms (see \S \ref{s:4dex}). Thus, the
first-order deformation induced by these discs is a
2-cocycle on the base of the fibration with values in the second cohomology
of the fiber, hence dually in $\Lambda^2 T_{X^\vee}$, whose value on the
relevant overlap of coordinate charts is
\begin{equation}\label{eq:maslov-2ex}
q'q''z_4\,\partial_{\log z_1} \wedge \partial_{\log z_2}.
\end{equation}
This element $w^{(2)}\in H^2(X^\vee,\Lambda^2 T_{X^\vee})$ is not closed under
$\iota_{dW}$, but it can be completed to a Hochschild cocycle in
$HH^{even}(MF(X^\vee,W))$ by adding to it a 1-cochain $w^{(1)}$ with values
in $T_{X^\vee}$, whose \v{C}ech coboundary cancels out
$\iota_{dW}(w^{(2)})$; meaning that the value of $\delta w^{(1)}$ on the overlap of coordinate charts
is the vector field
$$\iota_{dW}(q'q''z_4\,\partial_{\log z_1} \wedge \partial_{\log z_2})=
q'q''z_4(z_1\partial_{\log z_2}-z_2\partial_{\log z_1}),$$
which is exactly the inconsistency in \eqref{eq:wallcrossmain}.
Specifically, we can match \eqref{eq:wallcrossmain} by setting
$$w^{(1)}_{-0}=w^{(1)}_{0-}=0,\quad w^{(1)}_{+0}=q'q''z_1z_4\partial_{\log
z_2},\quad \text{and}\quad w^{(1)}_{0+}=q'q''z_2z_4\partial_{\log z_1}.$$
Finally, cancelling out $\iota_{dW}(w^{(1)})$ in turn forces one to also 
add a 0-cochain $w^{(0)}$ with values in $\O_{X^\vee}$, namely we take
$$w^{(0)}_{--}=w^{(0)}_{-+}=w^{(0)}_{+-}=0 \quad \text{and}\quad
w^{(0)}_{++}=q'q''z_1z_2z_4.$$

\subsection{A family Floer perspective}
The above example shows that the construction of 
SYZ mirrors in the presence of discs of negative Maslov index requires a
change of perspective from the usual approach.
In Section \ref{s:counting} we begin a general (but informal) exploration of the geometry
of SYZ mirror symmetry in the setting considered here from the perspective of family Floer
homology.

Consider as before a Lagrangian torus fibration $\pi:X^0\to B$ on
$X^0=X\setminus D$ whose fibers
$F_b=\pi^{-1}(b)$ have vanishing Maslov class in $X^0$ and are weakly unobstructed in $X$. 
Denote by $X^{\vee0}$ the
{\em uncorrected}\/ SYZ mirror of the smooth locus, a rigid analytic
space whose points correspond to unitary rank 1 local systems on the smooth 
fibers of $\pi$. More precisely, we restrict ourselves to a simply connected subset
$B^0$ of the smooth locus $B\setminus \mathrm{critval}(\pi)$, so as to ignore
the issues of compactification over the singular 
fibers of $\pi$ and consistency around the singular fibers, which are
largely orthogonal to our discussion.

The pushforward of the sheaf of analytic functions on $X^{\vee0}$ under
the rigid analytic torus fibration $\pi^\vee:X^{\vee0}\to B^0$ 
defines a sheaf $\O_{an}$ over $B^0$, which
is a certain completion of a local system over $B^0$ whose fiber at $b$ is $\K[H_1(F_b)]$.
(This is just a fancy way of saying that rigid analytic functions on
affinoid domains in $X^{\vee 0}$ are given by Laurent series which satisfy
appropriate convergence conditions.)

Moduli spaces of pseudo-holomorphic discs in $X$ with boundary in the
fibers $F_b$ (where $b$ is allowed to vary over $B^0$) determine 
$A_\infty$-operations $\{\m_k\}_{k\geq 0}$ not just on Floer cochains of a fixed fiber $F_b$
with coefficients in $\K[H_1(F_b)]$, but also on cochains on
$\pi^{-1}(B^0)\subset X^0$ with
coefficients in the pullback of $\O_{an}$, or equivalently, via K\"unneth
decomposition,
on 
\begin{equation}\label{eq:famfloercomplex}
\mathfrak{C}=\bigoplus_{i,j} \CC^{i,j}:=\bigoplus_{i,j} C^i(B^0\,;\,C^j(F_b)\,\hat\otimes\,\O_{an}).
\end{equation} The precise nature of these
cochains depends on the chosen model for Lagrangian Floer theory. Under
very strong transversality assumptions on evaluation maps, a convenient model
consists of an enlargement of differential forms to include currents of
integration along smooth submanifolds (cf.\ \S \ref{sss:singdiffforms}). While it is likely that these
assumptions can be lifted by working with Kuranishi structures, it seems
technically easier to work with a Morse-theoretic model, such as the one we describe in \S\S
\ref{sss:morse}--\ref{sss:adapted} (see also Keeley Hoek's thesis \cite{Hoek} for a more detailed
treatment).

\begin{definition}\label{def:wfunobs}
The Floer complex
$\mathfrak{C}$ is {\em weakly family unobstructed} if $\m_0$ 
can be expressed as a sum of degree $i$ cochains with values in degree $i$
cocycles,
\begin{equation}\label{eq:alphai}
\alpha^{(i)}\in C^i(B^0\,;\,Z^i(F_b)\,\hat\otimes\,\O_{an}), \quad
i=0,1,\dots
\end{equation}
\end{definition}

This definition is quite restrictive (exactly how much depends on the precise
model chosen for cochains) and clearly not satisfied by all SYZ fibrations, 
but we conjecture that weak family unobstructedness should arise
in SYZ mirror symmetry from the existence of a degeneration of the complex
structure on $X$ to the tropical limit, possibly after correction by a
suitably chosen weak family bounding cochain (see below). Indeed, in the 
tropical limit one expects that
the moduli spaces of holomorphic discs which
sweep loci of real codimension $2i$ inside $X$ 
should concentrate along ``walls'' of codimension
$i$ inside $B^0$. 

There is a natural bracket of degree $-1$ on
$H^*(F_b)\otimes \K[H_1(F_b)]$, defined by
\begin{equation}\label{eq:HF-bracket}
\{z^\gamma\, \alpha, z^{\gamma'} \alpha'\}=z^{\gamma+\gamma'}\,
\bigl(\alpha\wedge (\iota_\gamma \alpha')+(-1)^{|\alpha|}
(\iota_{\gamma'}\alpha)\wedge \alpha'\bigr)
\end{equation}
for all $\alpha,\alpha'\in H^*(F_b)$ and $\gamma,\gamma'\in H_1(F_b)$.
Extending to the completion $H^*(F_b)\,\hat{\otimes}\,\O_{an}$ and 
combining with the cup-product on $B^0$, this determines a bracket on
$C^*(B^0;H^*(F_b)\,\hat{\otimes}\,\O_{an})$, symmetric on even degree
elements, which we again denote by $\{\cdot,\cdot\}$.

\begin{conj}\label{conj:mastereq}
For SYZ fibrations on log Calabi-Yau varieties near the tropical limit, 
there exists a model of the family Floer complex $\CC$ for which
$\m_0$ can be expressed as an element in $\bigoplus
C^i(B^0;H^i(F_b)\,\hat\otimes\,\O_{an})$ and satisfies, up to sign, the 
{\em master equation}
\begin{equation}\label{eq:mastereq_m0}
\delta \m_0=\frac12 \{\m_0,\m_0\},
\end{equation}
where $\{\cdot,\cdot\}$ is the bracket defined by \eqref{eq:HF-bracket}
and $\delta$ is the differential on cochains on $B^0$.
\end{conj}

\begin{remark}
We expect \eqref{eq:mastereq_m0} to hold whenever
the moduli spaces of holomorphic discs underlying $\m_0$ behave like closed
manifolds, as a consequence of (a family version of) the master equation in Floer theory with
free loop space coefficients \cite{Floop,Irie}; see Section~\ref{ss:mastereq}.
So  the expectation that this happens for SYZ fibrations near 
the tropical limit is the geometric content of the conjecture. 
However, it is quite possible that the conjecture is too strong as stated,
and that the $A_\infty$-structure on $\CC$ may need to be deformed
by a suitable ``weak family bounding cochain'' $\b\in\CC_{>0}\subset \CC$ (the
subspace of elements whose components have positive Novikov valuation
everywhere), in order for the deformed $\m_0$ term
$$\m_0^\b:=\m_0+\m_1(\b)+\m_2(\b,\b)+\dots\in \CC_{>0}$$
to satisfy the requirements of the conjecture. We note that, even when
$\m_0^\b$ satisfies weak family unobstructedness, there is no geometric
reason for the master equation to hold; rather, it needs to be imposed as
an extra requirement on $\b$. It also seems natural to require $\b$ to
vanish outside of a neighborhood of the 
walls in $B^0$ (when there are infinitely many walls, this statement should
be understood order by order). 

On the other hand, we sketch in Section \ref{ss:boosted} a possible approach
to the master equation in a Morse-theoretic setup, via a deformation of
the moduli space of treed holomorphic discs.
\end{remark}

Via the isomorphism $H^i(F_b,\R)\simeq \Lambda^i H^1(F_b,\R)\simeq\Lambda^i T_B$, 
an element of $H^i(F_b)\,\hat\otimes\,\O_{an}$ naturally determines a
section of $\Lambda^i T_{X^{\vee0}}$ over $(\pi^\vee)^{-1}(b)$, where we again
denote by $X^{\vee0}$ the {\em uncorrected} mirror, equipped with the
rigid analytic torus fibration
$\pi^\vee:X^{\vee0}\to B^0$ (locally modelled on the valuation map, and dual to
$\pi$). Hence, under the assumption of weak family unobstructedness,
the components $\alpha^{(i)}$ of $\m_0$ determine elements
\begin{equation}\label{eq:Wi}
W^{(i)}\in C^i(X^{\vee0},\Lambda^i T_{X^{\vee0}}),
\end{equation}
which encode the instanton corrections to the geometry of $X^{\vee0}$; see
Section \ref{ss:instantoncorr}.
We denote by $\mathbb{W}=\sum_{i\geq 0} W^{(i)}\in C^*(X^{\vee0},\Lambda^*
T_{X^{\vee0}})$ the sum of these terms.
The master equation for $\m_0$ can be transcribed (by an easy argument,
cf.\ \S \ref{ss:instantoncorr}) into an analogous
identity for $\mathbb{W}$:

\begin{proposition}\label{prop:Wi}
If $\m_0$ satisfies \eqref{eq:mastereq_m0}, then
$\mathbb{W}=W^{(0)}+W^{(1)}+\dots \in C^*(X^{\vee0},\Lambda^*T_{X^{\vee0}})$
satisfies
\begin{equation}\label{eq:mastereq}
\delta \mathbb{W}+\frac12 [\mathbb{W},\mathbb{W}]=0,
\end{equation}
where $\delta$ is the differential on cochains and $[\cdot,\cdot]$ is the
bracket induced by the cup-product
and the Schouten-Nijenhuis bracket.
\end{proposition}

Equation \eqref{eq:mastereq} is equivalent to the property that the operator 
$\delta+[\mathbb{W},\cdot]$ squares to zero; in particular, the components
of $\mathbb{W}$ satisfy the equations
\begin{eqnarray}
\label{eq:mastereq0} &(\delta+[W^{(1)},\cdot]) \,W^{(0)} &= \ 0,\\
\label{eq:mastereq1} &(\delta+[W^{(1)},\cdot])^2 &= \ 
[\iota_{dW^{(0)}}(W^{(2)}),\cdot],\\
\label{eq:mastereq2} &(\delta+[W^{(1)},\cdot])\, W^{(2)} &= \ \iota_{dW^{(0)}}(W^{(3)}),
\end{eqnarray}
and so on. 
The geometric interpretation of these equations depends on
the chosen model for cochains in the above discussion, though in
all cases $W^{(1)}$ can be viewed as a deformation of the analytic
structure of $X^{\vee0}$, and \eqref{eq:mastereq0} states that $W^{(0)}$ 
is analytic with respect to the deformed structure, even as
\eqref{eq:mastereq1} measures the failure of $\delta+W^{(1)}$ to genuinely
equip the corrected mirror with an analytic structure.

If we view $\m_0$ as an element of the de Rham complex
$(\Omega^*(B^0,H^*(F_b)\,\hat\otimes\,\O_{an}),d)$
of differential forms on $B^0$
with coefficients in the sheaf $H^*(F_b)\,\hat\otimes\,\O_{an}$, 
then $\mathbb{W}$
ends up being an element of $(\Omega^{0,*}(X^{\vee0},\Lambda^* T_{X^{\vee0}}),d'')$,
the tropical Dolbeault complex of differential forms 
on $X^{\vee0}$ with coefficients in polyvector fields.
(See \cite{CLD,Jell} for a general construction; the version we need here
is significantly simpler because our forms are pulled back from the fixed
tropicalization $\pi^\vee:X^{\vee0}\to B^0$.) 

Assuming convergence, we can
view $X^{\vee0}$ as a family of complex manifolds over a punctured disc,
degenerating to the tropical limit. The tropical Dolbeault complex
specializes to the usual Dolbeault complex, and we can then view
$W^{(1)}\in \Omega^{0,1}(X^{\vee0},T_{X^{\vee0}})$ as a deformation of
the complex structure on $X^{\vee0}$ (deforming $\dbar$ to $\dbar+W^{(1)}$).
The equation \eqref{eq:mastereq0} then states that the function
$W^{(0)}:X^{\vee0}\to \C$ is holomorphic with respect to this deformed complex structure;
and \eqref{eq:mastereq1} states that the deformation in general fails
to be integrable, i.e.\ $\dbar+W^{(1)}$ is only an {\em almost-complex
structure}, whose Nijenhuis tensor is required to be equal to $\iota_{dW^{(0)}}(W^{(2)})\in
\Omega^{0,2}(X^{\vee0},T_{X^{\vee0}})$.

On the other hand, if we work with \v{C}ech cochains rather than
differential forms, then we end up with a picture similar to that discussed
above for our main example: $W^{(1)}$ can be viewed as a deformation of 
the gluing transformations used to assemble $X^\vee$ from local affinoid 
charts, \eqref{eq:mastereq0} states that the expressions for $W^{(0)}$ in
these local charts match under the deformed gluing transformations,
and \eqref{eq:mastereq1} states that $\iota_{dW^{(0)}}(W^{(2)})$ measures
the amount by which the deformed gluing transformations fail to
satisfy the cocycle condition. 

These two perspectives on deformed Landau-Ginzburg models ought to be
equivalent; for example it is readily apparent from both
viewpoints that the critical locus of the superpotential $W^{(0)}$ remains an honest analytic
space (since the right-hand side of \eqref{eq:mastereq1} vanishes along it),
even when the deformed total space fails to be one, so that it still makes
sense to try and relate the symplectic geometry of $X$ to the algebraic
geometry of $\mathrm{crit}(W^{(0)})$ in order to establish homological
mirror symmetry.

\begin{remark}
The $A_\infty$-structure on the family Floer complex $\mathfrak{C}$
is a curved deformation (induced by holomorphic discs) of
the classical algebraic structure on
$C^*(B^0;C^*(F_b)\hat\otimes \O_{an})$, which we have just seen can be
compared to the dg-algebra $C^*(X^{\vee 0};\Lambda^*T_{X^{\vee 0}})$
of cochains with values in polyvector fields on the uncorrected mirror.
Even though we expect that the curvature $\m_0$ of the family Floer complex 
$\mathfrak{C}$ determines the required instanton corrections to the geometry of 
$X^{\vee 0}$ (see also Remark
\ref{rmk:gluingphilosophy}), it is not true (even in the simplest
examples) that $\mathfrak{C}$ itself, as constructed in \S \ref{ss:famfloer}, 
describes cochains with values in polyvector fields on the corrected mirror.
Indeed, unlike the latter algebra, $\mathfrak{C}$ has nonzero curvature;
and its differential $\m_1$ does not match with
the desired expression $\delta+\{\m_0,\cdot\}$. On the other hand, work in progress of the author with
Keeley Hoek suggests that a variant of the construction in \S \ref{ss:boosted} can be
used to define an {\em uncurved} algebraic structure on the family Floer
complex that appears to describe the geometry of the corrected mirror.
\end{remark}

\subsection*{Acknowledgements}

I am heavily indebted to Ludmil Katzarkov and Maxim Kontsevich for many
stimulating discussions about deformations of Landau-Ginzburg models,
which directly led to this investigation of their geometric origin in
Lagrangian Floer theory. I am also grateful for the hospitality of
IH\'ES, where most of this work was carried out. This work was
partially supported by NSF grant DMS-2202984 and by the Simons Foundation
(grant \#\,385573, Simons Collaboration on Homological Mirror Symmetry).

\section{A 4-dimensional example}\label{s:4dex}

This section is devoted to the geometric construction of our main example
and proof of Theorem \ref{thm:main}. The geometric setup is
similar to \cite[Section 3]{Au09} and \cite[Sections 3-5]{AAK},
which also deal with SYZ mirror symmetry for blowups of toric varieties.

\subsection{The geometric setup}

Let $K_{\CP^1}=\O_{\CP^1}(-2)$ be the total space of the canonical bundle of
$\CP^1=\C\cup\{\infty\}$, and denote by $L_0$ and $L_\infty$ the fibers of
$K_{\CP^1}$ over $0$ and $\infty$ in $\CP^1$.
 We equip the product $\C^2\times K_{\CP^1}$ with coordinates
$(x_1,x_2,x_3,x_4)$, where $x_1,x_2$ are the standard coordinates
of $\C^2$, $x_3\in \C\cup\{\infty\}$ is a coordinate on $\CP^1$, and $x_4$
is a coordinate in the fibers of $K_{\CP^1}$ in the trivialization given by
the 1-form $d\log x_3$ over $\C^*$. In other terms, the affine chart $\{x_3\neq
\infty\}\subset K_{\CP^1}$ is isomorphic to $\C^2$ with coordinates $(x_3,x_3^{-1}x_4)$, 
while the affine chart $\{x_3\neq 0\}$ is isomorphic to $\C^2$ with coordinates
$(x_3^{-1},x_3x_4)$.

We denote by $X$ the blowup of $\C^2\times K_{\CP^1}$ along
$H_0=\C\times\{1\}\times L_0$, i.e.\ the locus where $x_2=1$ and
$x_3=0$, and along $H_\infty=\{1\}\times\C\times L_\infty$, i.e.\ the
locus where $x_1=1$ and $x_3=\infty$; we denote again by $x_1,\dots,x_4$ the
pullbacks of the coordinates of $\C^2\times K_{\CP^1}$ under the blowup map
$p:X\to \C^2\times K_{\CP^1}$. The $T^2$-action on $\C^2\times
K_{\CP^1}$ rotating the $x_3$ and $x_4$ coordinates leaves $H_0$ and
$H_\infty$ invariant, and hence lifts to $X$.

We equip $X$ with a $T^2$-invariant K\"ahler form $\omega$ constructed as in
\cite[Section 3.2]{AAK}, symplectomorphic to a toric K\"ahler form on
$\C^2\times K_{\CP^1}$ away from a neighborhood of the exceptional divisors $E_0=p^{-1}(H_0)$
and $E_\infty=p^{-1}(H_\infty)$. For example, one can take
$$\omega=p^*(\omega_{\C^2}\oplus \omega_{K_{\CP^1}})
+\frac{i\epsilon'}{2\pi} \partial\bar{\partial}
\left(\chi
\log(|x_2-1|^2+|x_3|^2) \right)
+\frac{i\epsilon''}{2\pi} \partial\bar{\partial}
\left(\chi
\log(|x_1-1|^2+|x_3^{-1}|^2) \right),$$
where $\omega_{\C^2}\oplus \omega_{K_{\CP^1}}$ is a product toric K\"ahler form on
$\C^2\times K_{\CP^1}$ (standard along the first factor), 
$\epsilon',\epsilon''>0$ are the areas of the fibers of the exceptional
divisors $E_0$ and $E_\infty$, and $\chi\log:\R_+\to \R$ is the product of
the logarithm with a suitable cut-off function.

We denote by $\mu=(\mu_3,\mu_4)$ the moment map of the $T^2$-action on $X$
rotating the $x_3$ and $x_4$ coordinates;
away from $E_0\cup E_\infty$ it coincides with the
pullback of the moment map of the chosen toric K\"ahler form on $K_{\CP^1}$,
and they have the same moment polytope $\Delta\subset \R^2$. We normalize
the moment map so that $$\Delta=\{(\xi_3,\xi_4)\in \R^2\,|\,
\xi_4\geq \max(0,|\xi_3|-a)\},$$ where $a>0$ is half the symplectic area of
the zero section of $K_{\CP^1}$. 

For every $(\xi_3,\xi_4)\in \Delta$, 
the reduced space $\mu^{-1}(\xi_3,\xi_4)/T^2$ is canonically identified
with $\C^2$ via projection to the $x_1$ and $x_2$ coordinates.
The reduced K\"ahler form $\omega_{red,(\xi_3,\xi_4)}$ is a product form
in the $(x_1,x_2)$ coordinates,
and coincides with the standard K\"ahler form of $\C^2$ whenever
$\mu^{-1}(\xi_3,\xi_4)$ lies sufficiently far away from $E_0\cup E_\infty$.
Near $E_0$ (which maps to the region of $\Delta$ where $\xi_4\leq -\xi_3-a+
\epsilon'$), the $x_2$-component of the reduced K\"ahler form
differs from the standard area form near $x_2=1$, and similarly near
$E_\infty$ (where $\xi_4\leq \xi_3-a+\epsilon''$), the $x_1$-component
differs from the standard area form near $x_1=1$. The reduced K\"ahler form
is singular along $x_2=1$ for $\xi_4=-\xi_3-a+\epsilon'$ (this corresponds
to a stratum of points with $S^1$ stabilizers where $E_0$
meets the proper transform of $\C\times \{1\}\times K_{\CP^1}$), and 
similarly along $x_1=1$ for $\xi_4=\xi_3-a+\epsilon''$ (where $E_\infty$
meets the proper transform of $\{1\}\times\C\times K_{\CP^1}$).
(The arguments are similar to \cite[Section 4.1]{AAK} and we omit the
details.)

Since the reduced K\"ahler form on $\mu^{-1}(\xi_3,\xi_4)/T^2$ is a product
form on $\C^2$, the product tori $\{|x_1|=r_1,\ |x_2|=r_2\}$ are Lagrangian
in the reduced space, hence their lifts to $\mu^{-1}(\xi_3,\xi_4)\subset X$ 
are $T^2$-invariant Lagrangian submanifolds of $X$, singular when the
lift contains degenerate $T^2$-orbits and smooth otherwise. Hence, we have:

\begin{defprop}\label{def:fibration}
For $(r_1,r_2,\xi_3,\xi_4)\in B:=\R_+^2\times \mathrm{int}(\Delta)$,
denote by $F_{(r_1,r_2,\xi_3,\xi_4)}$ the Lagrangian submanifold of 
$X$ defined by the equations 
$$|x_1|=r_1,\ |x_2|=r_2,\ \mu_3=\xi_3,\ \mu_4=\xi_4.$$
Denoting by $D\subset X$ the proper transform of the union of the toric
divisors of $\C^2\times K_{\CP^1}$, 
$$\pi=(|x_1|,|x_2|,\mu_3,\mu_4):X\setminus D\to B$$
defines a Lagrangian torus fibration on $X\setminus D$, with singular fibers over
\begin{equation}\label{eq:Bsing}
B^{sing}=\{r_2=1,\ \xi_4=-\xi_3-a+\epsilon'\}\cup 
\{r_1=1,\ \xi_4=\xi_3-a+\epsilon''\}\subset B.
\end{equation}
\end{defprop}

The fibers of $\pi$ which lie sufficiently far from the exceptional
divisors, 
i.e., away from 
\begin{equation}\label{eq:Bexc} 
B^{exc}=\{r_2=1,\ \xi_4\leq -\xi_3-a+\epsilon'\}\cup 
\{r_1=1,\ \xi_4\leq \xi_3-a+\epsilon''\}\subset B,
\end{equation}
are lifts to $X$ of product tori in $\C^2\times
K_{\CP^1}$, hence special Lagrangian with respect to the holomorphic volume form
$p^*(\prod d\log x_i)$ on $X\setminus D$ with simple poles along $D$. This implies immediately:

\begin{lemma}\label{l:maslovindex}
The fibers of $\pi$ have vanishing Maslov class in $X\setminus D$, and the
Maslov index of a disc in $X$ with boundary on a fiber of $\pi$ is 
twice its algebraic intersection number with~$D$.
\end{lemma}

\subsection{Discs and spheres}
The next few sections are devoted to the enumerative geometry of stable
holomorphic discs in $X$ with boundary on the fibers of $\pi$. We start
with two lemmas describing the relevant discs and spheres.

\begin{lemma}\label{l:blaschke}
Let $F$ be a fiber of $\pi$ which is the lift to $X$ of a product torus
$\{|x_i|=r_i\}$ in
$\C^2\times K_{\CP^1}$, and let $u:D^2\to X$ be a holomorphic disc with
boundary on $F$. Then:

$(1)$ The components of $p\circ u$ have Blaschke product 
expansions
\begin{align}\label{eq:blaschke}x_1(z)&=e^{i\theta_1}\,r_1\,
\prod_{i=1}^{n_1}\frac{z-\alpha_{i,1}}{1-\bar\alpha_{i,1} z},\quad&
x_2(z)&=e^{i\theta_2}\,r_2\,
\prod_{i=1}^{n_2}\frac{z-\alpha_{i,2}}{1-\bar\alpha_{i,2} z},\\
\nonumber x_3(z)&=e^{i\theta_3}\,r_3\,\prod_{i=1}^{n_3}\left(\frac{z-\alpha_{i,3}}{1-\bar\alpha_{i,3}
z}\right)^{\epsilon_{i,3}},\quad & 
x_4(z)&=e^{i\theta_4}\,r_4\,\prod_{i=1}^{n_3}\frac{z-\alpha_{i,3}}{1-\bar\alpha_{i,3}
z}\,\prod_{i=1}^{n_4} \frac{z-\alpha_{i,4}}{1-\bar\alpha_{i,4}z},
\end{align}
where $e^{i\theta_k}\in S^1$, $\alpha_{i,k}\in D^2$, and 
$\epsilon_{i,3}\in\{\pm 1\}$. 

$(2)$ $u$ is regular, except possibly
if $x_1(z)$ or $x_2(z)$ is constant and equal to $1$ for all $z$.

$(3)$ The Maslov index of $u$ is 
\begin{equation}\label{eq:maslovindex}
\mu(u)=2(n_1+n_2+n_3+n_4-k_0-k_\infty),
\end{equation}
where $k_0$ and $k_\infty$ are the total contact orders of $p\circ u$ with 
$H_0$ and $H_\infty$; in the absence of multiple roots $k_0$ is the number
of $i\in \{1,\dots,n_3\}$ such that
$\epsilon_{i,3}=+1$ and $x_2(\alpha_{i,3})=1$, and $k_\infty$ is the number
of $i$ such that $\epsilon_{i,3}=-1$ and $x_1(\alpha_{i,3})=1$.

\end{lemma}

\proof (1) The Blaschke product expansions follow from the general
classification of holomorphic discs with boundary on $T^n$-orbits in
toric manifolds, see e.g.\ \cite[Theorem 5.3]{ChoOh}; the only
specific feature in our case is that, given our choice of coordinates
on $K_{\CP^1}$, $x_3$ is allowed to have poles,
and $x_4$ must vanish at the zeroes and poles of $x_3$.

(2) Via the projection $p$, moduli spaces
of holomorphic discs in $X$ with boundary on $F$ correspond to moduli spaces
of holomorphic discs in $\C^2\times K_{\CP^1}$ with prescribed contact orders
with $H_0$ and $H_\infty$. For fixed $x_1(z)$ and $x_2(z)$,
requiring $x_3(z)$ to
vanish to given order at a certain roots of $x_2(z)-1$, and/or to have poles
of given order at certain roots of $x_1(z)-1$, cuts out a smooth subvariety 
of the space of possible Blaschke products of given degree for $x_3(z)$,
of the expected codimension except when $x_2(z)-1$ or $x_1(z)-1$ vanishes 
identically. (This is because the conditions amount to independent linear constraints on the
coefficients of the polynomials $\prod_{\epsilon_{i,3}=+1} (z-\alpha_{i,3})$
and $\prod_{\epsilon_{i,3}=-1} (z-\alpha_{i,3})$.)  
The regularity of $u$ then follows from a general
regularity result for holomorphic discs in the toric setting 
\cite[Theorem 6.1]{ChoOh} and from the fact that the prescribed incidence
conditions with $H_0$ and $H_\infty$ define a transversely cut out, smooth
submanifold of the expected codimension.

(3) By Lemma \ref{l:maslovindex}, the Maslov index of $u$ is twice its
intersection number with the divisor $D$; since $D+E_0+E_\infty$ is the
pullback of the toric anticanonical divisor of $\C^2\times K_{\CP^1}$,
the intersection number is given by counting the zeroes and poles of
$x_1,\dots,x_4$, and excluding intersections with the exceptional divisors
$E_0$ and $E_\infty$, which correspond to the intersections of $p\circ u$
with $H_0$ and $H_\infty$.
\endproof

\begin{lemma}
The only simple holomorphic spheres in $X$ are $(1)$ the spheres $S_{(x_1,x_2)}$ given by
the product of a point $\{(x_1,x_2)\}\in\C^2$ with the zero section of
$K_{\CP^1}$, or their proper transforms when $x_1=1$ and/or $x_2=1$, and
$(2)$ the fibers of the projection $p:X\to\C^2\times K_{\CP^1}$ above the
points of $H_0\cup H_\infty$.
\end{lemma}

\proof By the maximum principle, $x_1$, $x_2$ and $x_4$ are necessarily 
constant along any holomorphic map $u:S^2\to X$, and $x_4$ is necessarily zero
since the nonzero levels of $x_4$ are biholomorphic to $\C^2\times \C^*$.
If $x_3$ is nonconstant then we end up with $S_{(x_1,x_2)}$ or a multiple
cover; otherwise the image of $u$ is contained in a fiber of $p$ over the
blown up locus $H_0\cup H_\infty$. \endproof

We note that $S_{(x_1,x_2)}$ has normal bundle $\O\oplus\O\oplus \O(-2)$
for $x_1,x_2\neq 1$, while $S_{(x_1,1)}$ has normal bundle $\O\oplus \O(-1)
\oplus \O(-2)$, similarly for $S_{(1,x_2)}$, and $S_{(1,1)}$ has normal
bundle $\O(-1)\oplus \O(-1)\oplus \O(-2)$. Since the
$\dbar$ operator on $\O(-2)$ fails to be surjective, the curves
$S_{(x_1,x_2)}$ 
are not regular. However, we will see
that the union of $S_{(x_1,x_2)}$ with a holomorphic 
disc that meets the zero section of $K_{\CP^1}$ transversely is regular as a stable disc.

\subsection{Regularity of stable nodal discs}\label{ss:regularity}

Let $C$ be a nodal Rieman surface with boundary (in our case, $C$ will be
a Riemann sphere glued to a disc at an interior point, e.g.\
the origin), and let
$u:C\to X$ be a holomorphic map with boundary on a Lagrangian submanifold
$L$ (in our case a fiber of $\pi$). The first-order deformations 
of $u$ and their obstructions can be analyzed by methods of algebraic geometry, 
following Behrend-Fantechi \cite{Behrend,BehrendFantechi}. Recall that, when
$C$ is smooth and $u$ is an immersion, the
deformations and obstructions are governed by
$H^0$ and $H^1$ of the normal bundle $N_u=u^*TX/TC$, i.e.\ first-order
deformations correspond to holomorphic sections of $N_u$ over $C$ 
with real boundary conditions given by $N_{\partial u}=u^*TL/T\partial C$, 
while obstructions live in $H^1(C,N_u)$.
In the presence of singularities, the dual of the normal bundle is replaced by a complex of
sheaves, and the deformations and obstructions are given by 
$$\mathrm{Ext}^i_C(\{u^*\Omega^1_X \stackrel{du^*}{\longrightarrow}
\Omega^1_C\},\O_C)$$ for $i=0$ and $i=1$;
see \cite[\S $1\frac12$]{PT} and the discussion before Lemma 2.6 in \cite{GHS}.

As noted in \cite{GHS}, things are
simpler if we assume that the two branches of $u$ near each node of $C$ are
immersed and their tangent lines at the node are distinct: then
$\mathcal{H}om_{\O_C}(\{u^*\Omega^1_X \to
\Omega^1_C\},\O_C)$ is isomorphic to a coherent sheaf $\mathcal{N}_u$, the
{\em normal sheaf} of $u$, whose global sections and first cohomology
determine the deformations and obstructions.
The sheaf $\mathcal{N}_u$ has an explicit
description if we assume moreover that the restriction of $u$ to
each component of $C$ is an immersion. In this case, the restriction of 
$\mathcal{N}_u$ to each component of $C$ is the sheaf of
meromorphic sections of the normal bundle with at most a simple pole at 
each node, whose normal direction must be the tangent space to the other
branch of $u$ through the node; the meromorphic sections over the two
branches must additionally satisfy a matching condition at the node, which
we state below.

\begin{lemma}\label{l:defobs}
Let $C=C'\cup_{p} C''$ be a curve with a single node $p$, and 
$u:(C,\partial C)\to (X,L)$ a holomorphic map whose restrictions
$u'=u_{|C'}$ and $u''_{|C''}$ are immersions; assume moreover that 
the tangent lines $du'(T_pC'),du''(T_pC'')\subset T_{u(p)}X$ are distinct,
and denote by $N_{u'}=u'{}^*TX/TC'$ and $N_{u''}=u''{}^*TX/TC''$ the normal
bundles to the two components. Denote by $z',z''$ local coordinates on
$C',C''$ near $p$.

Then the first-order deformations of $u$ (resp.\ the obstruction space)
are the global sections over
$C$ (resp.\ the first cohomology group) of the normal sheaf $\mathcal{N}_u$ (with real boundary conditions along
$\partial C$), which are pairs of sections $$(v',v'')\in
H^0(C',N_{u'}\otimes \O_{C'}(p))\oplus H^0(C'',N_{u''}\otimes \O_{C''}(p))$$
(i.e., meromorphic sections of $N_{u'}$ and $N_{u''}$ with at most simple
poles at $p$) satisfying the following matching conditions:
\begin{itemize}
\item there exists a constant $\lambda\in \C$ such that the polar parts 
of $v'$ and $v''$ are respectively
\begin{equation}\label{eq:deformpolarpart}
v'(z')\sim \frac{\lambda}{z'}\,\frac{\partial u''}{\partial z''}(p)
\qquad\text{and}\qquad
v''(z'')\sim \frac{\lambda}{z''}\,\frac{\partial u'}{\partial z'}(p);
\end{equation}
\item the projections of $v'$ and $v''$ onto $T_pX/(du'(T_pC')+
du''(T_pC''))$ coincide at $p$.
\end{itemize}
\end{lemma}

Lemma \ref{l:defobs} can also be understood from a differential geometric
perspective, since a node-smoothing deformation of $u$ can be viewed as
the restriction to the family of curves $C_t=\{z'z''=\gamma(t)\}$ of
a family of maps $\tilde{u}_t$ from $C'\times C''$ (in our case, $\CP^1\times D^2$)
to $X$. We then have, for $z',z''\neq 0$,
\begin{align*}
\frac{d}{dt}_{|t=0}(\tilde{u}_t(z',\gamma(t)/z'))&=
\frac{\partial \tilde{u}_t}{\partial t}_{|t=0}(z',0)+
\frac{\gamma'(0)}{z'}\frac{\partial \tilde{u}_0}{\partial z''}(z',0),
\quad\text{and}\\
\frac{d}{dt}_{|t=0}(\tilde{u}_t(\gamma(t)/z'',z''))&=
\frac{\partial \tilde{u}_t}{\partial t}_{|t=0}(0,z'')+
\frac{\gamma'(0)}{z''}\frac{\partial \tilde{u}_0}{\partial z'}(0,z'').
\end{align*}
The first term in these expressions, $\partial \tilde{u}_t/\partial t$, 
is a genuine section of $u^*TX$ over $C$
(i.e., a pair of sections of $u'{}^*TX$ and $u''{}^*TX$ whose
values at the node coincide); while the second term has a first-order
pole at the origin, where the leading order term is exactly as in
\eqref{eq:deformpolarpart} with $\lambda=\gamma'(0)$ (the rate at which
the node is getting smoothed by the deformation). Deformations of the map
$u$ are then governed by the kernel
and cokernel of the $\dbar$ operator on sections of $u^*TX$ with
the appropriate behavior at the node; or, after quotienting out by
vector fields on $C$ (reparametrizations by diffeomorphisms),
the kernel and cokernel of the $\dbar$ operator on pairs of sections of the normal 
bundles $N_{u'}$ and $N_{u''}$ (allowed to have a simple pole at $p$ and
satisfying the matching conditions described in Lemma \ref{l:defobs}).

While the above suffices for our purposes, we also refer the reader to
\cite[Section 6]{SSZ} for a related discussion in the framework
of polyfolds.

We now use Lemma \ref{l:defobs} to prove the regularity of certain 
nodal configurations in $X$ with boundary on the fibers of $\pi$.

\begin{lemma} \label{l:z4stablediscs}
Given $(x_1,x_2)\in (\C^*)^2$, 
let $u:C=\CP^1\cup D^2\to X$ be a
stable map with boundary on $F_{(|x_1|,|x_2|,\xi_3,\xi_4)}$ whose
restriction to $\CP^1$ parametrizes the sphere $S_{(x_1,x_2)}$,
and whose restriction to $D^2$ parametrizes a disc of suitable radius in the
$x_4$ coordinate, with constant values of $x_1,x_2,x_3$.

$(1)$ If $x_1,x_2\neq 1$, then $u$ is regular as a stable disc in $X$
with boundary in $F_{(|x_1|,|x_2|,\xi_3,\xi_4)}$.

$(2)$ If $x_1=1$ and $x_2\neq 1$, then $u$ is regular as a stable disc
in $X$ with boundary in the family of fibers $F_{(r_1,|x_2|,\xi_3,\xi_4)}$
where $r_1$ is allowed to vary.

$(2')$ If $x_1\neq 1$ and $x_2=1$, then $u$ is regular as a stable disc
in $X$ with boundary in the family of fibers $F_{(|x_1|,r_2,\xi_3,\xi_4)}$
where $r_2$ is allowed to vary.

$(3)$ If $x_1=x_2=1$, then $u$ is regular as a stable disc
in $X$ with boundary in the family of fibers $F_{(r_1,r_2,\xi_3,\xi_4)}$
where $r_1$ and $r_2$ are allowed to vary.
\end{lemma}

\proof The normal sheaf splits into a direct sum
$\mathcal{N}_u=\mathcal{N}_{u,1}\oplus \mathcal{N}_{u,2}\oplus
\mathcal{N}_{u,34}$, where the first two summands correspond to the 
$x_1$ and $x_2$ directions and $\mathcal{N}_{u,34}$ corresponds to 
deformations inside
$\{(x_1,x_2)\}\times K_{\CP^1}$ (or its proper transform if $x_1$ or $x_2$
is 1). To prove the vanishing of $H^1(C,\mathcal{N}_u)$ (or equivalently
the surjectivity of the appropriate $\dbar$ operator) we consider each
summand separately.

The first summand $\mathcal{N}_{u,1}$ is a holomorphic line bundle over $C$, whose restriction
to $\CP^1$ is $\O$ if $x_1\neq 1$ and $\O(-1)$ if $x_1=1$, and trivial
over $D^2$, with a trivial real line subbundle as boundary condition.
The $\dbar$ operator on the disc component is surjective, with a real
1-dimensional kernel corresponding to constant sections; meanwhile,
the $\dbar$ operator on the $\CP^1$ component is surjective, and for
$x_1\neq 1$ it remains surjective if we
restrict the domain to sections of $\O$ which have a prescribed value at the
node. Thus, $H^1(C,\mathcal{N}_{u,1})=0$ when $x_1\neq 1$. However, for $x_1=1$ the
$\dbar$ operator on the $\CP^1$ component is only surjective if
we consider all sections of $\O(-1)$, without imposing a value at the
node; and the $\dbar$ operator on the disc component is no longer surjective
if we restrict its domain to functions which take a prescribed value at the
node. Thus, regularity fails if we consider $u$ as a disc with boundary on
a fixed fiber of $\pi$. Instead, we relax the boundary condition and
consider deformations of $u$ among discs with boundary on fibers
$F_{(r_1,|x_2|,\xi_3,\xi_4)}$ where $r_1$ is allowed to vary. Modifying
the problem in this way enlarges
the domain of the $\dbar$ operator on the disc component to the space of
complex-valued functions whose imaginary part is constant (rather than zero)
at the boundary, so that surjectivity holds even if we restrict to functions
that take a prescribed value at the node.

The situation is identical for $\mathcal{N}_{u,2}$: we find that
$H^1(C,\mathcal{N}_{u,2})=0$ when $x_2\neq 1$, and for $x_2=1$ we achieve
regularity by relaxing the boundary condition and allowing $r_2$ to vary.

Finally, $\mathcal{N}_{u,34}$ is a sheaf of sections of the normal bundles
to the components of $u(C)$ inside $K_{\CP^1}$ with at most simple poles
at the node and matching residues. The normal bundle to the $\CP^1$
component is $\O(-2)$, so its $\dbar$ operator has a one-dimensional
cokernel; however, the
$\dbar$ operator becomes surjective if we enlarge the domain to allow a
simple pole at the nodal point. Meanwhile, the normal bundle to the
disc component is trivial, with trivial real boundary condition, so the
corresponding $\dbar$ operator is surjective (on honest sections, and
hence also on sections with a fixed polar part). Thus
$H^1(C,\mathcal{N}_{u,34})=0$.
\endproof

\begin{lemma} \label{l:z1z2z4stablediscs}
Assume $F_{(r_1,r_2,\xi_3,\xi_4)}$ is the lift to $X$ of a product torus in
$\C^2\times K_{\CP^1}$.

$(1)$ For $r_2>1$, let $u:C=\CP^1\cup D^2\to X$ be a
stable map with boundary on $F_{(r_1,r_2,\xi_3,\xi_4)}$ such that
$u_{|\CP^1}$ parametrizes the sphere $S_{(x_1,1)}$ for some $x_1$
such that $|x_1|=r_1$,
and $u_{|D^2}$ is 
as in \eqref{eq:blaschke} with $n_1=n_3=0$ and $n_2=n_4=1$
(i.e., $x_1$ and $x_3$ are constant while $x_2$ and $x_4$ have degree one), with
$x_2(z)=1$ at the unique point where $x_4(z)=0$. If $x_1\neq 1$ then 
$u$ is regular as a stable disc with boundary on
$F_{(r_1,r_2,\xi_3,\xi_4)}$. If $x_1=1$ then $u$ is regular as a stable
disc with boundary on a family of fibers where $r_1$ is allowed to vary.

$(1')$ Similarly for a stable map with boundary on
$F_{(r_1,r_2,\xi_3,\xi_4)}$ $(r_1>1)$ which is the union of $S_{(1,x_2)}$ with a disc on which
$x_2$ and $x_3$ are constant, $x_1(z)$ and $x_4(z)$ have degree $1$,
and $x_1(z)=1$ at the unique point where $x_4(z)=0$.

$(2)$ For $r_1,r_2>1$, the
union of $S_{(1,1)}$ with a disc on which $x_3$ is constant
while $x_1(z)$, $x_2(z)$ and $x_4(z)$ have degree $1$, and $x_1(z)=x_2(z)=1$
at the unique point where $x_4(z)=0$, is regular as a stable map with boundary on
$F_{(r_1,r_2,\xi_3,\xi_4)}$.
\end{lemma}

\proof The normal sheaf $\mathcal{N}_u$ has a subsheaf $\mathcal{N}_{u,1}\oplus
\mathcal{N}_{u,2}$ corresponding to deformations which take place purely along the $x_1$ and $x_2$
directions. We establish vanishing of the first cohomology separately for
$\mathcal{N}_{u,1}$, $\mathcal{N}_{u,2}$, and the quotient 
$\mathcal{N}_{u,34}:=\mathcal{N}_u/(\mathcal{N}_{u,1}\oplus
\mathcal{N}_{u,2})$.

When $x_1$ is constant along the map $u$, the situation for $\mathcal{N}_{u,1}$ is
exactly as in Lemma \ref{l:z4stablediscs}, and the same argument proves the
vanishing of $H^1(C,\mathcal{N}_{u,1})$ if $x_1\neq 1$, and the regularity
once we allow $r_1$ to vary if $x_1=1$. When $x_1$ has degree 1 on the
disc component of $u$, the restriction of $\mathcal{N}_{u,1}$ to the disc
component is still a trivial holomorphic line bundle,
but now the boundary
condition is given by a family of real lines which rotates by one turn in the
positive direction along the unit circle, namely the line spanned by
$i\,x_1(z)$ at every $z\in \partial D^2$. The $\dbar$ operator on this space
of sections is surjective, and remains surjective even after we restrict the
domain to sections which take a prescribed value at the node. (This can be
checked e.g.\ by comparing the index of the operator and the dimension of
its kernel, which consists of infinitesimal automorphisms of the disc).
This implies the vanishing of $H^1(C,\mathcal{N}_{u,1})$ even when $x_1=1$.
The argument for $\mathcal{N}_{u,2}$ is identical.

Finally, $\mathcal{N}_{u,34}$ can be identified (via projection to the
$(x_3,x_4)$ coordinates) with the normal sheaf of the projection of $u(C)$
to $K_{\CP^1}$, which is exactly as in Lemma \ref{l:z4stablediscs} and whose
first cohomology vanishes by the same argument.
\endproof


Next, we give some constraints on nodal configurations which contribute to
the enumerative geometry of discs in $X$ with boundary on fibers of $\pi$. 

Consider a stable disc $u:C\to X$ with boundary on a fiber of $\pi$ which
does not meet $E_0\cup E_\infty$, i.e., $F=\pi^{-1}(b)$ for $b\in B\setminus
B^{exc}$. Denote by $N_1(u)$ the total intersection
number of $u(C)$ with $p^{-1}(\{0\}\times\C\times K_{\CP^1})$,
i.e., the sum of the degrees $n_1$ in \eqref{eq:blaschke} for the disc
components of $u$ , and by $N_2(u)$ the intersection number with
$p^{-1}(\C\times\{0\}\times K_{\CP^1})$, i.e.\ the sum of the values of
$n_2$ for the disc components. (Note that the $x_1$ and $x_2$ components of 
$p\circ u$ always have Blaschke product expansions, without needing to assume that $F$ is
the lift  of a product torus in $\C^2\times K_{\CP^1}$.) Let $N_{34}(u)$ be the total intersection number of
$u(C)$ with the preimages of the toric divisors of $K_{\CP^1}$; when $F$ is
the lift to $X$ of a product torus,
$N_{34}(u)$ is the sum of the quantities $n_3+n_4$ in \eqref{eq:blaschke} 
for the disc components of $u$. Finally, 
denote by $K_0(u)$ and $K_\infty(u)$ the intersection numbers of
$u(C)$ with $E_0$ and $E_\infty$; i.e.,
$K_0(u)$ is the sum of the quantities $k_0$ in
\eqref{eq:maslovindex} for the disc components, plus the degrees of the
sphere components mapping to the curves $S_{x_1,1}$, minus the degrees of
the sphere components mapping to fibers of the projection $p$ contained in
$E_0$; and similarly for
$K_\infty(u)$.
The Maslov index of $u$ is 
\begin{equation}\label{eq:maslovindexsum}
\mu(u)=2(N_1(u)+N_2(u)+N_{34}(u)-K_0(u)-K_\infty(u)).
\end{equation}

\begin{proposition}\label{prop:Kbounds}
There exist arbitrarily small deformations $J'$ of the complex structure on $X$
such that, given any holomorphic stable disc $u:C\to X$ 
with boundary on a fiber $F$ of $\pi$ that does not meet $E_0\cup E_\infty$,
if $u$ deforms to a $J'$-holomorphic
stable disc then:

$(1)$ The sum of the multiplicities of the spheres $S_{x_1,x_2}$ in $u(C)$
is at most $N_{34}(u)$;

$(2)$ $K_0(u)\leq N_{34}(u)$ and $K_\infty(u)\leq N_{34}(u)$;

$(3)$ $K_0(u)\leq N_2(u)$, except possibly if $x_2$ is constant and equal to 
$1$ on a disc component of $u(C)$; and $K_\infty(u)\leq N_1(u)$, except possibly if
$x_1$ is constant and equal to $1$ everywhere on a disc component of $u(C)$.

\end{proposition}

\proof The $x_4$ coordinate defines a Lefschetz fibration
$f=x_4:K_{\CP^1}\to \C$, whose two critical points both lie in the fiber
$x_4=0$ (which is the union of the zero section of $K_{\CP^1}$ and the lines 
$L_0$ and $L_\infty$). We deform this Lefschetz fibration
slightly so that its two critical values become distinct, and deform the
complex structure to some $J'$ so that the deformed fibration $f':K_{\CP^1}\to \C$ remains
$J'$-holomorphic. (This deformation can be viewed as an open subset of the
deformation of the Hirzebruch surface $\FF_2$ considered in \cite[Section~3.2]{Au09}, so that the deformed total space can be identified with the
complement of a curve of bidegree $(1,1)$ inside $\CP^1\times \CP^1$.)
The projection of $F$ to $K_{\CP^1}$ is disjoint from $f^{-1}(0)$, and so
by choosing the deformation to be sufficiently small we can ensure that it is
also disjoint from the preimage under $f'$ of a small disc containing both 
critical values. We may additionally assume that $L_0$ and
$L_\infty$ remain components of the singular fibers of $f'$ (now living over
the two distinct critical values), or denote by $L'_0,L'_\infty$ the
$J'$-holomorphic deformations of $L_0,L_\infty$ which arise in this manner.

We deform the complex structure on $X$ to the blowup of $(\C^2\times
K_{\CP^1},J_0\oplus J')$ along $H'_0=\C\times\{1\}\times L'_0$ and
$H'_\infty=\{1\}\times\C\times L'_\infty$. By abuse of notation, we again
denote by $J'$ the deformed complex structure on $X$, and by $f':X\to\C$ 
the pullback of $f'$ under the composition of the blowup map 
$p':X\to \C^2\times K_{\CP^1}$ and projection to the second factor.

Now, assume that a holomorphic stable disc $u:C\to X$ with boundary on $F$
deforms to a $J'$-holomorphic stable disc $u':C'\to X$.
The only rational curves in $(X,J')$ lie inside the fibers of $p'$ in the
exceptional divisors $E'_0=p'{}^{-1}(H'_0)$ and
$E'_\infty=p'{}^{-1}(H'_\infty)$, and have nonpositive intersection
numbers with $E'_0$ and $E'_\infty$. Meanwhile, the disc components
of $u'$ have total intersection number $N_{34}(u)$ with the fibers of $f'$
near the origin. Hence, by positivity of intersections with
the components of the singular fibers of $f'$, the intersection numbers of $p'\circ u'$ with
$\C^2\times L'_0$ and $\C^2\times L'_\infty$ are bounded by $N_{34}(u)$.
This has two consequences. First, the sum of the multiplicities of the
spheres $S_{x_1,x_2}$ in $u(C)$ (each of which contributes 1 to the
intersection numbers of $p\circ u$ with $\C^2\times L_0$ and $\C^2\times
L_\infty$) is at most $N_{34}(u)$. Second,
the total contact orders of $p'\circ u'$ with 
$H'_0$ and $H'_\infty$ are at most $N_{34}(u)$, so, after adding
the non-positive contributions of any sphere components,
$[u'(C')]\cdot [E'_0]=K_0(u)$ and
$[u'(C')]\cdot[E'_\infty]=K_\infty(u)$ are bounded by $N_{34}(u)$.

On the other hand, the disc components of $u'$ project to the $x_1$ coordinate as a multisection
of degree $N_1(u)$ over the disc of radius $r_1$, which implies that the
intersection number of $p'\circ u'$ with $\{x_1=1\}$ is bounded by $N_1(u)$.
Therefore, the total contact order of $p'\circ u'$ with $H'_\infty$ is
at most $N_1(u)$, unless $x_1$ is constant and equal to 1 on a component of
$p'\circ u'$; this in turn implies that $[u'(C')]\cdot
[E'_\infty]=K_\infty(u)$ is bounded by $N_1(u)$. The bound $K_0(u)\leq
N_2(u)$ (unless $x_2$ is constant and equal to 1 on a component) is proved
similarly by considering the intersection number of $p'\circ u'$ with
$\{x_2=1\}$. 
\endproof

\begin{remark}
Using basic methods of complex analysis to classify holomorphic
discs in conic bundles (arguing as in \cite{Au07,AuR6,AAK}), 
the deformation considered in the proof of Proposition \ref{prop:Kbounds}
can also be used to give an alternative proof of Lemmas
\ref{l:z4stablediscs}--\ref{l:z1z2z4stablediscs} by explicitly finding
the discs that the various nodal configurations deform to, as well as
another derivation of the superpotential formulas given in Section
\ref{ss:superpot}.
\end{remark}

\begin{corollary} \label{cor:Kbounds}
For $r_1\neq 1$ and $r_2\neq 1$, an arbitrarily small perturbation
of the complex structure ensures that all holomorphic stable discs in $X$ with boundary
on $F_{(r_1,r_2,\xi_3,\xi_4)}$ have Maslov index at least $2$.
\end{corollary}

\proof
Let $u:C\to X$ be a non-constant holomorphic stable disc with boundary on 
$F=F_{(r_1,r_2,\xi_3,\xi_4)}$.
Since $r_1,r_2\neq 1$, $F$ is disjoint from $E_0\cup E_\infty$, and
its projection to $\C^2\times
K_{\CP^1}$ is disjoint from the toric divisors. Therefore,
either one of $x_1,x_2$ is non-constant along $u$, in which case $N_1(u)$
or $N_2(u)$ is positive, or $p\circ u$ is a non-constant disc in
$\{(x_1,x_2)\}\times K_{\CP^1}$ with boundary on a product torus,
in which case $N_{34}(u)$ must be positive.

By Proposition \ref{prop:Kbounds}~(3), the stable discs which survive the
perturbation of the complex structure to $J'$ satisfy $K_0(u)\leq N_2(u)$
and $K_\infty(u)\leq N_1(u)$. (The other possibility, that $x_1$ or $x_2$
is constant and equal to 1 on a disc component, is excluded by the
assumption that $r_1,r_2\neq 1$.) Using \eqref{eq:maslovindexsum}, it
follows that $\mu(u)\geq 2N_{34}(u)$. If $N_{34}(u)>0$ the conclusion
follows. Otherwise, if $N_{34}(u)=0$ then Proposition \ref{prop:Kbounds}~(2)
implies that $K_0(u)=K_\infty(u)=0$, so $\mu(u)=2(N_1(u)+N_2(u))\geq 2$.
\endproof

\subsection{A brief review of SYZ mirror symmetry}\label{ss:syzreview}

Before proceeding further, we recall the construction of the SYZ mirror
of $X$ relative to the anticanonical divisor $D$; see
\cite[Section 2 and Appendix~A]{AAK} and \cite{Yuan} for details.

The construction of the mirror $X^\vee$ starts
from a moduli space of objects of the Fukaya category of $X^0=X\setminus D$ 
consisting of weakly unobstructed fibers of $\pi:X^0\to B$ equipped with
rank 1 unitary local systems. We work over the Novikov field over a field
$k$, say $k=\C$ for concreteness,
$$\K=\Lambda_k=\bigl\{\textstyle\sum a_i T^{\lambda_i}\,|\,a_i\in k,\ \lambda_i\in \R,
\ \lambda_i\to +\infty\bigr\},$$
and recall the unitary subgroup $U_\K=\val^{-1}(0)\subset \K^*$, where the 
valuation map $\val:\K^*\to \R$ is defined by $\val(\sum a_i
T^{\lambda_i})=\min\{\lambda_i\,|\,a_i\neq 0\}$. Unitary rank 1 
local systems over a Lagrangian torus $F_b=\pi^{-1}(b)$ are determined by
their holonomy $\mathrm{hol}\in \hom(\pi_1(F_b),U_\K)=H^1(F_b,U_\K)$, which
enters into the formulas for weighted counts of holomorphic discs in
Lagrangian Floer theory; specifically, a disc with boundary on $F_b$ 
representing the class $\beta\in \pi_2(X,F_b)$ is counted with a weight
$$z^\beta=T^{\omega(\beta)}\mathrm{hol}(\partial\beta)\in \K^*.$$

Over a simply connected subset $P \subset B$ where the fibers of $\pi$ are smooth and do not bound 
any holomorphic discs of Maslov index less than $2$, using isotopies between
the fibers to identify $\pi_2(X,F_b)\simeq \pi_2(X,F_{b'})$ for $b,b'\in P$,
there is a natural analytic structure on $$X^\vee_P:=\bigsqcup_{b\in P} H^1(F_b,U_\K)$$ 
for which the functions $z^\beta \in \O(X^\vee_P)$ are analytic. 
(Typically one might take $P$ to be a bounded rational convex polyhedral subset, so that $X^\vee_P$ is an affinoid domain.)
$X^\vee_P$ can be identified with a domain in $(\K^*)^n$ by considering the
coordinates $z_i=z^{\beta_i}$ for some choice of classes
$\beta_1,\dots,\beta_n$ such that $\partial \beta_1,
\dots,\partial \beta_n$ are a basis of $H_1(F_b,\Z)$; all other $z^\beta$
are then Laurent monomials in $z_1,\dots,z_n$. 
Moreover, the non-archimedean torus fibration defined by the natural projection 
$X_P^\vee \to P$ is modelled on the valuation map in these coordinates, in
the sense that the diagram
\begin{equation}\label{eq:localmirrorcoords}
\begin{tikzcd}[column sep=huge]
X_P^\vee \arrow[r,"(z_i)_{1\leq i\leq n}"] \arrow[d] & (\K^*)^n \arrow[d,"\val"]\\
P \arrow[r,"(\omega(\beta_i))_{1\leq i\leq n}"] & \R^n
\end{tikzcd}
\end{equation}
commutes. 

The superpotential $W\in \O(X_P^\vee)$ is the coefficient of identity in
the Floer-theoretic obstruction $\m_0$ for fibers of $\pi$ equipped with
unitary rank 1 local system, i.e.\ a weighted count of Maslov
index 2 holomorphic discs with boundary on $F_b$ passing through a generic point of
$F_b$. Namely,
\begin{equation}\label{eq:defW}
W=\sum_{\mu(\beta)=2} n_\beta\,z^\beta,
\end{equation} 
where $n_\beta\in\Z$ is the degree of
the evaluation map $ev:\mathcal{M}_1(F_b,\beta)\to F_b$ from
 the moduli space of holomorphic discs 
with boundary in $F_b$ representing the class $\beta$ and one boundary
marked point, $\mathcal{M}_1(F_b,\beta)$, to the Lagrangian $F_b$
(after fixing suitable orientations of both spaces, and possibly a perturbation to
achieve regularity of the moduli space).

The mirror $X^\vee$ is assembled from the subsets $X_P^\vee$ via
suitable gluing maps: the transition functions between the
affine coordinates in the bottom row of \eqref{eq:localmirrorcoords} are
given by elements of $GL(n,\Z)\ltimes \R^n$, and in the absence of discs 
of Maslov index less than two the local analytic coordinates $(z_i)$ on the
subsets $X_P^\vee\subset X^\vee$ transform by
the corresponding monomial automorphisms of $(\K^*)^n$. However, walls over
which the fibers of $\pi$ bound Maslov index 0 discs induce a modification of the
transition functions between between the portions of $X^\vee$ which correspond to subsets of 
$B$ lying on either side of the wall.  The existence of analytic
(valuation-preserving) coordinate changes which restore the analytic
dependence of Floer theory across the wall 
follows from the work of
Fukaya-Oh-Ohta-Ono on the invariance of Floer cohomology for Lagrangians with
weak bounding cochains \cite{FO3book} (in their language, the wall-crossing
coordinate transformation arises as an induced map on the moduli space of
weak bounding cochains); the new phenomenon we will evidence below, however,
is that in the presence of discs of negative Maslov index these coordinate 
changes need not be path-independent.

\begin{remark}
The above statements about gluing maps between local
charts of $X^\vee$ which lie over different
subsets of $B$ (and, a fortiori, the
composition of such transformations along paths in $B$) require further
explanation, since the local charts corresponding to disjoint subsets of $B$ do not actually overlap in $X^\vee$. 
The key observation, known as Fukaya's trick \cite{Fadic}, is that
when two fibers $F,F'$ of $\pi$ are close enough to be mapped to each other by an
isotopy $\psi$ such that the almost-complex structure $J'=\psi_*J$ is $\omega$-tame,
the Floer theory of $F'$ with respect to $J'$ is related to the Floer theory of $F$ with respect to
$J$ by analytic continuation. This implies that
Floer-theoretic expressions calculated for fibers
of $\pi$ over a given region of $B$ (possibly just a single fiber) can be analytically continued over a
slightly larger region of $B$. We can then form a cover of $B$ by rational
convex polyhedral subsets that overlap nontrivially, with the
understanding that, on the overlaps, the gluing maps
amount to a comparison of the Floer theory of given fibers of $\pi$ with
respect to almost-complex structures that are pulled back along different isotopies. 
In practice, for sufficiently simple examples, such as the one we consider
here, it is often the case that the wall-crossing 
transformations are given by birational maps, and can be composed freely
without worrying about convergence issues. For this reason we do not discuss
the details further, and instead refer the reader to 
\cite[Section 3]{AbICM} and \cite[Section 4]{Yuan} for details.
\end{remark}

A key property that can be used to determine the wall-crossing coordinate 
transformations between the
coordinate charts $X_P^\vee$ is that the local expressions
\eqref{eq:defW} for the superpotential must match and assemble to a
global analytic function $W\in \O(X^\vee)$ (to the extent that the corrected
mirror $X^\vee$ is globally well-defined). More generally, the same property
holds for weighted counts of holomorphic discs in moduli spaces 
that match under wall-crossing after accounting for disc bubbling. 
For example, expressions $\sum
n_\beta z^\beta$ where the sum ranges over classes
$\beta$ with fixed intersection numbers with certain divisors $D_i\subset X$ are also
invariant under wall-crossing, as long
as the intersection numbers with $D_i$ are zero for all Maslov index 0 bubbles.

\subsection{The superpotential: discs of Maslov index $2$}\label{ss:superpot}
We now return to our main example, and prove the first part of Theorem \ref{thm:main}, namely we
determine the superpotential on each chart of the mirror $X^\vee$.
By Corollary \ref{cor:Kbounds}, the Lagrangian tori
$F_{(r_1,r_2,\xi_3,\xi_4)}$ only bound discs of Maslov index at least 2 as
soon as $r_1\neq 1$ and $r_2\neq 1$; hence we work separately over each of
the four domains $P_{--}=\{r_1<1,\ r_2<1\}$, $P_{-+}=\{r_1<1,\ r_2>1\}$,
$P_{+-}=\{r_1>1,\ r_2<1\}$ and $P_{++}=\{r_1>1,\ r_2>1\}$, calculating the
superpotential on each chart and showing that it is given by \eqref{eq:Wmain}.

The very simplest holomorphic discs that we will encounter are the lifts to $X$
of ``standard'' Maslov index 2 holomorphic discs bounded by product tori in
$\C^2\times K_{\CP^1}$. Using the same notations as in Lemma
\ref{l:blaschke}, these are the discs for which $n_1+n_2+n_3+n_4=1$ (i.e.,
one of the $n_i$ is equal to 1 and the others are zero).
For $i\in \{1,2,4\}$, we denote by $\beta_i$ the class of a disc of the
appropriate radius along the $x_i$ coordinate axis, while
the other coordinates $x_j$ for $j\neq i$ are constant, and by $z_i=z^{\beta_i}$
the corresponding Floer-theoretic weight. 
For $i=3$ there are two different
classes of discs with $n_3=1$ and \hbox{$n_1=n_2=n_4=0$}, depending on whether the $x_3$
coordinate has a zero or a pole (i.e., whether the disc intersects
$\C^2\times L_0$ or $\C^2\times L_\infty$). We denote by $\beta_{3,\pm}$
their homotopy classes, and define $z_3=q^{-1}z_4^{-1}z^{\beta_{3,+}}$,
where $q=T^a$. Since $\beta_{3,+} +
\beta_{3,-}=2\beta_4+[S_{(x_1,x_2)}]$, and the Novikov
weight of the zero section of $K_{\CP^1}$ is $T^{2a}=q^2$, we find that
$z^{\beta_{3,\pm}}=q z_3^{\pm 1} z_4$.

We will use $z_1,z_2,z_3,z_4$ as coordinates on each of the four charts that
make up $X^\vee$; we observe that these are the weights of disc classes
$\beta_1,\dots,\beta_4$ (up to a factor of $q$ in the case of $z_3$) whose
boundaries $\partial \beta_i$ correspond to the standard basis of the
first homology of product tori in $\C^2\times K_{\CP^1}$. Moreover, using
the fact that the symplectic area of a disc which is invariant under a Hamiltonian
$S^1$-action is equal to the difference between the moment map values at its
boundary and at its center, one finds that
$\val(z_3)=\xi_3$ and $\val(z_4)=\xi_4$.

We are now ready to determine the formulas for the superpotential on each
of the four charts. Since the counts $n_\beta$ do not vary inside
each of the four regions $P_{\pm,\pm}$, it suffices to carry out the calculation for
fibers of $\pi$ which are lifts of product tori in $\C^2\times K_{\CP^1}$.

 Since we only consider fibers of $\pi$ for which
$r_1\neq 1$ and $r_2\neq 1$, Proposition \ref{prop:Kbounds} implies that all the
Maslov index 2 holomorphic stable discs $u:C\to X$ contributing to the sum \eqref{eq:defW} 
satisfy $K_0(u)\leq \min(N_2(u),N_{34}(u))$ and $K_\infty(u)\leq
\min(N_1(u),N_{34}(u))$. As noted in the proof of Corollary
\ref{cor:Kbounds}, plugging these bounds into \eqref{eq:maslovindexsum} 
it follows that 
\begin{equation}\label{eq:mulowerbound}
\mu(u)\geq 2\max(N_1(u),N_2(u),N_{34}(u)).
\end{equation}
Hence, we only need to consider stable discs for which each of
$N_1(u),N_2(u),N_{34}(u)$ is either $0$ or $1$. Moreover, since equality
must hold in \eqref{eq:mulowerbound}, necessarily
$K_0(u)=\min(N_2(u),N_{34}(u))$ and $K_\infty(u)=\min(N_1(u),N_{34}(u))$.
It is apparent from the proof of Proposition \ref{prop:Kbounds} that for
these equalities to hold, $u$ cannot have any sphere components contained
in the fibers of $\pi$; whereas by Proposition \ref{prop:Kbounds}(1) the
total multiplicities of the spheres $S_{(x_1,x_2)}$ add up to at most
$N_{34}(u)\leq 1$, i.e.\ $C$ contains at most one sphere component, and any
such component must map to some $S_{(x_1,x_2)}$ with degree one.
Furthermore, since each disc component of $u$ 
must separately satisfy the constraints of
Proposition \ref{prop:Kbounds}, the Maslov index 2 configurations we
consider have only one disc component.

With these constraints in hand, we can list all the possible homotopy
classes which may contribute to the superpotential.

{\em Case 1: $N_{34}(u)=0$.} Then there are no sphere components,
$K_0(u)=K_\infty(u)=0$, and $N_1(u)+N_2(u)=1$. Hence 
$x_3,x_4$ are constant along $u$, while one of $x_1,x_2$ is constant and the
other parametrizes a disc of radius $r_i$ parallel to the $x_i$ coordinate
axis. Thus $[u]$ is either $\beta_1$ or $\beta_2$, and its weight is either
$z_1$ or $z_2$. Both of these families of discs are regular and contain
exactly one disc through each point of $F=F_{(r_1,r_2,\xi_3,\xi_4)}$. The
orientation of the moduli space works out as in the classical toric case,
and $n_{\beta_1}=n_{\beta_2}=1$. Summarizing, the contributions of the discs
with $N_{34}(u)=0$ add up to
\begin{equation}\label{eq:Wtermsno34}
z_1+z_2.
\end{equation}

{\em Case 2: $N_{34}(u)=1$, $n_3=0$, $n_4=1$.} (Recall that $N_{34}(u)$ is 
the sum of the $n_3$ and $n_4$ degrees appearing in \eqref{eq:blaschke}.)
Then $x_3$ is constant along the disc component of $u$, while $x_4$ 
parametrizes a disc of the appropriate radius. 

When $N_1(u)=N_2(u)=0$ (i.e., the disc component of $u$ is parallel to 
the $x_4$ coordinate axis and represents the class $\beta_4$), we can
either have just the disc component, or consider its union with
the sphere $S=S_{(x_1,x_2)}$ (for the
same constant values taken by $x_1,x_2$ along the disc); the latter
nodal configuration is regular by Lemma \ref{l:z4stablediscs}(1).
In both cases there is one such configuration through each point of $F$,
and the orientations work out as in the toric case (the contributions of the
sphere component to the linearized Cauchy-Riemann problem amount to complex
linear operators and do not affect signs). Hence
$n_{\beta_4}=n_{\beta_4+[S]}=1$, contributing $(1+q^2)z_4$ to the
superpotential. We will now see that in all other cases (when either
$N_1(u)$ or $N_2(u)$ is non-zero) a sphere component {\em must} be present.

$N_1(u)$ is either zero or one. If it is zero then $x_1$ is constant along
$u$. If $N_1(u)=1$, then we must have $K_\infty(u)=1$ as well; and since the
disc component does not meet $E_\infty$ ($x_3$ has no pole), this forces the
presence of a sphere component mapping to $S_{(1,x_2)}$ for some value of
$x_2$. This in turn implies that $x_1$ must equal 1 at the point of the disc
component where $x_4$ vanishes. Since $x_1$ takes values in the disc of
radius $r_1$, this is only possible if $r_1>1$. After a
suitable reparametrization, the $x_1$ and $x_4$ coordinates along the disc 
component of $u$ can be put in the form
$$x_1(z)=r_1z,\qquad x_4(z)=e^{i\theta}r_4\frac{r_1z-1}{r_1-z}$$
for some $e^{i\theta}\in S^1$, and the sphere $S_{(1,x_2)}$ is attached to
the disc at $z=1/r_1$.

Similarly, $N_2(u)$ is either zero or one; if it is zero then $x_2$ is
constant; if $N_2(u)=1$ then necessarily $r_2>1$, the $x_2$ component
takes the value $1$ at the point of the disc where $x_4$ vanishes, and a 
sphere component $S_{(x_1,1)}$ is attached to the disc at that point.

The case where $N_1(u)$ and $N_2(u)$ are both equal to $1$ does
occur; this requires the sphere component to map to
$S_{(1,1)}$. Hence we must have $(x_1,x_2)=(1,1)$ at the point of the disc 
where $x_4$ vanishes (and necessarily $r_1,r_2$ are both greater than $1$).

The various configurations we have found are precisely those covered by
Lemma \ref{l:z1z2z4stablediscs}; hence they are regular, and one easily 
checks that there is one disc in each family through each point of $F$.
As before, the incidence constraints and contributions from sphere components modify
the linearized Cauchy-Riemann problem by complex linear operators, so that
the evaluation maps again have degree 1.
Thus,
\begin{align*}
n_{\beta_1+\beta_4+[S_{(1,x_2)}]}&=1\qquad \text{for}\ r_1>1\qquad (\text{and}\ 0\ \text{otherwise}),\\
n_{\beta_2+\beta_4+[S_{(x_1,1)}]}&=1\qquad \text{for}\ r_2>1\qquad (\text{and}\ 0\ \text{otherwise}),\\
n_{\beta_1+\beta_2+\beta_4+[S_{(1,1)}]}&=1\qquad \text{for}\ r_1,r_2>1\ \ (\text{and}\ 0\ \text{otherwise}).
\end{align*}
Setting $q'=T^{a-\epsilon'}$ and $q''=T^{a-\epsilon''}$, the Floer weights of 
$S_{(1,x_2)}, S_{(x_1,1)}, S_{(1,1)}$
are respectively $T^{2a-\epsilon''}=qq''$, $T^{2a-\epsilon'}=qq'$, and
$T^{2a-\epsilon'-\epsilon''}=q'q''$. Hence, the contributions of discs with 
$n_3=0$ and $n_4=1$ add up to
\begin{equation}\label{eq:Wtermsonly4}
\begin{cases}
(1+q^2)z_4 & \text{if } r_1<1 \text{ and } r_2<1\\
(1+q^2)z_4+qq'z_2z_4 & \text{if } r_1<1 \text{ and } r_2>1\\
(1+q^2)z_4+qq''z_1z_4 & \text{if } r_1>1 \text{ and } r_2<1\\
(1+q^2)z_4+qq''z_1z_4+qq'z_2z_4+q'q''z_1z_2z_4 & \text{if } r_1>1 \text{ and } r_2>1
\end{cases}
\end{equation}

{\em Case 3: $N_{34}(u)=1$, $n_3=1$, $n_4=0$.} Because the disc component of
$u$ does not meet the preimage of the zero section of $K_{\CP^1}$, there
cannot be any sphere component, and $u$ is as in \eqref{eq:blaschke}.
Moreover, $x_3$ has either a zero or a pole along $u$, but not both, so
$u$ meets at most one of $E_0$ or $E_\infty$. However, for this to happen,
$x_1$ or $x_2$ needs to be non-constant and take the value 1 at the point
where $x_3$ has its pole or zero. There are therefore three subcases.

If $N_1(u)=N_2(u)=0$, then $x_1$ and $x_2$ are constant along $u$, and
$u$ represents one of the classes $\beta_{3,\pm}$ discussed above; arguing
as in the toric case, $n_{\beta_{3,+}}=n_{\beta_{3,-}}=1$.

If $N_1(u)=1$, then $K_\infty(u)=1$, forcing $x_3$ to have a pole and not a
zero; this in turn forces $K_0(u)=0$ and $N_2(u)=0$, i.e.\ $x_2$ is constant
along $u$. Moreover, $x_1$ needs
to take the value $1$ at the pole of $x_3$, which can only happen if
$r_1>1$. Assuming this is the case, after a suitable reparametrization we
can write $$x_1(z)=r_1 z,\quad x_3(z)=e^{i\theta_3}r_3\frac{r_1-z}{r_1z-1},\quad
x_4(z)=e^{i\theta_4}r_4\frac{r_1z-1}{r_1-z}$$
for some $e^{i\theta_3},e^{i\theta_4}\in S^1$. 
There is one such disc through every point of $F$.

If $N_2(u)=1$, then $K_0(u)=1$, forcing $x_3$ to have a zero and not a pole;
hence $K_\infty(u)=0$, $N_1(u)=0$, and $x_1$ is constant along $u$. Moreover
$x_2$ takes the value $1$ at the zero of $x_3$. Such discs can only exist if
$r_2>1$; after a suitable reparametrization they are of the form
$$x_2(z)=r_2 z,\quad x_3(z)=e^{i\theta_3}r_3\frac{r_2z-1}{r_2-z},\quad
x_4(z)=e^{i\theta_4}r_4\frac{r_2z-1}{r_2-z},$$
and there is one such disc through each point of $F$.

Summarizing, the contributions of discs with $n_3=1$ and $n_4=0$ add up to
\begin{equation}\label{eq:Wtermsonly3}
\begin{cases}
qz_3z_4+qz_3^{-1}z_4 & \text{if } r_1<1 \text{ and } r_2<1\\
qz_3z_4+qz_3^{-1}z_4+q'z_2z_3z_4 & \text{if } r_1<1 \text{ and } r_2>1\\
qz_3z_4+qz_3^{-1}z_4+q''z_1z_3^{-1}z_4 & \text{if } r_1>1 \text{ and } r_2<1\\
qz_3z_4+qz_3^{-1}z_4+q''z_1z_3^{-1}z_4+q'z_2z_3z_4 & \text{if } r_1>1 \text{ and } r_2>1.
\end{cases}
\end{equation}

Adding \eqref{eq:Wtermsno34}, \eqref{eq:Wtermsonly4} and
\eqref{eq:Wtermsonly3}, we arrive at the expressions \eqref{eq:Wmain} for
the superpotential on the various charts of $X^\vee$.

\subsection{Wall-crossing: discs of Maslov index $0$ and $-2$}
In this section we study the wall-crossing transformations along which
the coordinate charts $X^\vee_{\pm,\pm}$ corresponding to the domains 
$P_{\pm,\pm}\subset B$ are glued to each other. Our first observation is that,
after a small perturbation of the complex structure as in Proposition
\ref{prop:Kbounds}, Maslov index 0 discs only exist along the walls
$r_1=1$ and $r_2=1$, and are entirely contained in the divisors $\{x_1=1\}$
and $\{x_2=1\}$, while negative Maslov index discs can only exist at
$r_1=r_2=1$. 

\begin{proposition}\label{prop:maslovzero}
Every holomorphic stable disc $u:C\to X$ with boundary on a
smooth fiber $F=F_{(r_1,r_2,\xi_3,\xi_4)}$ of $\pi$ which deforms to a
stable disc for arbitrarily small perturbations of the complex structure 
on $X$ chosen as in Proposition \ref{prop:Kbounds} satisfies the following:

(1) if $u$ has negative Maslov index, then $r_1=r_2=1$;

(2) if $(r_1,r_2)\neq (1,1)$ and $u$ has Maslov index zero, then either
$r_1=1$, in which case $x_1=1$ at every point of $u(C)$, or
$r_2=1$, in which case $x_2=1$ at every point of $u(C)$.
\end{proposition} 

\proof We prove, equivalently, that if $r_1$ and $r_2$ are not both equal to 
1 then $\mu(u)\geq 0$, and if $\mu(u)=0$ then the conclusion of (2) holds. 
There are two cases: either $r_1\neq 1$ or $r_2\neq 1$. The argument is the 
same for both; we give the proof for $r_2\neq 1$.

As in Proposition \ref{prop:Kbounds}, we deform slightly the Lefschetz
fibration $f=x_4:K_{\CP^1}\to\C$ to $f':K_{\CP^1}\to \C$ with two distinct
singular fibers. Denote by
$\Delta=f^{-1}(0)$ the singular fiber of $f$, which is the union of the toric 
divisors of $K_{\CP^1}$, and by $\Delta'_0,\Delta'_\infty\subset K_{\CP^1}$ the two
singular fibers of $f'$, labelled so that $L'_0$ (i.e., $L_0$ or a small
deformation thereof) is a component of $\Delta'_0$ and $L'_\infty$
($L_\infty$ or a small deformation) is a component of $\Delta'_\infty$.
Recall that we deform $X$ to the blowup of $(\C^2\times
K_{\CP^1},J_0\oplus J')$ along $H'_0=\C\times\{1\}\times L'_0$ and
$H'_\infty=\{1\}\times\C\times L'_\infty$.

Because $r_2\neq 1$, the fiber $F$ is disjoint not only from the
anticanonical divisor $D\subset X$ but also from the exceptional divisor
$E_0$. This implies that, for a small enough deformation, it is
also disjoint from the proper transform $Z'_\infty$ of 
$\C^2\times \Delta'_\infty$ under the blowup at $H'_\infty$.
(Indeed, $Z'_\infty$ is a small deformation of the union $Z_\infty$ of the 
proper transform of $\C^2\times \Delta$ and the exceptional divisor $E_0$.)

Furthermore, the anticanonical divisor $D\subset X$ is homologous in the complement of $F$ to the
(non-effective) divisor $D'_-=(\{0\}\times \C\times K_{\CP^1}) + (\C\times \{0\}\times
K_{\CP^1}) + Z'_\infty - E'_0$ (since the proper transform of $\C^2\times
\Delta$ is homologous in $X\setminus F$ to $Z'_\infty-E'_0$). Hence, the Maslov index of the
$J'$-holomorphic deformation
$u':C'\to X$ of the holomorphic disc $u$ is equal to twice its intersection
number with $D'_-$.

If $r_2<1$, then the maximum principle for $|x_2|$ implies that the image of
$u'$ is disjoint from $E'_0$. It then follows from positivity of intersections
between the $J'$-holomorphic curve $u'(C')$ and the other components of $D'_-$
(and the absence of any rational curves intersecting those components
negatively) that $\mu(u')=2[u'(C')]\cdot [D'_-]\geq 0$, and if $\mu(u')=0$ then
$u'(C')$ is disjoint from every component of $D'_-$. 

If $r_2>1$, then we can deform $D'_-$ inside the complement of $F$ to an effective
divisor $D'_+$ which is the sum of three components: 
$\{0\}\times \C\times K_{\CP^1}$, the proper transform of $\C\times \{1\}\times
K_{\CP^1}$ under the blowup at $H'_0$, and $Z'_\infty$. 
It then follows from positivity of
intersections (and the lack of rational curves intersecting $D'_+$
negatively) that $\mu(u')=2[u'(C')]\cdot [D'_+]\geq 0$, and if $\mu(u')=0$
then $u'(C')$ is disjoint from every component of $D'_+$.

Since $\mu(u')=\mu(u)$, we have proved that $\mu(u)$ is non-negative, and if
it is zero then the image of $u'$ is disjoint from the components of $D'_+$
or $D'_-$ depending on the value of $r_2$. 

From now on we assume that $\mu(u)=\mu(u')=0$. Since $u'(C')\cap
D'_\pm=\emptyset$, the image of $u'$ is disjoint from $\{0\}\times \C\times
K_{\CP^1}$ and from $Z'_\infty$; since the deformation from $Z_\infty$ to
$Z'_\infty$ does not cross $F$, the intersection numbers of $u(C)$ with
$\{0\}\times\C\times K_{\CP^1}$ and $Z_\infty$ also vanish.
A first consequence is that $x_1\circ u$ and $x_1\circ u'$ are nowhere vanishing holomorphic
functions on $C$ and $C'$, taking values in the circle of radius $r_1$ at the boundary; this 
implies that $x_1$ is constant along $u(C)$ and $u'(C')$.

Now assume, in addition to $\mu(u)=0$, that the constant value of 
$x_1$ along $u(C)$ is not equal to 1.
Thus, $u(C)$ is disjoint from $E_\infty$, and its total intersection
number with $Z_\infty\cup E_\infty$ (the total transform of $\C^2\times
\Delta$) is zero. Since the boundary of $u(C)$ lies away from $Z_\infty\cup
E_\infty=\{x_4=0\}$, the intersection number of $u(C)$ with the levels of
$x_4$ near zero is also zero. The nonzero
levels of $x_4$ do not contain any rational curves, so positivity of
intersection implies that $u(C)$ is disjoint from those levels of $x_4$, and
hence also from $x_4=0$.  This in turn implies that $u(C)$ is disjoint from
$Z_\infty$, hence from the proper transform of $\C^2\times \Delta$ and from
the exceptional divisor $E_0$. 

We have now shown that $u(C)$ is disjoint from all components of the
anticanonical divisor $D\subset X$, except possibly
$p^{-1}(\C\times\{0\}\times K_{\CP^1})$. The vanishing of $\mu(u)$ then
implies that $u(C)$ is also disjoint from that divisor.
(Or, slightly abusing the notation introduced before Proposition
\ref{prop:Kbounds}: having shown that
$N_1(u)=N_{34}(u)=K_0(u)=K_\infty(u)=0$, we deduce from $\mu(u)=0$ that
$N_2(u)=0$ as well.)  The non-vanishing of $x_2$ in turn implies that $x_2$
is constant on $u(C)$. Arguing as in the proof of Corollary
\ref{cor:Kbounds}, we now have that $u$ is a stable disc in
$\{(x_1,x_2)\}\times K_{\CP^1}\subset X$ with boundary on a product torus
(since $F$ is $T^2$-invariant), and disjoint from all the toric divisors of
$K_{\CP^1}$. Such a disc is necessarily constant. 

Summarizing: if $\mu(u)=0$ and $r_2\neq 1$ then $x_1$ is constant along $u$,
and if moreover $u$ is not a constant disc then the value of $x_1$ along
$u(C)$ is necessarily equal to 1 (which also implies that $r_1=1$).
This completes the proof in the case where $r_2\neq 1$. The argument for the
case $r_1\neq 1$ is identical (up to exchanging the roles of $x_1$ and
$x_2$, $E_0$ and $E_\infty$, etc.).
\endproof

\begin{corollary}\label{cor:wallcrossing}
The wall-crossing coordinate transformations $\varphi_{0-}:X^\vee_{+-}\to
X^\vee_{--}$ and $\varphi_{0+}:X^\vee_{++}\to X^\vee_{-+}$ across the
walls at $r_1=1$ preserve the coordinates $z_2,z_3,z_4$.
The wall-crossing coordinate transformations $\varphi_{-0}:X^\vee_{-+}\to
X^\vee_{--}$ and $\varphi_{+0}:X^\vee_{++}\to X^\vee_{+-}$ across the
walls at $r_2=1$ preserve the coordinates $z_1,z_3,z_4$.
\end{corollary}

\proof The wall-crossing transformations $\varphi_{0-}$ and $\varphi_{0+}$
are determined by the bubbling phenomena that occur in moduli spaces of holomorphic discs 
with boundary on fibers $F_{(r_1,r_2,\xi_3,\xi_4)}$ of $\pi$ as the value of $r_1$
passes through $1$ (after regularization by a small perturbation of the
complex structure as in Proposition \ref{prop:Kbounds}). We focus our attention on
Maslov index 2 discs representing the classes
$\beta_2$, $\beta_{3,\pm}$ and $\beta_4$, whose boundary passes through a 
generic point of $F_{(r_1,r_2,\xi_3,\xi_4)}$.   The value of $x_1$ 
along each of these discs is constant, and equal to the value of the $x_1$
coordinate at the
chosen boundary point constraint. Thus, as long as the family of
point constraints we choose for varying $r_1$ avoids $x_1=1$ as the value of
$r_1$ crosses 1, it follows from Proposition \ref{prop:maslovzero} that
none of these discs can participate in any disc bubbling phenomena. (Indeed, given that all
Maslov indices are non-negative, wall-crossing for Maslov index 2 discs
only involves Maslov index 0 bubbles, but by Proposition \ref{prop:maslovzero}
those all live inside the divisor $\{x_1=1\}$.)
It follows that the portions of the superpotential which count those discs
must match under the wall-crossing transformations. 
As noted in Section
\ref{ss:superpot}, $n_{\beta_2}=n_{\beta_3,\pm}=n_{\beta_4}=1$ in all four
coordinate charts. Hence, the
terms $z_2$, $qz_3^{\pm 1}z_4$, and $z_4$ in the expressions for
$W_{\pm,\pm}$ must match under $\varphi_{0\pm}$; it follows that $\varphi_{0-}$ and
$\varphi_{0+}$ preserve each of the coordinates $z_2,z_3,z_4$.

The argument for $\varphi_{-0}$ and $\varphi_{+0}$ is identical: we
consider the contributions to the superpotential from Maslov index 2
discs representing the classes $\beta_1$, $\beta_{3,\pm}$ and $\beta_4$,
along which $x_2$ is constant, so that disc bubbling across $r_2=1$ 
can be excluded by considering a family of point constraints that avoid
$x_2=1$; this implies the invariance of $z_1,z_3,z_4$ under the
wall-crossing transformations.
\endproof

Theorem \ref{thm:main} now follows directly from the calculations of the
superpotentials $W_{\pm,\pm}$ carried out in Section \ref{ss:superpot}, the fact that
the expressions \eqref{eq:Wmain} must match under the wall-crossing coordinate transformations, and
Corollary \ref{cor:wallcrossing}.

To be more explicit, the ``basic'' {\bf stable discs of Maslov index 0} that arise along
the walls at $r_2=1$ for $\xi_4>-\xi_3-a+\epsilon'$ (i.e., away from the exceptional divisor
$E_0$) belong to three families, one of which only exists for $r_1>1$:

\begin{enumerate}
\item The proper transform of a ``standard'' disc in $\C^2\times K_{\CP^1}$
with $n_3=1$ and $n_1=n_2=n_4=0$, with $x_2=1$, and where $x_3$ has a zero
rather than a pole. These discs have Maslov index 2 in $\C^2\times
K_{\CP^1}$, but intersect the
toric divisor $\C^2\times L_0$ at a point of $H_0$, so that their lift to
$X$ is disjoint from the divisor $D$ and has Maslov index zero. (These are
the ``typical'' Maslov index 0 discs that arise in blowups of toric varieties
along codimension 2 subvarieties contained in a toric divisor; compare \cite{AAK}.)
These discs represent the class $\beta_{3,+}-[\ell_{0}]$, where $[\ell_0]$ is the
class of the fiber of $p$ above a point of $E_0$, and their Floer-theoretic
weight is $q'z_3z_4$.\smallskip
\item The union of a standard disc along the
$x_4$ coordinate axis (representing the class $\beta_4$) at $x_2=1$ and a rational curve
$S_{(x_1,1)}$. These stable discs are regular by Lemma \ref{l:z4stablediscs}(2'), and
their weight is $qq'z_4$.\smallskip
\item For $r_1>1$: the union of $S_{(1,1)}$ with a disc on which $x_2$ and $x_3$
are constant, with $x_2=1$, while $x_1$ and $x_4$ have degree 1, and $x_1=1$ at
the unique point where $x_4$ vanishes. The disc component can be
parametrized by $x_1(z)=r_1\,z$, $x_4(z)=e^{i\theta}r_4(r_1z-1)/(r_1-z)$,
and represents the class $\beta_1+\beta_4$.
These stable discs are regular by Lemma \ref{l:z1z2z4stablediscs}(1'),
and their weight is $q'q''z_1z_4$.
\end{enumerate}

There are of course other Maslov index 0 discs, representing classes which
are linear combinations (with non-negative integer coefficients) of these three,
including multiple covers as well as discs built from
unions of the above configurations. The proof of Theorem \ref{thm:main}
shows that the various Maslov index 0 discs present along the walls at
$r_2=1$ altogether amount to the
wall-crossing transformations $\varphi_{-0}$ and $\varphi_{+0}$ described by
\eqref{eq:wallcrossmain}. A similar analysis can be carried out for the walls at $r_1=1$.
\medskip

To complete our discussion, we briefly
consider the {\bf stable discs of negative Maslov index} which occur at $r_1=r_2=1$;
for simplicity we only consider the fibers of $\pi$ which lie away from the
exceptional divisors $E_0$ and $E_\infty$, and only aim to identify the
``basic'' negative Maslov index discs from which all others may be
constructed.

Assume that a stable disc $u:C\to X$ of negative Maslov
index deforms to a $J'$-holomorphic stable disc $u':C'\to X$  under arbitrarily small
deformations of the complex structure as in Proposition \ref{prop:Kbounds}.
Restricting to a subset of the components of $u$, we may assume that $C'$
has only one disc component. (When decomposing $u$ according to the 
components of $C'$, at least one of the resulting pieces must
still have negative Maslov index.) It then follows from Proposition
\ref{prop:Kbounds} that $x_1$ and $x_2$ are constant and equal to $1$ along $u(C)$.
Indeed, if $x_1\not\equiv 1$ then Proposition \ref{prop:Kbounds} gives
$K_0(u)\leq N_{34}(u)$ and $K_\infty(u)\leq N_1(u)$, so using
\eqref{eq:maslovindexsum} we conclude that $\mu(u)\geq 2N_2(u)\geq 0$; and similarly
if $x_2\not\equiv 1$ then $K_0(u)\leq N_2(u)$ and $K_\infty(u)\leq N_{34}(u)$
so that $\mu(u)\geq 2N_1(u)\geq 0$. This in turn implies that $u(C)$ is a
stable disc with boundary on a product torus in $p^{-1}(\{(1,1)\}\times
K_{\CP^1})$; we can restrict our attention to the
proper transform of $\{(1,1)\}\times K_{\CP^1}$, since sphere components inside $E_0$ or
$E_\infty$ have positive Chern number. We are thus left with a disc in
$K_{\CP^1}$, whose $x_3$ and $x_4$ components admit Blaschke product
expressions with $n_3$ and $n_3+n_4$ factors as in \eqref{eq:blaschke}, together with one or more sphere
components mapping to $S_{(1,1)}$ with total multiplicity $m$. 

The Maslov index in $X$ of such a stable disc
is $\mu(u)=2n_4-4m$. Moreover, positivity of intersection of $u'(C')$ with
the divisors $Z'_0$ and $Z'_\infty$ (and careful consideration of the local
contributions to these intersections) implies that $m\leq n_4$.%
\footnote{The intersection number of $u'(C')$ with $Z'_0$ (resp.\ $Z'_\infty$) is the number
 of poles (resp.\ zeroes) of $x_3$ plus $n_4$ minus $m$.
 Considering the local contributions to these intersection numbers
 over the regions of $C'$ which correspond to clusters of sphere components of
 $C$, non-negativity of the local intersection numbers implies that the total 
 multiplicity of the sphere components attached at any point of a disc
 component of $u(C)$ is at most the order of contact of the disc component
 with the zero section of $K_{\CP^1}$. Thus, near every point of the domain
 the local contribution to $m$ is bounded by the local contribution to
 $n_4$.}
Hence, the very simplest configuration with $\mu(u)=-2$
corresponds to the case where $n_3=0$ and $n_4=m=1$, i.e.\ the
union of a standard disc along the $x_4$ coordinate axis (representing the
class $\beta_4$) and the rational curve $S_{(1,1)}$. This configuration is
regular by Lemma \ref{l:z4stablediscs}\,(3) (in the sense described there), and its Floer-theoretic weight
is $q'q''z_4$.

The next case to consider is when $n_3>0$ and $n_4=m=1$. These
configurations arise in families that have excess dimension along the
$x_3,x_4$ factors (as the disc component is the proper
transform of a disc of higher Maslov index in $K_{\CP^1}$), but 
carry nontrivial obstruction bundles along the $x_1$ and/or
$x_2$ coordinate axes ($\mathcal{N}_{u,1}$ or $\mathcal{N}_{u,2}$ in the
terminology of Lemma \ref{l:z4stablediscs}) depending on whether $x_3$ has
poles and/or zeroes. We conjecture that these discs do not contribute to
the enumerative geometry of $X$. Specifically, it
seems that a suitable deformation of the complex structure on $X$
would ensure that the walls of Maslov index 0 discs with $n_3=1$ propagating
from the exceptional divisors $E_0$ and $E_\infty$ live at slightly
different values of $x_1$ and $x_2$ than the Maslov index $-2$ discs
with $n_4=1$ which propagate from $S_{(1,1)}$, preventing the occurrence of
configurations representing a linear combination of these classes.
More generally, we conjecture that the only stable
discs of negative Maslov index relevant to the enumerative geometry of $X$ are those we have discussed
above, representing the class $\beta_4+[S_{(1,1)}]$.

\subsection{A compact example}\label{ss:compactified}
Our main example is not very interesting from the perspective of homological
mirror symmetry, as the mirror superpotential does not have any critical
points in the geometrically relevant range of values of the coordinates
$z_i$ ($\val(z_i)\in \R_{\geq 0}^2\times \Delta$), and the wrapped Fukaya
category of $X$ is expected to be trivial. In this section we briefly
describe the analogous result for a compactified example. 

Let $\bar{X}$ be the blowup of $\CP^1\times \CP^1\times \FF_2$ at 
$\bar{H}_0=\CP^1\times\{1\}\times \bar{L}_0$ and
$\bar{H}_\infty=\{1\}\times\CP^1\times \bar{L}_\infty$, equipped with a
suitable $T^2$-invariant K\"ahler form; here
$\FF_2=\PP(\O_{\CP^1}\oplus \O_{\CP^1}(-2))$ is the second Hirzebruch
surface, and $\bar{L}_0$ and $\bar{L}_\infty$ are the fibers of the
projection from $\FF_2$ to $\CP^1$ over $0$ and $\infty$. 
The proper transform $\bar{D}$ of the toric anticanonical divisor of
$\CP^1\times\CP^1\times \FF_2$ is an anticanonical divisor in $\bar{X}$.
We construct a Lagrangian torus fibration on $\bar{X}\setminus\bar{D}$ 
with fibers
$$F_{(r_1,r_2,\xi_3,\xi_4)}=\{|x_1|=r_1,\ |x_2|=r_2,\ \mu_3=\xi_3,\
\mu_4=\xi_4\}$$
exactly as in Definition-Proposition \ref{def:fibration}, with the only
difference that $(\xi_3,\xi_4)$ now take values in the interior of the
moment polytope of $\FF_2$, i.e.\ 
$$\bar{\Delta}=\{(\xi_3,\xi_4)\in \R^2\,|\,\max(0,|\xi_3|-a)\leq \xi_4\leq
b\}.$$ Here $a$ is again half the symplectic area of the exceptional section
of $\FF_2$, and $b$ is the symplectic area of the fibers of the projection
to $\CP^1$. Let $A_1,A_2$ be the symplectic areas of the two
$\CP^1$ factors, and denote by $\epsilon'$ and $\epsilon''$ the sizes of the
blowups as previously.

The derivation of the SYZ mirror of the log Calabi-Yau pair
$(\bar{X},\bar{D})$ equipped with this Lagrangian torus fibration runs
along the same lines as the argument presented above for $(X,D)$; in
particular, it is again the case that Maslov index zero discs only
arise along walls at $r_1=1$ and $r_2=1$, and negative Maslov index discs
only arise at $r_1=r_2=1$.

\begin{proposition}
The SYZ mirror of $(\bar{X},\bar{D})$ is built out of four charts which
are domains in $(\K^*)^4$, with superpotentials
\begin{eqnarray}\label{eq:Wcompactified}
\nonumber W_{--}&=&z_1+q_1z_1^{-1}(1+qq''z_4+q''z_3^{-1}z_4)
+z_2+q_2z_2^{-1}(1+qq'z_4+q'z_3z_4)\\
\nonumber && \quad +\,q_1q_2q'q''z_1^{-1}z_2^{-1}z_4+(1+q^2+qz_3+qz_3^{-1})z_4+q_4z_4^{-1},\\
\nonumber W_{-+}&=&z_1+q_1z_1^{-1}(1+qq''z_4+q''z_3^{-1}z_4)
+z_2(1+qq'z_4+q'z_3z_4)+q_2z_2^{-1}\\
\nonumber && \quad +\,q_1q'q''z_1^{-1}z_2z_4+(1+q^2+qz_3+qz_3^{-1})z_4+q_4z_4^{-1},\\
W_{+-}&=&z_1(1+qq''z_4+q''z_3^{-1}z_4)+q_1z_1^{-1}
+z_2+q_2z_2^{-1}(1+qq'z_4+q'z_3z_4)\\
\nonumber && \quad +\,q_2q'q'' z_1z_2^{-1}z_4+(1+q^2+qz_3+qz_3^{-1})z_4+q_4z_4^{-1},\\
\nonumber W_{++}&=&z_1(1+qq''z_4+q''z_3^{-1}z_4)+q_1z_1^{-1}
+z_2(1+qq'z_4+q'z_3z_4)+q_2z_2^{-1}\\
\nonumber &&\quad +\,q'q''z_1z_2z_4+(1+q^2+qz_3+qz_3^{-1})z_4+q_4z_4^{-1},
\end{eqnarray}
where $q=T^a$, $q'=T^{a-\epsilon'}$, $q''=T^{a-\epsilon''}$,
$q_1=T^{A_1}$, $q_2=T^{A_2}$, and $q_4=T^b$. These charts are
glued pairwise by coordinate transformations which preserve $z_3,z_4$ and
act on $z_1,z_2$ by
\begin{align}\label{eq:wallcrosscompactified}
\nonumber \varphi_{-0}(z_1,z_2)&=(z_1,z_2(1+qq'z_4+q'z_3z_4+q_1q'q''z_1^{-1}z_4)),& \varphi_{-0}^*(W_{--})=W_{-+},\\
\nonumber \varphi_{+0}(z_1,z_2)&=(z_1,z_2(1+qq'z_4+q'z_3z_4+q'q''z_1z_4)),& \varphi_{+0}^*(W_{+-})=W_{++},\\
\varphi_{0-}(z_1,z_2)&=(z_1(1+qq''z_4+q''z_3^{-1}z_4+q_2q'q''z_2^{-1}z_4),z_2),& \varphi_{0-}^*(W_{--})=W_{+-},\\
\nonumber \varphi_{0+}(z_1,z_2)&=(z_1(1+qq''z_4+q''z_3^{-1}z_4+q'q''z_2z_4),z_2),& \varphi_{0+}^*(W_{-+})=W_{++}.
\end{align}
\end{proposition}

The proof is essentially identical to that of Theorem \ref{thm:main}, except
the case analysis is more tedious as the $x_1$ and $x_2$ coordinates can now
have poles as well as zeroes (as does $x_4$, though this doesn't matter
nearly as much, as the standard discs hitting the section at infinity of
$\FF_2$, with weight $q_4z_4^{-1}$, do not participate in any of the
wall-crossing). It is helpful to note, as a consistency check, that the
symmetry $x_1\leftrightarrow x_1^{-1}$ of $\bar{X}$ induces a symmetry of
the mirror, which exchanges $z_1$ and $q_1z_1^{-1}$ while swapping the
chambers with $r_1<1$ and those with $r_1>1$. Similarly, $x_2\leftrightarrow
x_2^{-1}$ induces a symmetry of the mirror which exchanges $z_2$ and
$q_2z_2^{-1}$ while swapping the chambers with $r_2<1$ and those with $r_2>1$.

\section{Deformed Landau-Ginzburg models from family Floer theory}
\label{s:counting}

\subsection{Family Floer theory}\label{ss:famfloer}

As before, we consider a Lagrangian torus fibration $\pi:X^0\to B$ on the
complement $X^0=X\setminus D$ of an anticanonical divisor $D$ in a K\"ahler
manifold $X$, whose fibers $F_b=\pi^{-1}(b)$ have vanishing Maslov class in
$X^0$. Let $B^0$ be a simply connected open subset of $B$ which is disjoint from
the critical values of $\pi$. We consider the {\em uncorrected mirror}
$$X^{\vee 0}=X^\vee_{B^0}:=\bigsqcup_{b\in B^0} H^1(F_b,U_\K),$$
with its natural analytic structure
for which the Floer-theoretic weights of disc classes $\beta\in
\pi_2(X,F_b)$ define analytic functions $z^\beta\in \O(X^{\vee 0})$; we
denote by $\pi^\vee:X^{\vee 0}\to B^0$ the natural projection map.

Fixing a base point $b_0\in B^0$ and a basis $\gamma_1,\dots,\gamma_n$ of
$H_1(F_{b_0},\Z)$ (hence of the first homology of every fiber over $B^0$), 
we can consider the Floer-theoretic weights $z_i$ ($1\leq i\leq n$) 
of cylinders with boundary on $F_{b_0}\cup F_b$, obtained by transporting a 
loop in the class $\gamma_i$ in the fibers of $\pi$ over a path connecting
$b_0$ to $b$ inside $B^0$. The coordinates $(z_1,\dots,z_n)$ allow us to
identify $X^{\vee 0}$ with a domain in $(\K^*)^n$. The functions $z^\beta$
are then Laurent monomials in $z_1,\dots,z_n$ (with exponents determined by the
coefficients of $\partial \beta$ in the basis $(\gamma_1,\dots,\gamma_n)$).

Given a subset $P$ of $B^0$, analytic functions on
$X^\vee_P=(\pi^\vee)^{-1}(P)$ are Laurent series in $z_1,\dots,z_n$ which
converge adically at all points of $P$; these are a certain completion of
the ring of Laurent polynomials $\K[z_1^{\pm 1},\dots,z_n^{\pm 1}]=\K[H_1(F_{b})]$. 
The collection of these completions as $P$ ranges over suitable subsets of
$B^0$ (e.g.\ polyhedral subsets whose faces have rational slopes with
respect to the natural affine structure of $B^0$, whose inverse images
are affinoid domains in $(\K^*)^n$) then determines a sheaf $\O_{an}=\pi^\vee_*(\O_{X^{\vee
0}})$ on $B^0$.

\begin{remark}
The main reason why we restrict ourselves to a simply connected subset of
$B$ is to be able to treat the uncorrected mirror $X^{\vee 0}$ as a
single space, rather than as a collection of local charts to be assembled in
a manner that is inconsistent (until appropriately corrected) around the
singular fibers due to the monodromy of the affine structure on $B$.
This allows us to view Floer-theoretic corrections as geometric 
deformations of a single space. Another convenient feature is that, since
the abelian groups $\pi_2(X,F_b)$ form a local system over $B\setminus
\mathrm{critval}(\pi)$, they can be transported over paths in $B^0$ to
provide distinguished isomorphisms between the groups $\pi_2(X,F_b)$ for
all $b\in B^0$; we use this repeatedly in the discussion below in order to
treat the classes of 
discs with boundary in arbitrary fibers of $\pi$ over $B^0$ as elements
of a single relative homotopy group.

However, by essence our constructions are local over (the smooth part of)
$B$, and the Floer-theoretic structures on cochains with coefficients in $\O_{an}$ we
introduce below can be defined over all of $B\setminus
\mathrm{critval}(\pi)$. If one works with the
Morse-theoretic model of family Floer theory we describe below, the
corrections to the mirror geometry naturally come out to be \v{C}ech
cochains, and it is not particularly difficult to upgrade the construction
to work over all of $B\setminus \mathrm{critval}(\pi)$ by reformulating the
output in a way that only refers to the local pieces $X^\vee_P$ rather than to
the whole of $X^{\vee 0}$.
\end{remark}

We consider Floer-theoretic operations induced by moduli spaces of
holomorphic discs with boundary on the fibers of $\pi$ on cochains on
$X^{00}=\pi^{-1}(B^0)$ with coefficients in the pullback of $\O_{an}$,
giving an $A_\infty$-deformation of the classical differential and cup-product.
There are various possible models; we describe two, of which the first one
is more intuitive but unlikely to be well-defined without further
foundational work, while the second one should be
viewed as a more realistic setup to develop the theory.
(Note in any case that our main discussion only focuses on $\m_0$ and its
properties.)

\subsubsection{Singular differential forms.}\label{sss:singdiffforms}
We denote by $C^k(X^{00},\pi^*\O_{an})$ the space of linear combinations of
differential forms of degree $j$ with coefficients in $\pi^*\O_{an}$ on 
smooth codimension $\ell$ submanifolds of $X^{00}$, for all 
$0\leq j,\ell\leq k$ such that $j+\ell=k$, i.e., the completion of
$\bigoplus_{j+\ell=k}\bigoplus_{\mathrm{codim}\, Y=\,\ell} \Omega^j(Y)\otimes
\pi^*\O_{an}$ with respect to the Novikov valuation. We regard these
cochains as an enlargement of differential forms of degree $k$ on $X^{00}$
which includes currents of integration along smooth submanifolds.

Given a nonzero class $\beta\in \pi_2(X,F_b)$ and $d\geq 0$, we denote by 
$$\Mbar_{d+1}(X^{00},\beta,J)=\bigcup_{b\in B^0} \Mbar_{d+1}(\pi^{-1}(b),\beta,J)$$
the moduli space of $J$-holomorphic stable maps from
nodal discs with $d+1$ boundary marked points $z_0,\dots,z_d$ (in order along the
boundary) to $X$, with boundary contained in some fiber of $\pi$ over a
point of $B^0$, possibly regularized by some perturbation. 
(As noted above, we use parallel transport over $B^0$
to identify the groups $\pi_2(X,F_b)$ with each other for all $b\in B^0$.)
This moduli space
carries $d+1$ evaluation maps $ev_{\beta,0},\dots,ev_{\beta,d}:\Mbar_{d+1}(X^{00},\beta,J)\to
X^{00}$ (all mapping to the same fiber of $\pi$ by construction). Assume
(rather optimistically) that $\Mbar_{d+1}(X^{00},\beta,J)$ is a smooth manifold 
with corners of dimension $2n+d-2+\mu(\beta)$, with 
\begin{equation}\label{eq:dMbar}
\partial \Mbar_{d+1}(X^{00},\beta,J)=
\bigcup_{\substack{\beta_1+\beta_2=\beta\\ d_1+d_2=d+1\\ 1\leq i\leq d_1}}
\Mbar_{d_1+1}(X^{00},\beta_1,J){\ }_{ev_{\beta_1,i}}\!\!\times_{ev_{\beta_2,0}}
\Mbar_{d_2+1}(X^{00},\beta_2,J).
\end{equation}
Assume moreover that, for given $\alpha_1,\dots,\alpha_d \in C^*(X^{00},\pi^*\O_{an})$
supported on submanifolds $Y_1,\dots,Y_d\subset X^{00}$, the evaluation map
$ev_{\beta,i}$ is transverse to $Y_i$ for $i=1,\dots,d$, the submanifolds
$ev_{\beta,i}^{-1}(Y_i)\subset\Mbar_{d+1}(X^{00},\beta,J)$ intersect transversely, and the restriction of
the evaluation map $ev_{\beta,0}$ to their intersection is a submersion onto 
a smooth submanifold of $X^{00}$ (or that a consistent
perturbation scheme can be used to achieve these properties). Then we define
$$\m_{d,\beta}(\alpha_1,\dots,\alpha_d)=(ev_{\beta,0})_*(ev_{\beta,1}^*\alpha_1\wedge\dots
\wedge ev_{\beta,d}^*\alpha_d).$$

For $\beta=0$ we set $\m_{1,0}(\alpha)=\delta\alpha$, the natural extension
to $C^k(X^{00},\pi^*\O_{an})$ of the de Rham differential (if $\alpha$ is supported
on $Y\subset X^{00}$ then $\delta\alpha=d\alpha+\alpha_{|\partial Y}$), and
$\m_{2,0}(\alpha_1,\alpha_2)=\alpha_1\wedge \alpha_2$ (as a form supported on the
intersection of the supporting submanifolds of $\alpha_1$ and $\alpha_2$, which 
are assumed to be transverse); $\m_{d,0}$ is zero for $d\neq 1,2$.

Finally, we set
\begin{equation}\label{eq:DR_md}
\m_d(\alpha_1,\dots,\alpha_d)=\sum_\beta z^\beta\,\m_{d,\beta}(\alpha_1,\dots,\alpha_d).
\end{equation}

In particular, $$\m_0=\sum_{\beta\neq 0} z^\beta\, (ev_{\beta,0})_*
1_{\Mbar_1(X^{00},\beta,J),}$$
where given our assumptions the nonzero terms correspond to currents of integration
along $ev_{\beta,0}(\Mbar_1(X^{00},\beta,J))$ when these are embedded submanifolds of $X^{00}$
(obviously an extremely restrictive setting).
Assuming the restriction of $\pi$ to each of these submanifolds
is a submersion onto a smooth submanifold of $B^0$, we can further rewrite
$\m_0$ as a sum of cochains on $B^0$ with coefficients in $\O_{an}$-valued
cochains on the fiber tori, i.e.\ elements of the bigraded complex
$\mathfrak{C}$ defined in \eqref{eq:famfloercomplex}.

It seems likely that deformation
by a suitable bounding cochain $\b\in \mathfrak{C}_{>0}$ can be used to
``smudge'' the support of $\m_0$ and turn it into a smooth differential
form, avoiding many of the pitfalls of working with currents. We will not
consider this further, and instead turn our attention to a Morse-theoretic
model whose technical foundations are easier to set up.

\subsubsection{Morse cochains and perturbed holomorphic treed discs.}\label{sss:morse}

We fix a Morse function $f$ and a Morse-Smale metric on $X^{00}$, and assume
that $\nabla f$ is transverse to the boundary of $X^{00}$ over $\partial
B^0$. (See the next section for a particularly convenient class of Morse functions for
our purposes.) We now denote by $C^k(X^{00},\pi^*\O_{an})$ the space of
linear combinations of index $k$ critical points of $f$, with coefficients
in $\O_{an}$; the coefficient of $p\in \mathrm{crit}(f)$ is typically expressed as a
sum of monomials $z^\beta$, $\beta\in \pi_2(X,F_{\pi(p)})$, and lies in a
suitable completion of $\K[H_1(F_{\pi(p)})]$ (more on this below). 
We use a family version of the construction described in
\cite[Chapter 4]{CharestWoodward} for a single Lagrangian (itself an
elaboration on the work of Cornea and Lalonde \cite{CorneaLalonde}), and define
Floer operations in terms of counts of perturbed
$J$-holomorphic treed discs. (See \cite{Hoek} for details.)

Given $\beta\in \pi_2(X,F_ b)$, $d\geq 0$, and $p_0,\dots,p_d\in 
\mathrm{crit}(f)$, we denote by 
$$\Mbar_{d+1}(p_0,p_1,\dots,p_d;\beta,J)$$
the moduli space of {\em perturbed $J$-holomorphic treed discs} with inputs at
$p_1,\dots,p_d$ and output at $p_0$, representing the class $\beta$. 
These consist of:
\begin{itemize}
\item an oriented metric ribbon tree $T$ with $d+1$ semi-infinite edges 
($d$ inputs and one output);
\item for each $d_v+1$-valent vertex $v$ of $T$, a stable (perturbed) pseudo-holomorphic map $u_v$ from
a (nodal) disc $D_v$ with $d_v+1$ boundary marked points $z_{v,0},\dots,z_{v,d_v}$
(and possibly also some interior marked points) to $X$, with boundary in the fiber of $\pi$ over some point $b_v\in
B^0$;
\item for a finite edge $e$ of $T$ connecting the output of a vertex $v$ to the
$i$-th input of a vertex $v'$, a gradient flow line $u_e$ of (a perturbation of)
$f$ connecting $u_v(z_{v,0})$ to $u_{v'}(z_{v',i})$;
\item for a semi-infinite edge of $T$ connecting the $i$-th input of the
tree to the $j$-th input of a vertex $v$ (resp.\ the output of a vertex $v$
to the output of the tree), a gradient flow line connecting the critical point 
$p_i$ to $u_v(z_{v,j})$ (resp.\ $u_v(z_{v,0})$ to $p_0$).
\end{itemize}
(As a degenerate case, for $d=1$ and $\beta=0$ the moduli space consists of
gradient flow lines of $f$ connecting two critical points $p_1$ and $p_0$.)

Recalling that the abelian groups $\pi_2(X,F_b)$ form a local system over $B\setminus
\mathrm{critval}(\pi)$, we use the identifications given by parallel 
transport along the images under $\pi$ of the gradient flow lines
$u_e$ and define the total class of a treed disc 
to be the sum of the classes of its components, $\beta=\sum_v \beta_v$, where  
$\beta_v=[u_v]\in \pi_2(X,F_{b_v})$.

Transversality can be achieved as in \cite{CharestWoodward} by considering 
domain-dependent perturbations of the complex structure and of the 
Morse function, using interior intersections with Donaldson hypersurfaces to stabilize the domain discs.
The latter point requires some adjustment compared to the case of a single
Lagrangian, as we cannot arrange for a single stabilizing divisor to be 
disjoint from all the fibers of $\pi$ simultaneously. However, for each
rational point $b\in B^0$ we can find a stabilizing divisor $D_b$ which is disjoint from
$\pi^{-1}(b)$, and hence from $\pi^{-1}(U_b)$ for some neighborhood $U_b$ of
$b$. A finite number of these neighborhoods $U_{b_i}$, $i=1,\dots,N$ suffice
to cover an arbitrarily large compact subset of $B^0$ (containing the
projections of all the critical points of $f$ and connecting Morse flow
trees). Discs with boundary in $\pi^{-1}(b)$ can thus be equipped with
several collections of marked points, coming from the intersections with
the stabilizing divisors $D_{b_i}$ for all $i$ such that $b\in U_{b_i}$.
One then needs to choose consistent domain-dependent perturbation data for discs
equipped with several collections of interior marked points, in a manner
which depends continuously on $b$ and moreover factors through
the forgetful map which erases the marked
points coming from intersections with $D_{b_i}$ whenever $b$ gets sufficiently
close to $\partial U_{b_i}$.

With this understood, we define
\begin{equation}\label{eq:Morse_md} \m_d(p_1,\dots,p_d)=
\sum_{p_0,\beta}
\bigl(\#\Mbar_{d+1}(p_0,p_1,\dots,p_d;\beta,J)\bigr)\,z^\beta\,p_0
\end{equation}
for generators of the Morse complex, where the sum ranges over critical
points $p_0$ and classes $\beta$ such that the expected dimension of
$\Mbar_{d+1}(p_0,p_1,\dots,p_d;\beta,J)$ is zero. 
We then extend the definition of $\m_d$ to general inputs in 
$C^*(X^{00},\pi^*\O_{an})$ in an $\O_{an}$-linear manner; in particular,
$$\m_d(z^{\alpha_1}p_1,\dots,z^{\alpha_d}p_d):=z^{\alpha_1+\dots+\alpha_d}\,
\m_d(p_1,\dots,p_d),$$
where as before we implicitly use parallel transport in the local system
$\{\pi_2(X,F_b)\}_{b\in B^0}$ to make sense of the sum $\alpha_1+\dots+\alpha_d$.

As in the case of a single Lagrangian, the $A_\infty$-relations follow from the fact
that the boundary of $\Mbar_{d+1}(p_0,p_1,\dots,p_d;\beta,J)$ 
consists of configurations in which a gradient flow lines breaks through
a critical point of $f$, i.e.\ pairs of perturbed $J$-holomorphic treed disks.

Because of the manner in which the Floer-theoretic weights of holomorphic 
discs are transported along Morse gradient flow lines to different fibers of
$\pi$, the total symplectic areas of the $J$-holomorphic treed discs in the moduli space 
$\Mbar_{d+1}(p_0,\dots,p_d;\beta,J)$ do not coincide with the symplectic
area of the class $\beta\in \pi_2(X,F_{\pi(p_0)})$, which determines
the valuation of each term in \eqref{eq:Morse_md}; in fact the latter
quantity does not even need to be positive in general. The 
convergence of the sum \eqref{eq:Morse_md} is therefore not 
automatic. One possible solution is to choose
the Morse function $f$ so that its gradient flow trees are guaranteed to
remain within subsets of $B^0$ that are sufficiently small for Fukaya's
trick to apply. 

Specifically, every point $b\in B^0$ admits a neighborhood $V_b$ such that
the fibers of $\pi$ over points of $V_b$ can be mapped to $F_b$ by 
diffeomorphisms $\phi_{b'\to b}$ which are $C^1$-close to identity, ensuring
that $\phi_{b'\to b}^*\omega$ tames $J$ and that the symplectic areas of a
$J$-holomorphic disc with boundary on $F_{b'}$ with respect to $\omega$ and
$\phi_{b'\to b}^*\omega$ differ by at most a bounded multiplicative factor. Thus,
given critical points $p_1,\dots,p_d\in \pi^{-1}(V_b)$, and assuming that
the gradient flow lines appearing in any
treed disc with inputs $p_1,\dots,p_d$ are guaranteed to remain within 
$\pi^{-1}(V_b)$, the symplectic area of such a treed disc and the valuation
of its contribution to $\m_d(p_1,\dots,p_d)$ differ by at most a bounded
factor; hence the sum \eqref{eq:Morse_md} converges by
the same Gromov compactness argument as in the case of a single Lagrangian.
Moreover, convergence also holds for linear combinations of critical points
in $\pi^{-1}(V_b)$ with coefficients given by Laurent series
which converge adically at every point of $V_b$.
With this understood, we cover an arbitrarily large compact subset of $B^0$
by finitely many of the neighborhoods $V_{b_i},\ i=1\dots,M$, and choose
the Morse function $f$ in such a way that the gradient flow lines appearing
in any treed disc are guaranteed to be entirely contained within a single $V_{b_i}$.

\subsubsection{Adapted Morse functions}
\label{sss:adapted}
While the above construction can be carried out for fairly general Morse
functions (with the restrictions noted), the connection to family Floer
theory becomes clearer for specific classes of Morse functions,
constructed as follows. 

Start from a simplicial
decomposition $\mathcal{P}$ of a large compact subset onto which $B^0$ retracts, 
with every cell of $\mathcal{P}$ contained in a single open subset
$V_{b_i}$. Pick a Morse function $h:B^0\to \R$ and a Morse-Smale metric on $B^0$, such
that for every $k$-cell $\sigma\in \mathcal{P}^{[k]}$ the function $h$ has
a unique critical point $b_\sigma$ in the interior of $\sigma$, of index $k$, whose
descending manifold is $\sigma$ itself. (Such a function and metric
can be constructed e.g.\ from a barycentric subdivision of $\mathcal{P}$.)
Then construct the Morse function $f:X^{00}\to\R$ by combining the pullback of $h$
under the projection $\pi$ with Morse functions on the fibers of $\pi$, as
well as a Morse-Smale metric on $X^{00}$, in such a way that:
\begin{itemize}
\item all the critical points of $f$ project to critical points of $h$;
\item for each cell $\sigma\in \mathcal{P}^{[k]}$, the restriction of $f$
to $\pi^{-1}(b_\sigma)$ is a standard Morse function on the $n$-torus (i.e., it has
$2^n$ critical points, whose ascending and descending submanifolds represent
dual standard bases of $H_*(T^n)$), and every index $j$ critical point
of $f_{|\pi^{-1}(b_\sigma)}$ is also a critical point of $f$, of index
$k+j$; 
\item for each cell $\sigma$, the gradient flow of $f$ is tangent to 
$\pi^{-1}(b_\sigma)$, and the union of the descending submanifolds of the
critical points of $f$ which lie in $\pi^{-1}(b_\sigma)$ is $\pi^{-1}(\sigma)$.
\end{itemize}

\begin{definition} \label{def:adaptedmorse} We call a Morse function $f:X^{00}\to\R$ with these
properties {\em adapted} to the simplicial decomposition $\mathcal{P}$.
\end{definition}

(The assumption that $f_{|\pi^{-1}(b_\sigma)}$ has only $2^n$ critical
points and vanishing Morse differential is extraneous and might be best
left out of the definition, but it is convenient for the
rest of our discussion.)

The main advantage of adapted Morse functions for our purposes is that
Morse cochains can be expressed as Morse cochains for the function
$h$ on $B^0$ with coefficients in the Morse complexes of the functions
$f_{|\pi^{-1}(b_\sigma)}$. In this sense, for adapted $f$ we have
$$C^*(X^{00},\pi^*\O_{an})=C^*(B^0;C^*(F_b)\,\hat\otimes\,\O_{an});$$
denoting $\CC^{i,j}=C^i(B^0;C^j(F_b)\,\hat\otimes\,\O_{an})$, this
recovers the setting considered in \eqref{eq:famfloercomplex}.
Moreover, the assumption made on the restrictions of $f$ to the critical fibers
implies that the fiberwise Morse differential vanishes, so in fact we have
$$\CC^{i,j}=C^i(B^0;H^j(F_b)\,\hat\otimes\,\O_{an}).$$
By construction the Morse differential $\delta$ and the Floer differential
$\m_1$ on this complex are
filtered, in the sense that the Morse index $i$ on $B^0$ is non-decreasing;
and the only terms which preserve $i$ are the (trivial) Morse differential 
and the Floer differential on $C^*(F_{b_\sigma})=H^*(F_{b_\sigma})$ for each
$\sigma$. Meanwhile, the terms which increase $i$ by one correspond to
Morse, resp.\ Floer-theoretic continuation maps from $C^*(F_{b_\sigma})$ to
$C^*(F_{b_{\sigma'}})$ over a gradient flow line of $h$ from $b_\sigma$ to
$b_{\sigma'}$; and those which increase $i$ by more than one correspond to
homotopies between different compositions of such continuation maps.

\begin{remark}\label{rmk:morse_to_cech}
It is typically possible to arrange for the
latter homotopies to vanish in Morse theory (e.g., since we have assumed
$B^0$ to be simply connected and one also typically has $\pi_2(B^0)=0$, 
by trivializing $\pi$ over $B^0$ and taking $f$ to be the sum of the
pullback of $h$ and a fixed Morse function on $T^n$). The Morse
complex $(\CC,\delta)$ is then identified with the \v{C}ech complex 
$\check{C}^*(B^0;H^*(F_b)\,\hat\otimes\,\O_{an})$ for the polyhedral cover of (a
retract of) $B^0$ given by
the stars of the vertices of $\mathcal{P}$. We will use this fact below to
recast $\m_0$ as a \v{C}ech cochain (with values in polyvector fields) on the uncorrected mirror $X^{\vee 0}$.
\end{remark}

We finish this section by noting the manner in which the Floer-theoretic obstruction
$\m_0\in \CC$ encodes information not only about the holomorphic
discs bounded by individual fibers of $\pi$ but also about those bounded by families of fibers
over the simplices of $\mathcal{P}$. Namely, the part of $\m_0$ which lies 
over an index $0$ critical point of $h$ at a vertex of $\mathcal{P}$ counts 
(treed) holomorphic discs bounded by the fibers of $\pi$ over that point,
in the sense of Floer theory for a single Lagrangian; whereas the portion of
$\m_0$ which lives over an index $i$ critical point $b_\sigma$ of $h$
corresponds to (treed) family counts of holomorphic discs bounded by the fibers of
$\pi$ over the $i$-dimensional cell $\sigma\subset B^0$. 

In this sense, the
component of $\m_0$ in $\oplus_j \CC^{i,j}$ counts families of holomorphic discs that
occur along
(possibly thickened) codimension $i$ walls in $B^0$, i.e.\ those which can be
meaningfully counted along $i$-dimensional families of fibers of $\pi$.
The notion of weak family unobstructedness (Definition \ref{def:wfunobs}) 
expresses the requirement that all non-zero counts should live in fiberwise
cohomological degree $j=i$, i.e.\ correspond to discs of Maslov index $2-2i$.

\begin{remark}
Besides fleshing out the details of the construction of 
the curved $A_\infty$-algebra $\CC$ via perturbed $J$-holomorphic treed
discs, Hoek's thesis \cite{Hoek} also implements a key step of the family Floer program in this setting
by constructing a functor from the Fukaya category of Lagrangian sections of the
fibration $\pi$ to the category of $A_\infty$-modules over $\CC$.

\end{remark}

\subsection{A heuristic derivation of the master equation}\label{ss:mastereq}

The algebraic properties of $\m_0$ generally follow from the fact that the boundary
strata of moduli spaces of holomorphic discs are fibered products of 
moduli spaces of discs, as expressed in \eqref{eq:dMbar}. Most immediately,
this yields the identity $\m_1(\m_0)=0$, which is part of the
$A_\infty$-equations. Our goal, however, is to find (when possible)
a constraint involving only $\m_0$: the master equation
\eqref{eq:mastereq_m0}. 

In this section we give a heuristic derivation
of this equation under the assumption that the moduli spaces of holomorphic
discs entering into the definition of $\m_0$ are {\em fiberwise closed},
in order to provide motivation for Conjecture \ref{conj:mastereq}. (It seems
difficult, or in any case well beyond the scope of this paper, to make the
argument rigorous under realistic assumptions.)

One particularly convenient way to understand the origin of the master equation in Lagrangian
Floer theory
is at the level of loop spaces, as first proposed by Fukaya \cite{Floop}, and
further studied by Irie \cite{Irie}, even though the technical details
are daunting. (Working in families however does not bring much
additional complexity.) A very informal account is as follows.
The moduli space $\mathcal{M}_1(X^{00},\beta,J)$ carries
an evaluation map not only to $X^{00}$, but also to its free loop space
$\mathcal{L}X^{00}$ (in fact, to free loops contained in the fibers of
$\pi$). (This requires preferred parametrizations of the 
boundary loops, which can be done e.g.\ by stabilizing the domains or
by using arc length in $X^{00}$). 
Denote by $\m_{0,\beta}^\mathcal{L}\in
C_{2n-2+\mu(\beta)}(\mathcal{L}X^{00})$ the
evaluation pushforward of the fundamental chain of 
$\mathcal{M}_1(X^{00},\beta,J)$ 
(after a suitable regularization).
Summing over relative classes, we set $\m_0^\mathcal{L}=\sum_{\beta} \m_{0,\beta}^\mathcal{L}\,z^\beta\in
C_*(\mathcal{L}X^{00};\pi^*\O_{an})$. By analogy with \cite{Floop,Irie}, 
one expects that (up to sign)
\begin{equation}\label{eq:mastereq_loop}
\partial \m_0^\mathcal{L}=\frac12
\{\m_0^\mathcal{L},\m_0^\mathcal{L}\},
\end{equation}
where $\{\cdot,\cdot\}$ denotes a chain-level refinement of the Chas-Sullivan 
bracket on $H_*(\mathcal{L}X^{00})$; or rather, as shown by Irie, the chain-level master equation
also involves higher order terms due to the chain-level loop bracket actually
being only part of a homotopy Lie ($L_\infty$) structure on chains on the
loop space \cite{Irie}.
To avoid the inherent difficulties of chain-level string topology, we focus
on the main case of interest to us, and assume that the moduli spaces of
holomorphic discs under consideration are fiberwise closed manifolds. 
To further avoid the need to regularize the moduli spaces, we 
make the following (unrealistic) assumptions about $\Mbar_1(X^{00},\beta,J)$:
\begin{itemize}
\item (regularity) $\Mbar_1(X^{00},\beta,J)$ is a smooth manifold with corners, of the expected
dimension, whose boundary is as in \eqref{eq:dMbar};
\item (transversality) the projections $\pi_*:\Mbar_1(X^{00},\beta,J)\to B^0$
induced by $\pi$ are submersions onto
smooth submanifolds of $B^0$ with boundary and corners which meet
transversely;
\item (fiberwise closed) the fibers of $\pi_*:\Mbar_1(X^{00},\beta,J)\to B^0$
are closed manifolds, i.e.\ 
\begin{equation}\label{eq:Mbar1fiberwiseclosed}
\pi_*(\partial \Mbar_1(X^{00},\beta,J))\subset \partial(\pi_*(
\Mbar_1(X^{00},\beta,J)).\end{equation}
\end{itemize}
Then we can view $\m_0^\mathcal{L}$ as a chain on $B^0$ with coefficients in
$H_*(\mathcal{L}F_b)\,\hat\otimes\,\O_{an}$, and the master equation
\eqref{eq:mastereq_loop} expresses the boundary of $\m_0^\mathcal{L}$
(as a chain on $B^0$) in terms of the bracket induced by 
the classical cup-product on $B^0$ and the Chas-Sullivan loop bracket
\cite[Definition 4.1]{ChasSullivan} on
$H_*(\mathcal{L}F_b)$. 

\begin{lemma}
Under these assumptions, \eqref{eq:dMbar} implies that
$\m_0^\mathcal{L}\in C_*(B^0;H_*(\mathcal{L}F_b)\,\hat\otimes\,\O_{an})$
satisfies the master equation \eqref{eq:mastereq_loop}.
\end{lemma}

\proof[Sketch of proof] On one hand,
$\partial \m_{0,\beta}^\mathcal{L}$ is the image of the boundary of
$\Mbar_1(X^{00},\beta,J)$ under the loop space-valued evaluation map.
On the other hand, the evaluation image of $\Mbar_2(X^{00},\beta_1,J)
{\ }_{ev_{\beta_1,1}}\!\!\times_{ev_{\beta_2,0}}
\Mbar_{1}(X^{00},\beta_2,J)$ is 
the chain formed by inserting the loops that appear in
$\m_{0,\beta_2}^\mathcal{L}$ into the loops that make up
$\m_{0,\beta_1}^\mathcal{L}$ whenever the latter pass through the base 
points of the former, i.e.\ $\m_{0,\beta_2}^\mathcal{L}*
\m_{0,\beta_1}^\mathcal{L}$ in the notation of \cite[\S 3]{ChasSullivan}.
Summing over all $\beta_1,\beta_2$ such that $\beta_1+\beta_2=\beta$, 
we find that the evaluation image of the right-hand side of \eqref{eq:dMbar}
is equal to the coefficient of $z^\beta$ in $\m_0^\mathcal{L}*\m_0^\mathcal{L}=\frac12
\{\m_0^\mathcal{L},\m_0^\mathcal{L}\}$.
\endproof

Now we observe that each term $\m_{0,\beta}^\mathcal{L}$ consists of loops
representing the class $\partial\beta\in H_1(F_b)$, and recall that
each component of $\mathcal{L}F_b$ is homotopy equivalent to $F_b$ itself,
via evaluation at the base point. A simple calculation shows:

\begin{lemma}\label{l:bracket-loop-to-HF}
Denoting by $\mathcal{L}_\gamma
T^n$ the component of $\mathcal{L}T^n$ which consists of loops in
the class $\gamma\in H_1(T^n)$, and using evaluation at the base point and
Poincar\'e duality to identify
$H_*(\mathcal{L}_\gamma T^n)$ with $H_*(T^n)\simeq H^{n-*}(T^n)\simeq
\bigwedge^{n-*}H^1(T^n)$, up to sign
the Chas-Sullivan bracket
$\{\cdot,\cdot\}:H_*(\mathcal{L}_\gamma T^n)\otimes
H_*(\mathcal{L}_{\gamma'} T^n)\to H_*(\mathcal{L}_{\gamma+\gamma'}T^n)$
is given
by $$\{\alpha,\alpha'\}=\alpha\wedge (\iota_\gamma \alpha')+(-1)^{|\alpha|}
(\iota_{\gamma'}\alpha)\wedge \alpha'.$$
\end{lemma}
\proof
We can represent the classes $\alpha,\alpha'$ by cycles consisting of
straight line loops on a flat torus, with tangent vectors given by
$\gamma$ and $\gamma'$ respectively (under the identification of the first
homology of a flat torus with the lattice of integer tangent vectors).
The element $\iota_\gamma \alpha'$ is (Poincar\'e dual to) the cycle on
$T^n$ obtained by 
spreading the evaluation image of $\alpha'$ by translation along $\gamma$,
i.e.\ the set of points $p\in T^n$ such that the straight loop in the direction of
$\gamma$ based at $p$ hits the base point of one of the loops in the
cycle $\alpha'$. Thus, $\alpha\wedge (\iota_\gamma\alpha')$ corresponds to the
cycle formed by the base points of loops in the chain $\alpha$
which hit the base points of the loops in the chain $\alpha'$.
This is, up to sign, the operation denoted $\alpha'*\alpha$ in
\cite[Section~3]{ChasSullivan}, whose skew-symmetrization is the Chas-Sullivan bracket
\cite[Definition 4.1]{ChasSullivan}.
\endproof

\noindent
Using Lemma \ref{l:bracket-loop-to-HF} to rewrite the Chas-Sullivan bracket
$\{\m_0^\mathcal{L},\m_0^\mathcal{L}\}$ in terms of the bracket 
defined by \eqref{eq:HF-bracket} on $H^*(T^n)\otimes \K[H_1(T^n)]$, 
we arrive at:

\begin{corollary}
Still assuming regularity of moduli spaces and
transversality of evaluation maps, if $\m_0^\mathcal{L}\in C_*(B^0;H_*(\mathcal{L}F_b)\,\hat\otimes\,\O_{an})$
satisfies \eqref{eq:mastereq_loop} then $\m_0\in
C^*(B^0;H^*(F_b)\,\hat\otimes\,\O_{an})$ satisfies \eqref{eq:mastereq_m0}.
\end{corollary}

We can in fact give a
more direct derivation of \eqref{eq:mastereq_m0} without involving loop
spaces:

\begin{proposition}\label{prop:fiberwiseclosed_mastereq}
Assuming that family Floer theory can be set up using the singular 
differential forms model of Section \ref{sss:singdiffforms} and
that the moduli spaces of holomorphic discs are fiberwise closed in
the sense of \eqref{eq:Mbar1fiberwiseclosed},
the cochain $\m_0\in
C^*(B^0;H^*(F_b)\,\hat\otimes\,\O_{an})$ satisfies \eqref{eq:mastereq_m0}
up to sign.
\end{proposition}

\proof
Fixing a family of flat metrics on the fibers of $\pi$ over $B^0$, we
can deform by simultaneous homotopies the boundary loops of all holomorphic stable 
discs in $\Mbar_1(F_b,\beta,J)$ (parametrized e.g.\ by arc length)
into straight line geodesics representing the class $[\partial\beta]\in
H_1(F_b)$, for all $b\in B^0$. This produces a homotopy between the
chains represented by the evaluation map
$$(ev_{\beta,0},ev_{\beta,1}):\Mbar_2(X^{00},\beta,J)\to X^{00}\times
X^{00}$$ and by 
$$(ev_{\beta,0},t_{[\partial\beta]}\circ ev_{\beta,0}):
\Mbar_1(X^{00},\beta,J)\times S^1\to X^{00}\times X^{00},$$
where $t_{[\partial\beta]}:X^{00}\times S^1\to X^{00}$ denotes translation along the straight line
geodesics in the class $[\partial\beta]$ inside the fibers of $\pi$.
This implies that
$$(ev_{\beta_1,0})_*\left[\Mbar_2(X^{00},\beta_1,J)
{\ }_{ev_{\beta_1,1}}\!\!\times_{ev_{\beta_2,0}}
\Mbar_{1}(X^{00},\beta_2,J)\right]$$
and $$(ev_{\beta_1,0})_*\left[\Mbar_1(X^{00},\beta_1,J)\times S^1
{\ }_{t_{[\partial\beta_1]}\circ ev_{\beta_1,0}}\!\!\times_{ev_{\beta_2,0}}
\Mbar_{1}(X^{00},\beta_2,J)\right]$$
are equal as cochains on $B^0$ with coefficients in $H^*(F_b)$. Since the
fiber product expresses the condition that the output marked points of the
two discs line up along a straight line geodesic in the class
$[\partial\beta_1]$, the latter cochain can be expressed as
\begin{multline*}
(ev_{\beta_1,0})_*\left[\Mbar_1(X^{00},\beta_1,J)
{\ }_{ev_{\beta_1,0}}\!\!\times_{t_{[\partial\beta_1]}\circ ev_{\beta_2,0}}
(\Mbar_{1}(X^{00},\beta_2,J)\times S^1)\right]\\
=(ev_{\beta_1,0})_*\left[\Mbar_1(X^{00},\beta_1,J)\right] \cap
t_{[\partial\beta_1]}\left((ev_{\beta_2,0})_*[\Mbar_1(X^{00},\beta_2,J)]\times
S^1\right).
\end{multline*}
Since spreading a homology class along $[\partial\beta_1]$ corresponds under
Poincar\'e duality to interior product with
$[\partial\beta_1]$, this expression can be rewritten more concisely as
$$\m_{0,\beta_1}\wedge \iota_{[\partial\beta_1]}(\m_{0,\beta_2}).$$
It then follows from \eqref{eq:dMbar} that
$z^\beta\,\delta \m_{0,\beta}$ is, up to sign, equal to
$$\sum_{\beta_1+\beta_2=\beta}
z^{\beta_1+\beta_2}\,\m_{0,\beta_1}\wedge
\iota_{[\partial\beta_1]}\m_{0,\beta_2}=\frac12
\sum_{\beta_1+\beta_2=\beta}\,\{z^{\beta_1}\m_{0,\beta_1},z^{\beta_2}\m_{0,\beta_2}\}.$$
Summing over $\beta$ then gives \eqref{eq:mastereq_m0}.
\endproof

\subsection{Spliced treed $J$-holomorphic discs and the master equation}
\label{ss:boosted}

With some care, it seems likely that the argument of Proposition \ref{prop:fiberwiseclosed_mastereq} can be
transcribed into the language of Morse cochains and perturbed holomorphic
treed discs, to arrive at a similar result in that setup, still subject to
very strong assumptions about moduli spaces of discs. However, it is more
appealing to try to modify the model of Section \ref{sss:morse} to arrive at
a setup where the master equation holds in full generality.
In this section, we sketch such an approach. 

Despite the fairly detailed 
outline, the description we give here
is by no means complete: we skip over various limiting cases, and do not
attempt to check consistency, discuss orientations, or prove the existence of suitable perturbation data.
The details of the construction will appear elsewhere.

\begin{remark}
The approach we describe here using ``standard loops'' has some advantages
but also some notable drawbacks, chief among them the need to choose and
keep track of a number of homotopies between various types of loops. As of
this writing it is likely that the construction will eventually be modified
to rely on a suitable geometric flow and evolve families of fiberwise loops along 
the tree portions of spliced treed discs, rather than homotoping them to standard
loops.
\end{remark}

The boundary of the usual moduli space of treed holomorphic discs
$\Mbar_1(p_0;\beta,J)$ consists of configurations where the length of an
internal edge becomes infinite; these can be viewed as pairs of treed discs
where the output of one treed disc serves as an {\em input} for the other, 
giving rise to the identity $\m_1(\m_0)=0$. Our aim is to modify the moduli 
space so that its boundary consists of pairs
of configurations whose {\em outputs} are matched to each other 
via the bracket $\{\cdot,\cdot\}$ defined in
\eqref{eq:HF-bracket}. We do this by allowing the matching 
condition at the ends of broken (infinite length) gradient flow lines 
to deform towards the output $p_0$. More precisely, once the length of a
gradient flow line in a treed disc becomes infinite (i.e., the flow line
breaks through a Morse critical point), we first allow the incidence condition
for the end point of the flow line to deform along a homotopy from the boundary loop of the
appropriate disc component to a ``standard'' loop
in the same homotopy class, and then we allow the standard loop to slide along
the gradient flow tree towards the output of the treed disc. (Standard loops, defined below, are a class of loops
which are well-behaved with respect to the action of $H_1(F_b)$ on the Morse
cohomology of $f$ by interior product.) For simplicity, we assume that we work with an
adapted Morse function for some simplicial decomposition $\mathcal{P}$ of $B^0$ in the 
sense of Definition \ref{def:adaptedmorse}.

\subsubsection{Standard loops}

By Definition \ref{def:adaptedmorse}, the critical points of an
adapted Morse function $f:X^{00}\to\R$ lie in fibers
$F_{b_\sigma}=\pi^{-1}(b_\sigma)$ indexed by the cells $\sigma\in
\mathcal{P}^{[k]}$ of $\mathcal{P}$, and the restriction of $f$ to
$F_{b_\sigma}$ is a standard Morse function on $T^n$, i.e.\ 
there is a basis $e_{\sigma,1},\dots,e_{\sigma,n}$ of
$H_1(F_{b_\sigma})$ such that the ascending and descending submanifolds of the
critical points of $f_{|F_{b_\sigma}}$ represent exterior products of 
elements of the basis. For $I\subseteq \{1,\dots,n\}$ we denote by $p_{\sigma,I}$ the critical point
of index $|I|$ whose descending (resp.\ ascending) submanifold within
$F_{b_\sigma}$ represents
the homology class $e_{\sigma,I}=\bigwedge_{i\in I} e_{\sigma,i}$ (resp.\ 
$e_{\sigma,\overline{I}}$, where $\overline{I}=\{1,\dots,n\}-I$).

Given a homology class $[\gamma]=\sum n_i
e_{\sigma,i}\in H_1(F_{b_\sigma},\Z)$, interior product with $[\gamma]$ defines
an operator of degree $-1$ on $H^*(F_{b,\sigma},\Z)\simeq CM^*(f_{|F_{b_\sigma}},\Z)$, 
$$\iota_{[\gamma]}:CM^*(f_{|F_{b_\sigma}})\to CM^{*-1}(f_{|F_{b_\sigma}}),$$
which maps $p_{\sigma,I}$ to $\iota_{[\gamma]}(p_{\sigma,I})=\sum_{i\in I}
(-1)^{|I\cap\{1,\dots,i-1\}|}\,n_i\, p_{\sigma,I-\{i\}}$.
For $p_{\sigma,I}\in \mathrm{crit}(f)$, we denote by 
$[\overline{W}^+(p_{\sigma,I})]$ the fundamental chain of the (closure of the) ascending manifold of
$p_{\sigma,I}$ inside $X^{00}$, and for $[\gamma]\in H_1(F_{b_\sigma})$
we define $[\overline{W}^+(\iota_{[\gamma]}(p_{\sigma,I}))]$
to be the appropriate linear combination of the ascending manifolds
$\overline{W}^+(p_{\sigma,I-\{i\}})$, $i\in I$.

\begin{definition}\label{def:stdloops}
Let $\widetilde{X}^{00}$ be the 
space of pairs $([\gamma],x)$ where $x\in X^{00}$ and $[\gamma]\in
H_1(F_{\pi(x)},\Z)$.
A {\em system of standard loops} for $f$ is a smooth submersive map from
$\widetilde{X}^{00}\times S^1$ to $X^{00}$, $([\gamma],x,t)\mapsto
s_{[\gamma]}(x,t)=s_{[\gamma],x}(t)$, such that:
\begin{enumerate}
\item for all $[\gamma]$ and $x$, $s_{[\gamma],x}:S^1\to X^{00}$ is a loop in $F_{\pi(x)}$ based at $x$, 
representing  the homology class $[\gamma]$;
\smallskip
\item for every critical point $p_{\sigma,I}$ of $f$, and for every
$[\gamma]\in H_1(F_{b_\sigma})$, 
\begin{equation}\label{eq:stdloops}
p_*s_{[\gamma]}^{-1}([\overline{W}^+(p_{\sigma,I})])=
[\overline{W}^+(\iota_{[\gamma]}(p_{\sigma,I}))]
\end{equation}
as chains in $X^{00}$ modulo degenerate chains supported on the
lower-dimensional submanifold
$\overline{W}^+(p_{\sigma,I})$; here $p:X^{00}\times S^1\to X^{00}$ is the
projection to the first factor, and we implicitly identify $H_1(F_b)\simeq
H_1(F_{b_\sigma})$ for $b$ near $b_\sigma$.
\end{enumerate}
\end{definition}

For $[\gamma]=e_{\sigma,i}$, $i\in I$, condition (2) states that the loop
$s_{[\gamma],x}$ passes through the ascending manifold of $p_{\sigma,I}$
if and only if $x$ lies in the ascending manifold of $p_{\sigma,I-\{i\}}$, 
and in that case it does so just once (counting with appropriate signs). Likewise
for general $[\gamma]$ and linear combinations of ascending manifolds
of the critical points appearing in $\iota_{[\gamma]}(p_{\sigma,I})$.

As will be clear from the arguments below, it would in fact suffice for
the two sides of \eqref{eq:stdloops}
to be equivalent from the perspective of Morse theory, i.e.\ that they 
intersect in the same manner with the ascending and descending submanifolds 
of other critical points of $f$.

\begin{lemma} When $B^0$ is simply connected, there exists an adapted
Morse function $f$ which admits a system of standard loops.
\end{lemma}

\proof Since $\pi_1(B^0)=1$, the structure group
of the fibration $\pi:X^{00}\to B^0$ reduces to translations of the
$n$-torus $T^n=(S^1)^n$, i.e.\ we have well-defined fiberwise coordinates 
up to translation on the fibers of $\pi$. We can then choose the admissible
Morse function $f$ so that its restriction to each fiber of $\pi$ is 
the sum of standard Morse functions on the $S^1$ factors and a Morse
function on the base $B^0$, and choose the 
metric in a suitable manner, so that the ascending manifold 
$\overline{W}^+(p_{\sigma,I})$ is invariant under translation along the $i$-th
$S^1$ factor in the fibers of $\pi$ whenever $i\not\in I$, and 
translating it along the $i$-th $S^1$ factor for $i\in I$ yields exactly
$\overline{W}^+(p_{\sigma,I-\{i\}})$.

Denoting by $e_i$ the homology class of the $i$-th $S^1$ factor, and
given a class $[\gamma]=\sum n_ie_i$ and a point $x\in X^{00}$, 
we define the loop $s_{[\gamma],x}$ to be the concatenation of loops
based at $x$ which run successively $n_i$ times along each $S^1$ factor 
of $F_{\pi(x)}$ (with vanishing derivatives at the end points so that the
concatenation is a smooth loop). The identity \eqref{eq:stdloops} then
follows from the observation that a loop based at $x$ and running
along the $i$-th $S^1$ factor intersects $\overline{W}^+(p_{\sigma,I})$
if and only if $x$ lies in the image of $\overline{W}^+(p_{\sigma,I})$ under
translation along the $i$-th $S^1$ factor, i.e.\ $\overline{W}^+(p_{\sigma,I-\{i\}})$
if $i\in I$ and $\overline{W}^+(p_{\sigma,I})$ itself otherwise.
\endproof

\subsubsection{Spliced treed discs}\label{sss:spliced}

A {\em spliced treed disc} consists of a collection of $k+1$ treed
discs $\mathbb{T}_\alpha=(T_\alpha,\{D_v\}_{v\in \mathrm{Vert}(T_\alpha)})$,
$\alpha\in \{0,\dots,k\}$, inductively attached onto each other by $k$ 
semi-infinite edges $e^{spl}_\alpha$, $\alpha\in\{1,\dots,k\}$ (the {\em splicings}).
Each splicing $e^{spl}_\alpha$ connects the output
of the treed disc $\mathbb{T}_\alpha$ to some point 
$t^{spl}_\alpha$ (the {\em target} of the splicing) in $T^+_{<\alpha}:=\bigcup_{\alpha'<\alpha} T_{\alpha'}
\cup \bigcup_{\alpha'<\alpha} e^{spl}_{\alpha'}$, the underlying tree
of the configuration obtained by splicing the treed discs
$\mathbb{T}_{\alpha'}$ for $\alpha'<\alpha$.

The end result of this process differs from a stable treed disc in that
there are some broken (infinite length) internal edges, formed by the output
edges of the tree discs $\mathbb{T}_\alpha$ together with the splicing edges
$e^{spl}_\alpha$, and these broken edges do not attach to the boundary
of a disc, but rather onto the underlying tree $$T^+=\bigcup_\alpha
T_\alpha\,\cup\,\bigcup_\alpha e^{spl}_\alpha$$ of the spliced treed disc;
the manner in which this translates into an incidence condition for the
end point of a gradient flow line is governed by {\em splicing data} which
we describe below. (Note that $T^+$ is not a ribbon tree, as the
splicing data does not specify how the splicing edge fits into a
cyclic ordering at its target.)

The splicing data for a given splicing depends on whether its target 
lies on an edge of $T^+$ or at a vertex, and on the number of splicings which share
the same target $t=t_\alpha^{spl}$. The set of incidence conditions we impose
on the ends of the splicings with target $t$ is parametrized by a manifold with corners
$S^{spl}(t)$; when all the splicings have distinct targets,
$S^{spl}(t)$ is 
$S^1$ if $t$ lies on an edge of $T^+$, or $[0,1]\times S^1$ if 
it lies at a vertex. 
The incidence conditions we impose on the
ends of the splicing edges with target $t$ are described by maps
$\sigma_\alpha:S^{spl}(t)\to X^{00}$ for all $\alpha$ such that
$t^{spl}_\alpha=t$, defined below.

\begin{definition}
A {\em spliced treed $J$-holomorphic disc} $u:\mathbb{T}^+\to X$ with domain
$\mathbb{T}^+=\bigcup \mathbb{T}_\alpha\cup \bigcup e^{spl}_\alpha$ 
consists of:
\begin{itemize}
\item for each $\alpha$, a (perturbed) treed $J$-holomorphic disc
$u_\alpha:\mathbb{T}_\alpha\to X$, i.e., (perturbed) stable
$J$-holomorphic discs $u_v:D_v\to X$ with boundary in some fiber $F_{b_v}$ of $\pi$ 
for every vertex $v$ of $T_\alpha$, connected to each
other and to the output critical point $p_{0,\alpha}\in
\mathrm{crit}(f)$ by (perturbed) gradient
flow lines $u_e$ of $f$ for every edge $e$ of $T_\alpha$;
\smallskip
\item for each point $t$ of $T^+$ which is the target of one or more
splicings, a choice of splicing data $\theta_t\in S^{spl}(t)$;
\smallskip
\item for each splicing edge $e^{spl}_\alpha$, 
a semi-infinite (perturbed) gradient
flow line of $f$ whose negative end converges to the critical point
$p_{0,\alpha}$, and whose positive end maps to the point
$\sigma_\alpha(\theta_{t^{spl}_\alpha})\in X^{00}$.
\end{itemize}
\end{definition}

The description of the maps $\sigma_\alpha:S^{spl}(t)\to X^{00}$ involves
the standard loops introduced in the previous section, as well as the following
definition:

\begin{definition}
The {\em weight} $\beta_{t\in e}$ of an edge $e$ of the tree $T^+$ underlying a spliced
treed holomorphic disc $u:\mathbb{T}^+\to X$ at a point $t\in e$ is the
sum of the homotopy classes $\beta_v=[u_v]$ of all the disc components of $u$ 
which correspond to vertices $v\in \bigcup \mathrm{Vert}(T_\alpha)$ such that the
path in $T^+$ from $v$ to the output of $T^+$ passes through $t$, and reaches
$t$ via the edge $e$.
\end{definition}

As in \S \ref{sss:morse}, we use the identifications between 
the abelian groups $\pi_2(X,F_b)$, $b\in B^0$ along the images under $\pi$ of
the gradient flow lines $u_e$ to define the sum of the homotopy
classes $\beta_v$ and view the weight $\beta_{t\in e}$ as an element of $\pi_2(X,F_b)$ for
$b=\pi(u_e(t))$. We also introduce the homology class
\begin{equation}\label{eq:boundaryweightedge}
[\gamma_{t\in e}]=\partial\beta_{t\in e}\in H_1(F_b).
\end{equation}
With this understood, let $e^{spl}_\alpha$ be a splicing edge, with target $t=t^{spl}_\alpha \in T^+$.
\smallskip

{\em Case 1.} Assume $t$ is an interior point of an edge $e$ of $T^+$,
mapping to $x=u_e(t)\in X^{00}$, and
no other splicing has the same target $t$. Then we require the end point of
the splicing to lie on $s_{[\gamma_{t\in e}],x}$, the standard loop at $x$ 
in the homology class $[\gamma_{t\in e}]=\partial\beta_{t\in e}$. Namely, we set
$S^{spl}(t)=S^1$, and require the end point of $e^{spl}_\alpha$ to map
to $\sigma_\alpha(\theta_t):=s_{[\gamma_{t\in e}], x}(\theta_t)$.
\smallskip

{\em Case 2.} Assume $t$ is a vertex $v$ of one of the trees $T_\alpha$,
corresponding to a stable $J$-holomorphic disc $u_v:D_v\to X$ with boundary
in $F_{b_v}$, and no other
splicing has the same target. Denote by $e_i$ the edges of $T_\alpha$
that attach to the input boundary marked points $z_{v,i}\in \partial D_v$,
by $\beta_{v,i}=\beta_{v\in e_i}$ the weights of these edges at their
end points, and by $x_i=u_v(z_{v,i})$ the end points of the
gradient flow lines $u(e_i)$. Finally, let $\beta_{v,tot}=\beta_v+\sum \beta_{v,i}\in
\pi_2(X,F_{b_v})$, where $\beta_v=[u_v(D_v)]$ is the class of the stable
disc $u_v:D_v\to X$.

We use this data to define two loops in $F_{b_{v}}$, both based at 
the image of the output marked point of $D_v$, $x=u_{v}(z_{v,0})$.
On one hand, let $\sigma_1=s_{[\partial \beta_{v,tot}],x}$ be the standard loop at 
$x$ in the homology class $\partial \beta_{v,tot}$. On the other hand, let
$\sigma_0$ be the loop obtained by inserting the standard loop 
$s_{[\partial \beta_{v,i}],x_i}$ at each input marked point $z_{v,i}$ into the boundary loop $u_{v|\partial D_v}$ of
the disc $u_v$.
This loop does not have a canonical parametrization by $S^1$, but we
can choose one in a consistent manner, using the fact that the domain
$D_v$ is stable (possibly after adding interior marked points corresponding
to intersections with stabilizing divisors).

Denote by $\sigma:[0,1]\times S^1\to F_{b_v}$ a homotopy between $\sigma_0$ and 
$\sigma_1$ produced by some consistent method of interpolation between
based loops in the fibers of $\pi$; for example, after identifying all the
fibers with flat tori we can just use straight line
interpolation. We set $S^{spl}(t)=[0,1]\times S^1$, $\sigma_\alpha=\sigma$, and require the
end point of $e^{spl}_\alpha$ to map to $\sigma(\theta_t)$.

This choice is motivated by the observation that the boundary of the
homotopy $\sigma$ precisely
accounts for the various ways in which a splicing with target
$v$ can deform: the target can move into the output edge, whence
the required incidence condition becomes a standard loop 
in the class $\partial \beta_{v,tot}$ (cf.\ Case 1 above), or it can
move into one of the input edges $e_i$, and the incidence condition
becomes
a standard loop in the class $\partial \beta_{v,i}$; or the splicing can
disappear altogether by deforming to an honest gradient flow line attached
to $u_v(\partial D_v)$.

\smallskip

Things become more complicated when two or more splicings share the
same target. We describe the splicing data and incidence condition in
the next simplest case, to illustrate the general construction, which
will appear elsewhere. \smallskip

{\em Case 3.} An interior point $t$ of an edge $e$ of $T^+$ is the
target of exactly two splicings $e^{spl}_{\alpha_1}$ and $e^{spl}_{\alpha_2}$.
Denote by $\beta=\beta_{t\in e}$, $\beta_1=\beta_{t\in e^{spl}_{\alpha_1}}$
and $\beta_2=\beta_{t\in e^{spl}_{\alpha_2}}$ the weights of the different
edges which attach together at $t$, and let $x=u_e(t)$.

As in Case 2 above, the space $S^{spl}(t)$ and the maps
$\sigma_{\alpha_1},\sigma_{\alpha_2}:S^{spl}(t)\to X^{00}$ should describe a
homotopy between the incidence conditions imposed on the ends of the
splicing edges $e^{spl}_{\alpha_1}$ and $e^{spl}_{\alpha_2}$ after 
small deformations which make their targets
$t_1$ and $t_2$ distinct.
There are four manners in which such a configuration can deform to one
where the two splicings have distinct targets $t_1$ and $t_2$:
\begin{itemize}
\item Type I: $t_1$ lies before $t_2$ along $e$ (farther from the output), 
\item Type I': $t_1$ lies after $t_2$ along $e$ (closer to the output),
\item Type II: $t_1$ can move to the edge $e^{spl}_{\alpha_2}$, 
\item Type II': $t_2$ can move to $e^{spl}_{\alpha_1}$.
\end{itemize}
When $t_1$ lies before $t_2$ along $e$ (Type I), the
incidence condition at $t_1$ is given by a standard loop in the class $\partial
\beta$, while at $t_2$ it is a standard loop in the class
$\partial(\beta+\beta_1)$ (both based at points close to $x$). 
When $t_1$ lies on $e^{spl}_{\alpha_2}$ instead (Type II), the incidence condition at $t_2$ is a standard loop
in the class $\partial \beta$ (based near $x$), while at $t_1$ it is a standard loop in
the class $\partial \beta_2$ (now based near the end point of the second
splicing, rather than $x$). Similarly for Types I' and II', exchanging the indices 1 and 2.

We set $S^{spl}(t)$ to be the disjoint union of two copies of $[0,1]\times S^1\times
S^1$. On the first one, we define
$\sigma_{\alpha_1}(\tau,\theta_1,\theta_2)=s_{[\partial
\beta],x}(\theta_1)$, i.e.\ the incidence condition for $e^{spl}_{\alpha_1}$
is independent of $\tau\in [0,1]$ and lies along the standard loop at $x$ in
the class $\partial \beta$. Meanwhile, we pick for each value of $\theta_1$
a homotopy (chosen by some consistent process, e.g.\ straight line
interpolation after identifying $F_{\pi(x)}$ with a flat torus)
$$\varsigma_{[\partial\beta],[\partial\beta_1],x,\theta_1}:[0,1]\times S^1\to F_{\pi(x)}$$
between the standard loop $s_{[\partial\beta]+[\partial\beta_1],x}$ (for
$\tau=0$) and the
loop obtained by inserting $s_{[\partial \beta_1],y}$ 
into $s_{[\partial\beta],x}$ at the point
$y=s_{[\partial\beta],x}(\theta_1)$ (for $\tau=1$); and we set
$$\sigma_{\alpha_2}(\tau,\theta_1,\theta_2)=\varsigma_{[\partial\beta],[\partial\beta_1],x,\theta_1}(\tau,\theta_2).$$
This parametrizes a homotopy between the incidence conditions associated to
Type I deformations (for $\tau=0$) and the union of the incidence conditions for Type II'
deformations together with the ``symmetric'' incidence condition where both
splicings map to the standard loop $s_{[\partial\beta],x}$.
On the second copy of $[0,1]\times S^1\times S^1$, we set instead
$$\sigma_{\alpha_1}(\tau,\theta_1,\theta_2)=\varsigma_{[\partial\beta],[\partial\beta_2],x,\theta_2}(\tau,\theta_1)
\quad\text{and}\quad \sigma_{\alpha_2}(\tau,\theta_1,\theta_2)=s_{[\partial \beta],x}(\theta_2).$$
This yields a homotopy between the incidence conditions for the remaining
types of deformation (Types I' and II) and the symmetric incidence condition
$s_{[\partial\beta],x}\times s_{[\partial\beta],x}$, so that the two
components of $S^{spl}(t)$ taken together provide the desired homotopy
between incidence conditions.

\medskip

We expect that a similar construction of homotopies between different
incidence conditions can be used to deal with the remaining cases (when two
splicings have a vertex as common target, or when more than two splicings
have the same target). A detailed treatment will appear elsewhere.

\subsubsection{The master equation}

Returning to the notation of \S\S \ref{sss:morse}--\ref{sss:adapted}, we
define 
$$\m_0^{spl}\in C^*(X^{00},\pi^*\O_{an})=C^*(B^0;H^*(F_b)\,\hat\otimes\,\O_{an})$$
to be the same weighted sum as in \eqref{eq:Morse_md} (with \hbox{$d=0$}), except we
use moduli spaces of spliced treed holomorphic discs instead of their
ordinary counterparts.
(As usual, $\m_0^{spl}$ is a weighted count of {\em rigid}\/ spliced treed
holomorphic discs, i.e.\ those which arise in zero-dimensional moduli spaces,
while the master equation comes from considering the boundaries of
one-dimensional moduli spaces.)
The master equation is now expected to arise from the behavior of spliced
treed holomorphic discs at the boundary of the moduli space.

Boundary configurations where the output edge of the tree breaks through a 
critical point of the Morse function $f$ contribute $\delta \m_0^{spl}$, where
$\delta$ is the Morse differential. 

Otherwise, once the length
of an internal edge $e$ of a treed holomorphic disc becomes infinite, it
turns into a splicing edge, whose target can move up along the
remaining part of the tree, all the way to the 
output edge $e_{out}$. The boundary of the moduli space of spliced discs
is reached once the target of the splicing has moved ``to
infinity'' along $e_{out}$, at which point the gradient flow line
corresponding to $e_{out}$ must also break through a critical point of $f$
below the target of the splicing.
Thus, in the limit we have a pair of rigid treed 
holomorphic discs $u_\pm:\mathbb{T}_\pm\to X$ with outputs
$p_{\pm}\in\mathrm{crit}(f)$, together with gradient flow lines
of $f$ from $p_{\pm}$ to a pair of points $x_{\pm}$ with the property that
$x_-$ lies on the standard loop $s_{[\gamma_+],x_+}$ through $x_+$
representing the homology class $[\gamma_+]=[\partial \beta_+]$ of the
boundary of the treed disc $u_+$ (and then upward from $x_+$
to a critical point of $f$ which is the overall output of the limit configuration). 
The conditions involving $p_\pm$ and $x_\pm$ can be rewritten as:
$$x_+\in [\overline{W}^+(p_+)]\cap p_*s_{[\gamma_+]}^{-1}([\overline{W}^+(p_-)]),$$
which by \eqref{eq:stdloops} is the same as
$[\overline{W}^+(p_+)]\cap \bigl[\overline{W}^+(\iota_{[\gamma_+]}(p_-))\bigr]$.
In other terms, the gradient flow line from $p_-$ to $x_-$ can be replaced
by a gradient flow line from one of the critical points appearing in the
linear combination $\iota_{[\gamma_+]}(p_-)$ to $x_+$, and after this
modification the top-most portion of the limiting configuration amounts to
a gradient flow tree computing the cup-product $p_+\wedge
\iota_{[\gamma_+]}(p_-)$.

Summing over all such configurations, with appropriate weights, we therefore
arrive at $$\m_0^{spl}*\m_0^{spl}=\frac12 \{\m_0^{spl},\m_0^{spl}\},$$
where $*$ denotes the operation of degree $-1$ on
$C^*(X^{00},\pi^*\O_{an})=C^*(B^0;H^*(F_b)\,\hat\otimes\,\O_{an})$ defined
by 
\begin{equation}\label{eq:morsestar}
(z^{\gamma_+} p_+)*(z^{\gamma_-} p_-)=z^{\gamma_++\gamma_-}\, p_+\wedge
\iota_{\gamma_+}(p_-),
\end{equation}
whose skew-symmetrization is the bracket
\eqref{eq:HF-bracket}.

Unlike the case of ordinary treed holomorphic discs, configurations where
the targets of several spliced edges simultaneously escape to infinity can also contribute 
to the codimension 1 boundary of the moduli space of spliced treed discs.

When the targets of two different splicings $e_1$ and $e_2$ both escape
towards the output of the spliced treed disc in such a way that the distance
between the two targets also goes to infinity, we arrive at
codimension 2 strata consisting of three treed discs $u_1,u_2,u_+$ with
outputs $p_1,p_2,p_+$, together with gradient flow lines of $f$ which attach 
onto each other via standard loops. 
Using \eqref{eq:stdloops}, these can be recast by the same trick as above
as a broken gradient flow tree computing one of $p_+\wedge \iota_{[\gamma_+]}(p_1
\wedge \iota_{[\gamma_1]}(p_2))$, $p_+\wedge \iota_{[\gamma_+]}(p_2
\wedge \iota_{[\gamma_2]}(p_1))$, $(p_+\wedge \iota_{[\gamma_+]}(p_1))
\wedge \iota_{[\gamma_++\gamma_1]}(p_2)$, or $(p_+\wedge \iota_{[\gamma_+]}(p_2))
\wedge \iota_{[\gamma_++\gamma_2]}(p_1)$. (Here
$[\gamma_1],[\gamma_2],[\gamma_+]$ are the homology classes associated to the
boundary loops of the treed discs $u_1,u_2,u_+$.)

However, when the targets of $e_1$ and $e_2$ remain a finite distance
apart, the configuration of gradient flow lines gets recast as a gradient
flow tree of the sort used to define higher $A_\infty$-operations in Morse
theory, and/or 
the incidence conditions we impose on the end points of the gradient flow
lines involve homotopies between the various products of standard loops
appearing in the above expressions (see the discussion of Case 3 in \S
\ref{sss:spliced} above). Due to the extra degree of freedom afforded
by the various homotopies, these strata have codimension 1 rather than~2. 

Algebraically, these homotopies define an operation of
degree $-3$ on $C^*(X^{00},\pi^*\O_{an})^{\otimes 3}$, whose
skew-symmetrization $\ell_3$ is the next term in a shifted
$L_\infty$-structure whose first two operations are the Morse differential
$\delta$ and (up to sign) the bracket $\{\cdot,\cdot\}$. Configurations in which
two splicing targets escape to infinity then contribute an additional term
$\frac16 \ell_3(\m_0^{spl},\m_0^{spl},\m_0^{spl})$ to the master equation; and so on with
higher homotopies when
more than two edge lengths simultaneously become infinite.
To summarize:

\begin{conj}
(1) The Morse complex $C^*(X^{00},\pi^*\O_{an})$ carries 
operations $$\ell_m:C^*(X^{00},\pi^*\O_{an})^{\otimes m}\to
C^*(X^{00},\pi^*\O_{an})[3-2m],\quad m\geq 1,$$ 
defined in terms of counts of spliced configurations consisting of $m-1$
splicing edges and one infinite gradient flow line, without any disc
components; in particular
$\ell_1$ is the Morse
differential and $\ell_2$ is (up to sign) the bracket $\{\cdot,\cdot\}$.
These operations define a shifted $L_\infty$-structure on $C^*(X^{00},\pi^*\O_{an})$.

(2) Weighted counts of rigid spliced treed holomorphic discs define
an element $\m_0^{spl}\in C^*(X^{00},\pi^*\O_{an})$ which satisfies the $L_\infty$
master equation
\begin{equation}
\label{eq:mastereq_Linfty}
\sum_{m\geq 1} \frac{1}{m!} \ell_m((\m_0^{spl})^{\otimes
m})=\delta\m_0^{spl}\pm
\frac12 \{\m_0^{spl},\m_0^{spl}\}+\frac16
\ell_3(\m_0^{spl},\m_0^{spl},\m_0^{spl})+\dots=0.
\end{equation}
\end{conj}\smallskip

\begin{remark}
Since the product on the Morse complex is only homotopy associative, it shouldn't
be surprising that in general the bracket $\{\cdot,\cdot\}$ should only be the leading term of an
$L_\infty$ structure, and hence that the master equation should also involve higher
order terms. As noted in \S \ref{ss:mastereq} above, Irie's result on the chain-level
master equation with loop space coefficients \cite{Irie} also involves 
higher terms, due to the chain-level string bracket being part of
an $L_\infty$ structure on the space of chains on the free loop space.

However, we expect that, for suitable choices of the 
adapted Morse function $f$ and of the system of standard loops
$(s_{[\gamma],x})$, one
can arrange for the higher terms of the $L_\infty$-structure to vanish,
reducing \eqref{eq:mastereq_Linfty} to the ordinary master equation
\eqref{eq:mastereq_m0}. The reason for this expectation is that, choosing
$f$ as in Remark \ref{rmk:morse_to_cech}, the Morse complex we consider can be
recast as the \v{C}ech complex of a suitable cover of $B^0$ with coefficients in
$H^*(F_b)\,\hat\otimes \O_{an}$, which is a dg-algebra; and, as we shall
see below, the bracket $\{\cdot,\cdot\}$ corresponds under
mirror symmetry to the Schouten-Nijenhuis bracket for \v{C}ech cochains on
the uncorrected mirror $X^{\vee 0}$ with coefficients in polyvector fields,
which satisfies the Jacobi identity at chain level.
\end{remark}

\begin{remark} 
Spliced treed discs can be used to produce not only the element $\m_0^{spl}$ 
but also a wealth of algebraic operations on $C^*(X^{00}, \pi^*\O_{an})$.
Work in progress of the author with Keeley Hoek suggests that, by counting
spliced treed holomorphic discs whose inputs all lie on the tree that contains the
output (i.e., the treed discs that feed into splicings do not carry any
inputs), one can define operations $\tilde{\m}_k^{spl}$ for $k\geq 1$ which make
$C^*(X^{00},\pi^*\O_{an})$ into an {\em uncurved} $A_\infty$-algebra. 
The details will appear elsewhere. Beyond this, one might hope that more
general moduli spaces of spliced treed discs also endow
$C^*(X^{00},\pi^*\O_{an})$ with the structure of a framed $E_2$-algebra (or
homotopy BV algebra); compare with \cite{AbGrVa}.
\end{remark}

\subsection{From $\m_0$ to the geometry of the corrected mirror}\label{ss:instantoncorr}

We finally turn our attention to the passage from family Floer theory to the geometry of
instanton corrections on the mirror. As explained at the beginning of
\S \ref{ss:famfloer}, the uncorrected
mirror $X^{\vee 0}$ comes equipped with a rigid analytic torus fibration
$\pi^\vee:X^{\vee 0}\to B^0$, locally modelled on 
the valuation map $H^1(F_b,\K^*)\to H^1(F_b,\R)$, or more explicitly
after choosing a basis $(\gamma_1,\dots,\gamma_n)$ of $H_1(F_b,\Z)$, the valuation map $(\K^*)^n\to \R^n$. 

Under the assumption of weak family unobstructedness (Definition
\ref{def:wfunobs}), family Floer theory as described in the preceding sections determines an element 
$\m_0\in C^*(B^0,H^*(F_b)\,\hat\otimes\,\O_{an})$,
where $\O_{an}=\pi^\vee_*(\O_{X^{\vee 0}})$ is a completion of the
ring $\K[H_1(F_b)]$ of Laurent polynomials in the local coordinates of
$X^{\vee 0}$. Moreover, $\m_0$ can be expressed as a sum of elements $$\alpha^{(i)}\in
C^i(B^{0},H^i(F_b)\,\hat\otimes\,\O_{an})$$ which encode counts of
holomorphic discs of Maslov index $2-2i$ in $X$ bounded by $i$-dimensional families
of fibers of the Lagrangian torus fibration $\pi:X^{00}\to B^0$.

The natural isomorphism $H^1(F_b,\R)\simeq T_b B$ allows us to map elements
of $H^1(F_b)$ (resp.\ $H^*(F_b)=\Lambda^* H^1(F_b)$) to vector fields (resp.\ poly vector
fields) on $X^{\vee 0}$. Namely, denoting by $(z_1,\dots,z_n)$ the
local coordinates on $X^{\vee 0}$ induced by a choice of basis
$(\gamma_1,\dots,\gamma_n)$ of $H_1(F_b)$, and by 
$(\gamma_1^*,\dots,\gamma_n^*)$ the dual basis of $H^1(F_b)$,
we map each basis element $\gamma_j^*$ to the vector field 
$\partial_{\log z_j}=z_j\partial_{z_j}$, and extend this map to $H^*(F_b)\,\hat\otimes\,
\O_{an}$ by setting
\begin{equation}\label{eq:floertopolyvf0}
z^\gamma\,\gamma_{j_1}^*\wedge \dots\wedge \gamma_{j_k}^* \mapsto
z^\gamma\,\partial_{\log z_{j_1}}\wedge\dots\wedge \partial_{\log z_{j_k}}.
\end{equation}

Combining this with pullback of cochains under the projection $\pi^\vee:X^{\vee
0}\to B^0$, we obtain a (bigraded) map
\begin{equation}
\label{eq:floertopolyvf}
C^*(B^0,H^*(F_b)\,\hat\otimes\, \O_{an})\to C^*(X^{\vee 0},\Lambda^* T_{X^{\vee 0}}).
\end{equation}

If we use the Morse-theoretic model of \S \ref{sss:morse} for an adapted
Morse function in the sense of \S \ref{sss:adapted}, then by Remark
\ref{rmk:morse_to_cech} the Morse cochains on the left-hand side
of \eqref{eq:floertopolyvf} can be recast as \v{C}ech cochains for a
certain polyhedral cover of $B^0$ (by the stars $\mathcal{U}_v$ of the vertices of the
simplicial decomposition $\mathcal{P}$). The right-hand side of
\eqref{eq:floertopolyvf} should then be interpreted as \v{C}ech cochains for
a cover of $X^{\vee 0}$ by affinoid domains approximating the preimages
$(\pi^\vee)^{-1}(\mathcal{U}_v)$, $v\in\mathrm{vert}(\mathcal{P})$.

On the other hand, if we assume the existence of a model for family Floer theory
in which $\m_0$ is expressed as a differential form on $B^0$ with values in
$H^*(F_b)\,\hat\otimes\, \O_{an}$, then the right-hand side of \eqref{eq:floertopolyvf}
should be interpreted in terms of {\em tropical
differential forms} (also known as {\em superforms}) on $X^{\vee 0}$ (see e.g.\ \cite{CLD,Jell}). Choosing a basis of $H_1(F_b)$ as
above, and denoting by $(x_1,\dots,x_n)$ and $(z_1,\dots,z_n)$ the
corresponding local coordinates on $B^0$ and on $X^{\vee 0}$ (with
$\val(z_j)=x_j$), we define
$$(\pi^\vee)^*(dx_j)=d''\log|z_j|,$$
a superform of type $(0,1)$; and similarly for exterior products. By
definition, the pullback of differential forms intertwines
the de Rham differential $d$ on
$\Omega^*(B^0)$ and the tropical Dolbeault differential $d''$ on
$\Omega^{0,*}(X^{\vee 0})$. Consequently, the map \eqref{eq:floertopolyvf}
also intertwines the de Rham differential $d$ on $\Omega^*(B^0,H^*(F_b)\,\hat\otimes\,
\O_{an})$ and the tropical Dolbeault differential $d''$ on $\Omega^{0,*}(X^{\vee
0},\Lambda^*T_{X^{\vee 0}})$.

The key to the proof of Proposition \ref{prop:Wi} is the following lemma:

\begin{lemma}\label{l:brackets}
The map \eqref{eq:floertopolyvf} intertwines the bracket $\{\cdot,\cdot\}$
on $C^*(B^0;H^*(F_b)\,\hat\otimes\,\O_{an})$ defined by \eqref{eq:HF-bracket}
and the negative of the Schouten-Nijenhuis bracket
$-[\cdot,\cdot]$ on $C^*(X^{\vee 0},\Lambda^*T_{X^{\vee 0}})$.
\end{lemma}

\proof
Since the map \eqref{eq:floertopolyvf} is compatible with the cup-product of
cochains, it suffices to compare the brackets on
$H^*(F_b)\,\hat\otimes\,\O_{an}$ and on $\Lambda^*T_{X^{\vee 0}}$.

Recall that the Schouten-Nijenhuis bracket is a bracket of degree $-1$ on
polyvector fields, characterized by the following properties \cite{Marle}:
given a smooth function $f$, vector fields $X,Y$, and polyvector fields $P,Q,R$
of degrees $p,q,r$,
\begin{enumerate}
\item $[X,f]=L_X f=\iota_{df}(X)$;
\item $[X,Y]=L_X Y$;
\item $[P,Q]=-(-1)^{(p-1)(q-1)}[Q,P]$;
\item $[P,Q\wedge R]=[P,Q]\wedge R+(-1)^{(p-1)q}Q\wedge [P,R]$;
\item $[P\wedge R,Q]=P\wedge [R,Q]+(-1)^{(q-1)r}[P,Q]\wedge R$.
\end{enumerate}
Assume now that the polyvector fields $P$ and $Q$ have constant components
in some local coordinate system, and let $f,g$ be two smooth functions.
Then properties (2), (4) and (5) imply that $[P,Q]=0$, whereas (1) and (5)
imply that $[P,g]=(-1)^{p-1}\iota_{dg}(P)$. Hence, by (4) we have
$[P,gQ]=(-1)^{p-1}\iota_{dg}(P)\wedge Q$. Meanwhile, (3) and (4)
imply that $[f,gQ]=g[f,Q]=-g\iota_{df}(Q)$. Finally, using (5) once more
we arrive at
\begin{align}\label{eq:schouten}
\nonumber [fP,gQ]&=(-1)^{p-1}f\iota_{dg}(P)\wedge Q-(-1)^{(q-1)p}g\iota_{df}(Q)\wedge P
\\&=(-1)^{p-1}f\,\iota_{dg}(P)\wedge Q-g\,P\wedge \iota_{df}(Q).\end{align}

Recall that, for $\gamma\in H_1(F_b)$, $\alpha,\alpha'\in H^*(F_b)$, 
we define 
\begin{equation}\label{eq:HFbracketagain}
\{z^\gamma\, \alpha, z^{\gamma'} \alpha'\}=z^{\gamma+\gamma'}\,
\bigl(\alpha\wedge (\iota_\gamma \alpha')+(-1)^{|\alpha|}
(\iota_{\gamma'}\alpha)\wedge \alpha'\bigr).
\end{equation}
Denote by $V,V'\in \Lambda^*T_{X^{\vee 0}}$ the images of $\alpha,\alpha'$
under \eqref{eq:floertopolyvf0}, which are polyvector fields with constant
coefficients in terms of the basis formed by exterior products of
$\partial_{\log z_i}$. By \eqref{eq:schouten},
\begin{equation}\label{eq:schoutenpartialcalc}
-[z^\gamma\, V,z^{\gamma'}V']=
z^{\gamma'}\,V\wedge \iota_{d(z^{\gamma})}(V')+
(-1)^{|\alpha|}z^{\gamma}\,\iota_{d(z^{\gamma'})}(V)\wedge V'.
\end{equation}
Thus, in order to complete our comparison of the two brackets it suffices to show that \eqref{eq:floertopolyvf0} maps
$z^{\gamma}\,\iota_{\gamma}\alpha'$ to $\iota_{d(z^{\gamma})}(V')$,
and similarly for the other interior product appearing in \eqref{eq:HFbracketagain}.

Distributing the interior products into the expressions of $\alpha'$ and
$V'$ in the chosen bases, it is in fact
sufficient to compare $z^{\gamma}\,\iota_{\gamma}\alpha'$ to $\iota_{d(z^{\gamma})}(V')$
in the specific case where $\alpha'=\gamma_{j}^*$ is an element 
of the chosen basis of $H^1(F_b)$, and $V'=\partial_{\log z_j}$. Expressing
$\gamma$ in the chosen basis of $H_1(F_b)$ as
$\gamma=a_1\gamma_1+\dots+a_n\gamma_n$, we find that
$$z^{\gamma}\,\iota_{\gamma}(\gamma_j^*)=a_jz^{\gamma}=\partial_{\log
z_j}(z^{\gamma})=\iota_{d(z^{\gamma})}(\partial_{\log z_j}),$$
which completes the proof.
\endproof

\proof[Proof of Proposition \ref{prop:Wi}]
Recall that $\mathbb{W}$ (resp.\ $W^{(i)}$) is by definition the image of $\m_0$ (resp.\ its
components $\alpha^{(i)}$) under the map \eqref{eq:floertopolyvf}.
The compatibility of \eqref{eq:floertopolyvf} with the differentials
implies that it maps $\delta \m_0$ to $\delta \mathbb{W}$; while
Lemma \ref{l:brackets} implies that it maps $\{\m_0,\m_0\}$ to
$-[\mathbb{W},\mathbb{W}]$. The master equation \eqref{eq:mastereq_m0}
for $\m_0$ thus maps to the equation \eqref{eq:mastereq} for $\mathbb{W}$.
\endproof

We also give a derivation of equations \eqref{eq:mastereq0}--\eqref{eq:mastereq2}
for completeness. Recall that the Schouten-Nijenhuis bracket satisfies the
Jacobi identity
$$[[P,Q],R]=[P,[Q,R]]-(-1)^{(p-1)(q-1)}[Q,[P,R]]$$
(where $p=\deg P$, $q=\deg Q$). Since $\mathbb{W}$ is even, this yields
$$[[\mathbb{W},\mathbb{W}],\cdot]=2[\mathbb{W},[\mathbb{W},\cdot]].$$
Moreover, the differential $\delta$ on cochains does not interact with the
Schouten-Nijenhuis bracket (this is manifest in the case of \v{C}ech cochains,
and for tropical differential forms it follows from the fact that $\delta$
is the Dolbeault differential $d''$ on $(0,*)$-forms while the 
Schouten-Nijenhuis bracket involves differentiation along analytic vector
fields). Hence,
$$\delta([\mathbb{W},\cdot])=[\delta\mathbb{W},\cdot]-[\mathbb{W},\delta(\cdot)].$$
Thus, assuming \eqref{eq:mastereq}, we have
\begin{equation}\label{eq:mastereq_diffl}
(\delta+[\mathbb{W},\cdot])^2=\delta[\mathbb{W},\cdot]+[\mathbb{W},\delta(\cdot)]+
[\mathbb{W},[\mathbb{W},\cdot]]=
[\delta\mathbb{W},\cdot]+\frac12 [[\mathbb{W},\mathbb{W}],\cdot]=0.
\end{equation}
Writing $\mathbb{W}=W^{(0)}+W^{(1)}+\dots$ with $W^{(i)}\in C^i(X^{\vee
0},\Lambda^i T_{X^{\vee 0}})$, the component of
\eqref{eq:mastereq} in bidegree $(1,0)$ (i.e., in $C^1(X^{\vee 0},\O_{X^{\vee 0}})$) is
$$\delta W^{(0)}+[W^{(1)},W^{(0)}]=0,$$
which gives \eqref{eq:mastereq0}. The component of \eqref{eq:mastereq_diffl}
in bidegree $(2,1)$ is
$$(\delta+[W^{(1)},\cdot])^2+[W^{(2)},[W^{(0)},\cdot]]+[W^{(0)},[W^{(2)},\cdot]]=0,$$
which can be rewritten using the Jacobi identity as
$$(\delta+[W^{(1)},\cdot])^2+[[W^{(2)},W^{(0)}],\cdot]=0.$$
Since $[W^{(2)},W^{(0)}]=-\iota_{dW^{(0)}}(W^{(2)})$, this yields
\eqref{eq:mastereq1}.
Next, the component of \eqref{eq:mastereq} in bidegree $(3,2)$ is
$$\delta W^{(2)}+[W^{(1)},W^{(2)}]+[W^{(0)},W^{(3)}]=0,$$ which yields
\eqref{eq:mastereq2}; and so on.

\begin{remark}\label{rmk:gluingphilosophy}
The components of the \v{C}ech cochain groups $C^*(B^0;H^*(F_b)\hat\otimes\O_{an})$
obtained from a simplicial decomposition of $B^0$ are certain
completions of $H^*(F_b)\otimes \K[H_1(F_b)]$, which can be viewed as
the symplectic cohomologies (in the classical limit, before any 
instanton corrections) of the local pieces of the corresponding
decomposition of $X^{00}$. (See also Groman and Varolgunes' work
introducing the relative symplectic cohomology sheaf of an SYZ fibration
\cite{GromanVarolgunes}.) In this sense, the map \eqref{eq:floertopolyvf}
can be viewed as a local version of the expected isomorphism between
the Hochschild cohomologies of the Fukaya category of $X^{00}$ (without
instanton corrections) and the derived category of the uncorrected
mirror $X^{\vee 0}$, with the latter recast in terms of polyvector fields
via the Hochschild-Kostant-Rosenberg isomorphism. In this language,
Lemma \ref{l:brackets} expresses the fact that homological
mirror symmetry intertwines the BV-structures on these local Hochschild
cohomologies.

The claim that the element $\mathbb{W}\in C^*(X^{\vee 0},\Lambda^*T_{X^{\vee
0}})$ obtained by applying \eqref{eq:floertopolyvf} to $\m_0$ represents
the instanton corrections to be applied to the geometry of the uncorrected
mirror $X^{\vee 0}$ also follows naturally from this perspective.
Namely, we expect that the Fukaya category of $X$ can be recovered from
the local relative wrapped Fukaya categories of the pieces of $X^{00}$ 
by a pullback diagram. Prior to any instanton corrections,
the ``classical limit'' of this pullback diagram (using the ``classical''
local wrapped categories rather than their relative deformations inside $X$)
matches the description of the derived category of the uncorrected mirror 
$X^{\vee 0}$ in terms of local pieces.
The components of the Floer-theoretic obstruction
$\m_0\in C^*(B^0;H^*(F_b)\hat\otimes\O_{an})$ describe the deformations
of the classical local wrapped categories and of their gluing data 
needed to account for holomorphic discs in $X$ and arrive at the
relative wrapped Fukaya categories and the diagram via which they recover
the Fukaya category of $X$.  Thus, the corresponding element
$\mathbb{W}\in C^*(X^{\vee 0},\Lambda^*T_{X^{\vee 0}})$ should similarly be used
to deform the local pieces of the uncorrected mirror $X^{\vee 0}$ and 
their gluing data in order to arrive at the correct mirror.

(In the Calabi-Yau setting, the above ideas are also
likely related to Chan, Leung and Ma's work 
on the construction of mirror spaces by using Maurer-Cartan elements to
deform gluings of BV algebras \cite{ChanLeungMa}.)
\end{remark}

\end{document}